\documentclass[a4paper, 11pt]{article}

\usepackage[latin1]{inputenc}
\usepackage{makeidx}
\usepackage[francais,english]{babel}

\usepackage{amsmath,amsthm,amsfonts,amssymb}
\usepackage{newlfont}
\usepackage{showidx}
\usepackage[dvips]{graphicx}
\newtheorem{theo}{Théorème}
\newtheorem{prop}{Proposition}
\newtheorem{cor}{Corollaire}

\newtheorem{defi}{Définition}
\newtheorem{lem}{Lemme}

\def\g{\mathfrak{g}}
\def\G{\Gamma}
\def\ov{\overline}
\def\tr{\mathrm{tr}}
\def\ad{\mathrm{ad}}
\def\p{\mathfrak{p}}\def\b{\mathfrak{b}}
\def\h{\mathfrak{h}}\def\r{\mathfrak{r}}
\def\q{\mathfrak{q}}\def\n{\mathfrak{n}}
\def\k{\mathfrak{k}}
\def\s{\mathfrak{s}}
\def\Exp{\mathrm{Exp}}
\newcommand{\R}{\mathbb{R}}
\newcommand{\C}{\mathbb{C}}

\newcommand{\semi}{\,\triangleright\!\!\! <}
\newcommand{\fin}{\hfill $\blacksquare$ \\}
\numberwithin{equation}{section}

\newcommand{\bbR}{{\mathbb{R}}}

\newcommand{\bx}{\mathbf{x}}

\newcommand{\hs}{{\sqsubset}}
\newcommand{\hsb}{{\sqsubseteq}}
\newcommand{\hpi}{{\frac\pi 2}}
\newcommand{\ndash}{\nobreakdash-\hspace{0pt}}

\newcommand{\dd}{{\mathrm{d}}}

\newcommand{\ee}{{\mathrm{e}}}

\begin{document}

\bibliographystyle{alpha}
\author{Alberto S.
Cattaneo\footnote{Institut für Mathematik, Universität Zürich,
Winterthurerstr. 190 CH-8057 ZÜRICH, alberto.cattaneo@math.unizh.ch}
\and Charles Torossian\footnote{New address : Institut Mathématiques
de Jussieu, Équipe de Théorie des Groupes, Université Paris 7, Case
7012, 2 place Jussieu,  75251 Paris Cedex 05 FRANCE, torossian@math.jussieu.fr}}

\title{Quantification pour les paires symétriques et diagrammes de Kontsevich}

\maketitle

{\noindent \bf{Résumé: }} Dans cet article nous appliquons les
méthodes de bi-quantification décrites dans \cite{CF1}
 au cas des espaces symétriques. Nous introduisons
une fonction $E(X,Y)$, définie pour toutes paires symétriques, en
termes de graphes de Kontsevich. Les propriétés de cette fonction
permettent de démontrer de manière unifiée des résultats importants
dans le cas des paires symétriques résolubles ou quadratiques. Nous
montrons que le star-produit décrit dans \cite{CF1} coïncide, pour
toute paire symétrique, avec celui de Rouvière. On généralise un
résultat de Lichnerowicz sur la commutativité d'algèbres
d'opérateurs différentiels invariants et on résout un problème de M.
Duflo sur l'écriture, en coordonnées exponentielles,  des opérateurs différentiels invariants sur
tout espace symétrique. On décrit
l'homomorphisme d'Harish-Chandra en termes de graphes de
Kontsevich. On développe une théorie nouvelle pour construire des
caractères des algèbres d'opérateurs différentiels invariants. On
applique ces méthodes dans
le cas des polarisations $\sigma$-stables.\\

\selectlanguage{english}
\begin{center}
{\Large  Quantization for symmetric pairs and Kontsevich's diagrams}
\end{center}

 {\noindent
\bf{Abstract:}} In this article we use the expansion for biquantization
described in \cite{CF1} for the case of symmetric spaces. We
introduce a function of two variables $E(X,Y)$ for any symmetric
pairs. This function has an expansion in terms of Kontsevich's
diagrams. We recover most of the known results though in a more
systematic way by using some elementary properties of this $E$
function. We prove that Cattaneo and Felder's star product coincides
with Rouvière's for any symmetric pairs. We generalize some  of
Lichnerowicz's results for the commutativity of the algebra of
invariant differential operators and solve a long standing problem
posed by M. Duflo for the expression of invariant differential
operators on any symmetric spaces in exponential coordinates. We
describe the Harish-Chandra homomorphism in the case of symmetric
spaces by using all these constructions. We develop a new method to
construct characters for algebras of invariant differential
operators. We apply these methods in the case of $\sigma$-stable
polarizations.
\\


\noindent {\bf Acknowledgments :}  ASC acknowledges partial support
of SNF, Grant No. 200020-107444/1. 53D55. and thanks IHES for kind
hospitality during the completion of part of this work. This work
has been partially supported by the European Union through the FP6
Marie Curie RTN ENIGMA (Contract number MRTN-CT-2004-5652). Research of C.T. was supported by CNRS.

ASC and CT thank G. Felder, D. Indelicato, B. Keller and  D. Manchon
for useful discussions. We thank Jim Stasheff and F. Rouvière for their comments on
the draft version of this paper. We are indebted to the referees of this paper for valuable remarks.


\section*{Introduction}

Let $(\g, \sigma)$ be a symmetric pair: viz., $\g$ is a
finite-dimensional Lie algebra over $\R$, while $\sigma$ is an
involution and a Lie algebra automorphism of $\g$. We denote by
$\g=\k\oplus \p$ the decomposition relative to $\sigma$,
with $\k$ and $\p$ the $+1$- and $-1$-eigenspaces, respectively (\textit{Cartan's decomposition}).\vspace{0.3cm}

The Poincaré--Birkhoff--Witt (PBW) theorem ensures the following
decomposition of the universal enveloping algebra $U(\g)$:
\[U(\g)=U(\g)\cdot \k \oplus \beta
(S(\p))
\]
with $\beta$ the symmetrization map from $S(\g)$ into
$U(\g)$.\footnote{We have $\beta(X_1\ldots
X_n)=\frac1{n!}\sum_{\sigma \in S_n} X_{\sigma(1)}\ldots
X_{\sigma(n)}$.} One can then identify, as vector spaces, the
symmetric algebra $S(\p)$ of $\p$ with $U(\g)/U(\g)\cdot \k$ via
$\beta$. Though $U(\g)/U(\g)\cdot \k$ is not an algebra in general,
its $\k$-invariant subspace $\big(U(\g)/U(\g)\cdot \k\big)^\k$ is an algebra.\vspace{0.3cm}

This is of fundamental importance as it is the algebra of invariant
differential operators on the symmetric space $G/K$ associated to
the symmetric pair $(\g, \sigma)$, and as such it occurs in harmonic
analysis on
symmetric spaces in a crucial way.\vspace{0.3cm}

This algebra is commutative \cite{lich, du79}, and the PBW theorem
ensures that
\[\Big(U(\g)/U(\g)\cdot \k\Big)^\k \quad \mathrm{and}\quad  S(\p)^\k\]
are isomorphic as vector spaces. Observe that $S(\p)^\k$ is just the
associated graded of
$\Big(U(\g)/U(\g)\cdot \k\Big)^\k$.\vspace{0.3cm}

It is conjectured \cite{To3} that these two algebras are isomorphic \textit{as
algebras}, what the second author has called the \textit{polynomial
conjecture}. This is a generalization for symmetric pairs of Duflo's
isomorphism \cite{du77} between the center $U(\g)^\g$ of the
universal enveloping algebra of a Lie algebra $\g$
and the invariant subalgebra $S(\g)^\g$ of its symmetric algebra.\vspace{0.3cm}

Recall that for every homogeneous space $G/H$ the algebra of
$G$-invariant differential operators may be identified, thanks to a
result of Koornwinder \cite{koorn}, with the algebra
$\big(U(\g)/U(\g)\cdot \h\big)^\h$ (where $\h$ denotes the Lie
algebra of $H$). The latter is not commutative in general. Its
associated graded is then a Poisson subalgebra of
$\big(S(\g)/S(\g)\cdot \h\big)^\h$ with a natural Poisson structure.\vspace{0.3cm}

If $\h$ admits a complement $\q$ which is invariant under the
adjoint action of $\h$, then an easy consequence of PBW is that
$\big(U(\g)/U(\g)\cdot \h\big)^\h$ and $\big(S(\g)/S(\g)\cdot
\h\big)^\h$ are still isomorphic as vector spaces. It is not known
whether this holds in general, for there is no natural map
from $\big(S(\g)/S(\g)\cdot \h\big)^\h$ to $U(\g)$ (or a quotient thereof).\vspace{0.3cm}

It has been conjectured anyway by M.~Duflo \cite{duflo-japon} that
the center of $\big(U(\g)/U(\g)\cdot \h\big)^\h$ and the Poisson
center of $\big(S(\g)/S(\g)\cdot\h\big)^\h$ are always
isomorphic as algebras.\vspace{0.3cm}

Little is known in such generality. In case  $\g$ is a nilpotent Lie
algebra, appreciable advances have been achieved  in the last few
years by Corwin--Greenleaf \cite{CG},
Fujiwara--Lion--Magneron--Mehdi \cite{FLMM}, Baklouti-Fujiwara
\cite{BF} and Baklouti-Ludwig \cite{BL}. In case where $G$ and $H$
are reductive groups, F.~Knop \cite{knop} gives a satisfying and
remarkable answer to the conjecture. In case where $H$ is compact
and $G=H\semi N$ is the semidirect product of $H$ with a Heisenberg
group $N$, Rybnikov
 \cite{ryb} makes use of F.~Knop's result to prove Duflo's conjecture.\vspace{0.3cm}

In this paper we propose a novel approach to these questions based
on Kontsevich's construction \cite{Kont} and its extension to the
case of coisotropic submanifolds by Cattaneo and Felder \cite{CF1},
\cite{CF2}. We only treat the problem of symmetric pairs here, but
we think that our methods have a wider scope, namely in the
nilpotent homogeneous case \footnote{One of our students is working on this case.}.\vspace{0.3cm}

One may regard the present work as a link between the methods of
deformation quantization
and the orbit method in Lie theory.\vspace{0.3cm}

\paragraph{Plan of the paper:}

In Section~\ref{sectionrappels} we recall Kontsevich's construction
for the deformation quantization of Poisson manifolds and its
extension by Cattaneo and Felder to the case of coisotropic
submanifolds. We discuss in details the compatibility in cohomology.
This Section should be useful for the Lie algebra experts who are
not familiar
with the deformation quantization constructions.\vspace{0.3cm}

Next we study the dependency of this construction on the choice of a
complement in the linear case. We show that the reduction spaces are
isomorphic and describe the isomorphism
(Proposition~\ref{propdependance} and Theorem~\ref{theodependance})
which is an element of the gauge group
obtained by solving a differential equation.\vspace{0.3cm}

In Section~\ref{ExempleReduction} we describe the graphs appearing
in the linear case and present three fundamental examples of
reduction spaces occurring in Lie theory: symmetric pairs
(Proposition~\ref{propreductionsym}), Iwasawa's decompositions
(Proposition~\ref{propreductioniwasawa}), and polarizations
(Proposition~\ref{propreductionpolarisation}). These examples are
new and show that this novel quantization methods
are well-suited for Lie theory.\vspace{0.3cm}

Recall that for symmetric pairs F.~Rouvière \cite{Rou90, Rou91,
Rou94}, following Kashi\-wara and Vergne \cite{KV}, introduced a
mysterious function $e(X,Y)$ defined for $X, Y \in \p$. Up to
conjugation this function computes the star product
\begin{equation}
P\underset{Rou}\sharp Q= \beta^{-1}\Big(\beta(P)\cdot
 \beta(Q) \quad \mathrm{modulo} \quad U(\g)\cdot \k \Big)\end{equation}
for $P, Q$ in $S(\p)^\k$.\vspace{0.3cm}

In Section~\ref{sectionfonctionE} we define a function $E(X,Y)$ for
$X, Y \in \p$ in terms of graphs. This function will behave as
Rouvière's function $e(X,Y)$. The comparison of these two functions
is a key point of this paper. The function $E(X, Y)$ expresses the
Cattaneo--Felder star product $\underset{CF} \star$ in the case of
symmetric pairs. By inspection of the graphs appearing in its
construction, we obtain a symmetry property (Lemma~\ref{symetrieE})
together with some additional properties in the solvable case
(Proposition~\ref{propE=1resoluble}), in the case of
Alekseev--Meinrenken symmetric pairs (Proposition~\ref{propE=1AM})
as well as in the case of very symmetric quadratic pairs
(Proposition~\ref{propE=1tressym}). In all these cases we show that
the function $E$ is identically equal to $1$. These elementary but
remarkable properties yield new and uniform proofs of results
obtained by Rouvière in the solvable case
(Theorem~\ref{theoResoluble},
Section~\ref{sectionpropritesupplementaires} and
Proposition~\ref{propordre4}) and generalize a theorem by Alekseev
and Meinrenken (Theorem~\ref{theoAM}) in the (anti-invariant)
quadratic case.\vspace{0.3cm}

In Section~\ref{sectionopdiff} we show that the Cattaneo--Felder and
the Rouvière star products coincide (Theorem~\ref{calculCF}). This
is a new result. {}From this we deduce (Theorem~\ref{theoz}) the
commutativity for all $z\in\mathbb{R}$ of the algebras
$\Big(U(\g)/U(\g)\cdot \k^{z\tr_k}\Big)^\k$ of invariant
differential operators on $z$-densities, thus generalizing a result
by Duflo \cite{du79} and
 Lichnerowicz \cite{lich}.\vspace{0.3cm}

Another problem posed by M.~Duflo in \cite{duflo-japon} is solved in
Section~\ref{sectionopdiff}: the expression in exponential
coordinates of invariant differential operators
(Theorem~\ref{theoecritureexp}).
Our solution is given in terms of Kontsevich's graphs.\vspace{0.3cm}

We define at the end of Section~\ref{sectionopdiff}  a deformation along the axis of  the Campbell-Hausdorff formula for symmetric pairs in the spirit of the Kashiwara-Vergne conjecture. We proved (Theorem \ref{theoKVsym} and Proposition~\ref{propKVsym}) that this deformation in the case of quadratic Lie algebras, considered as very symmetric quadratic pairs, implies the Kashiwara--Vergne conjecture. We conjecture that in the case of Lie algebras considered as symmetric pairs our $E$ function is equal identically to $1$. This conjecture would solve the Kashiwara-Vergne conjecture.\vspace{0.3cm}

In Section~\ref{sectionHomomorphismeHC} we consider the
Harish-Chandra homomorphism for symmetric pairs. Actually there are
two natural choices for a complement of $\k^\perp$, one by Cartan's
decomposition and the other by Iwasawa's. We show that these two
choices together lead to the Harish-Chandra homomorphism for general
symmetric pairs. In this language the Harish-Chandra homomorphism
consists of the restriction to the little symmetric pair in
Iwasawa's decomposition. The former decomposition and the
intertwinement then yield a formula for the Harish-Chandra
homomorphism in terms of graphs. It follows from this expression
that the Harish-Chandra homomorphism is invariant under the action
of the generalized Weyl group. We hope that this formula
will allow a resolution of the polynomial conjecture for symmetric pairs.\vspace{0.3cm}

Finally, in Section~\ref{sectionconstructioncaractere} we apply the
principles of bi-quantization \cite{CF1} to the case of triplets $f+
\b^\perp, \g^*, \k^\perp$, where $\b$ is a polarization for $f\in
\k^\perp$. These constructions produce characters for the algebras
of invariant differential operators
(Proposition~\ref{propconstructioncaractere}). It is a novel method
that we hope will be promising in other situations as
well\footnote{These methods may be applied
in certain cases of homogeneous spaces.}.\vspace{0.3cm}

In the case of polarizations in normal position we show, by a
homotopy on the coefficients (which makes use of an $8$ color form),
that the characters are independent of the choice of polarization
(Proposition~\ref{propindependance}). Thus we recover some classical
results
of the orbit method for Lie algebras.\vspace{0.3cm}

It follows that for symmetric pairs admitting $\sigma$-stable
polarizations
Rou\-vière's isomorphism computes the characters of the orbit method (Theorem~\ref{theosigmastable}).\vspace{0.3cm}

One can regard these new methods as a replacement for the orbit
method.

\selectlanguage{francais}

\section*{Introduction}

 Soit $(\g, \sigma)$ une paire symétrique, c'est-à-dire $\g$ est
une algèbre de Lie (quelconque) de dimension finie sur $\R$ et
$\sigma$ est une involution qui est un automorphisme d'algèbres de
Lie.  On note alors $\g=\k\oplus \p$ la décomposition relative à
$\sigma$, où $\k$ désigne l'espace propre associé à la valeur propre $+1$
et $\p$ l'espace
propre associé à la valeur propre $-1$.
Cette décomposition est aussi appelée \textit{décomposition de Cartan}.\\

Le théorème de Poincaré-Birkhoff-Witt (PBW) assure  la décomposition
de l'algè\-bre enveloppante $U(\g)$ :

\[U(\g)=U(\g)\cdot \k \oplus \beta
(S(\p)),\]où $\beta$ désigne la symétrisation de $S(\g)$ dans
$U(\g)$.\footnote{On a $\beta(X_1\ldots X_n)=\frac1{n!}\sum_{\sigma
\in S_n} X_{\sigma(1)}\ldots X_{\sigma(n)}$.} On peut alors
identifier $S(\p)$  l'algèbre symétrique de~$\p$ et
$U(\g)/U(\g)\cdot \k$ via la symétrisation $\beta$. En général
$U(\g)/U(\g)\cdot \k$ n'est pas une algèbre, mais les
$\k$-invariants $\big(U(\g)/U(\g)\cdot \k\big)^\k$ forment une algèbre.\\

Cette algèbre est un objet central car c'est l'algèbre des
opérateurs différentiels invariants sur l'espace symétrique $G/K$
associé à la paire symé\-trique $(\g, \sigma)$. Elle
 intervient de manière cruciale dans l'analyse harmonique sur les espaces
 symétriques.\\

 Cette algèbre est commutative \cite{lich, du79} et le théorème de
Poincaré-Birkhoff-Witt nous assure que
\[\Big(U(\g)/U(\g)\cdot \k\Big)^\k \quad \mathrm{et}\quad  S(\p)^\k\]
sont isomorphes comme espaces vectoriels. Remarquons que $S(\p)^\k$ est tout simplement le
gradué associé de $\Big(U(\g)/U(\g)\cdot \k\Big)^\k$.\\

 On  conjecture \cite{To3} que ces
deux algèbres sont isomorphes \textit{comme algèbres}, ce que le
second auteur a nommé la \textit{conjecture polynomiale}. Cette
conjecture généralise pour les paires symétriques l'isomorphisme de
Duflo pour les algèbres de Lie \cite{du77} entre le centre de
l'algèbre enveloppante $U(\g)^\g$
et les invariants dans l'algèbre symétrique $S(\g)^\g$.\\

Précisons dans cette introduction que pour un espace homogène $G/H$
quelconque (on note $\h$ l'algèbre de Lie de $H$) l'algèbre des
opérateurs différen\-tiels invariants sous l'action de $G$
s'identifie grâce à un résultat de Koornwinder \cite{koorn} à
l'algèbre $\big(U(\g)/U(\g)\cdot \h\big)^\h$. Cette algèbre n'est
pas commutative en général. Son gradué associé est alors une
sous-algèbre de Poisson de $\big(S(\g)/S(\g)\cdot \h\big)^\h$ (cette
dernière possède une structure de Poisson
naturelle).\\

S'il existe un supplémentaire de $\h$ qui soit invariant sous
l'action adjointe de $\h$,  notons   le  $\q$, alors une conséquence
facile de PBW est que $\big(U(\g)/U(\g)\cdot \h\big)^\h$ et
$\big(S(\g)/S(\g)\cdot \h\big)^\h$ sont encore isomorphes comme
espaces vectoriels. En général, on ne sait pas si cette propriété
reste vraie, pour la simple raison qu'il n'existe pas d'application
naturelle de $\big(S(\g)/S(\g)\cdot \h\big)^\h$ dans $U(\g)$ (ou un
quotient).\\

On conjecture toutefois que le centre de $\big(U(\g)/U(\g)\cdot
\h\big)^\h$ et le centre de Poisson de $\big(S(\g)/S(\g)\cdot
\h\big)^\h$ sont toujours isomorphes comme algèbres : c'est une
conjecture de M. Duflo \cite{duflo-japon}.\\

Peu de choses sont connues dans cette généralité. Dans le cas où
$\g$ est une algèbre nilpotente, des progrès sensibles ont été faits
ces dernières années par Corwin-Greenleaf \cite{CG},
Fujiwara-Lion-Magneron-Mehdi \cite{FLMM}, Baklouti-Fujiwara
\cite{BF} et Baklouti-Ludwig \cite{BL}. Dans le cas où $G$ et $H$
sont des groupes réductifs, F. Knop \cite{knop} donne une réponse
satisfaisante et remarquable à cette conjecture. Dans le cas ou $H$
est compact et $G=H\semi N$ est un produit semi-direct de $H$ par un
groupe d'Heisenberg~$N$, Rybnikov
 \cite{ryb} utilise  le   résultat de F. Knop pour conclure
 positivement à la conjecture de~Duflo.\\\\

Cet article  propose une approche nouvelle sur ces questions basée
sur la construction de Kontsevich \cite{Kont} et ses extensions  aux
cas des sous-variétés co-isotropes par Cattaneo-Felder \cite{CF1},
\cite{CF2}. Nous abordons ici la problématique des paires
symétriques mais nous pensons que nos méthodes ont un champ
d'applications plus vaste, notamment dans le cas nilpotent
homogène\footnote{Un des nos étudiants travaille sur ce cas.}.\\

On peut voir ce mémoire comme un pont  entre les méthodes de
quantification par  déformation  et la méthode des orbites en
théorie de Lie.\\

\paragraph{Résultats détaillés de l'article :}

La Section~\ref{sectionrappels} rappelle les constructions de
quantification de Kontsevich et l'extension au cas des variétés
co-isotropes par  Cattaneo-Felder.  On détaillera l'argumentation
sur la compatibilité en cohomologie. Cette section est utile pour
les experts en algèbre de Lie qui ne sont pas familiers avec les
constructions de quantification par déformation.\\

 On étudie ensuite dans le cas linéaire la dépendance de ces
constructions par rapport au choix du supplémen\-taire. Nous
montrons que les espaces de réduction sont isomorphes et nous
décrivons l'isomorphisme (Proposition~\ref{propdependance} et
Théorème~\ref{theodependance}): c'est un élément du
groupe de jauge qui s'obtient via la résolution d'une équation différentielle.\\
 \\

Dans la Section~\ref{ExempleReduction}, on précisera les graphes
qui interviennent dans le cas linéaire et on donnera trois exemples
fondamentaux d'espaces de réduction  que l'on peut rencontrer en
théorie de Lie :  le cas des paires symétriques (Proposition~\ref{propreductionsym}), le cas des décompositions d'Iwasawa
(Proposition~\ref{propreductioniwasawa}) et le cas des polarisations
(Proposition~\ref{propreductionpolarisation}). Ces exemples sont
nouveaux et démontrent que ces nouvelles méthodes de
quantification sont adaptées à la théorie de~Lie.\\

 Rappelons que dans une série d'articles F. Rouvière \cite{Rou90,
Rou91, Rou94} a introduit pour les paires symétriques, suivant les
méthodes de Kashiwara-Vergne \cite{KV},
 une fonction mystérieuse $e(X,Y)$ définie pour $X, Y \in \p$.  Cette fonction calcule
 à une conjugaison près le star-produit :

\begin{equation}P\underset{Rou}\sharp Q= \beta^{-1}\Big(\beta(P)\cdot
 \beta(Q) \quad \mathrm{modulo} \quad U(\g)\cdot \k \Big)\end{equation}
  pour
  $P, Q$ dans $S(\p)^\k$.\\

Dans la Section~\ref{sectionfonctionE}, on définira une fonction
 $E(X,Y)$ pour $X, Y \in \p$ en termes de graphes, qui aura des
 propriétés analogues à la fonction $e(X,Y)$ de Rouvière. La
  comparaison de ces deux fonctions est un point important de cet
  article.    Cette fonction $E(X, Y)$ exprime le star-produit
  $\underset{CF} \star$ de
Cattaneo-Felder dans le
 cas des paires symétriques.
 On déduira,
 à partir des graphes intervenant dans la construction de cette
 fonction, une propriété de symétrie (Lemme~\ref{symetrieE}) et des propriétés
 supplémentaires
 dans le cas résoluble (Proposition~\ref{propE=1resoluble}),
 dans le cas des paires symétriques
 d'Alekseev-Meinrenken (Proposition~\ref{propE=1AM}) ou
  le cas très symétrique quadratique  (Proposition~\ref{propE=1tressym}).
  Dans tous ces cas on montre que la fonction $E$ vaut identiquement
  $1$.
Ces propriétés remarquables et élémentaires
 donnent des démonstrations nouvelles et unifiées de résultats démontrés par
 Rouvière dans le cas résoluble
 (Théorème~\ref{theoResoluble} et  Proposition~\ref{propordre4})
 et généralisent un théorème d'Alekseev-Meinrenken (Théorème~\ref{theoAM}) dans le cas quadratique (anti-invariant).\\

Dans la Section~\ref{sectionopdiff}, on montrera que le star-produit
$\underset{CF} \star$ de Cattaneo-Felder et celui de Rouvière
coïncident (Théorème~\ref{calculCF}). Ce résultat est nouveau. On en
déduira  (Théorème~\ref{theoz}) la commutativité  pour tout $z\in
\mathbb{R}$ des algèbres d'opérateurs différentiels invariants sur
les $z$-densités
 $\Big(U(\g)/U(\g)\cdot \k^{z\tr_k}\Big)^\k$ généra\-lisant
 ainsi un résultat de Duflo~\cite{du79} et
 Lichnerowicz~\cite{lich}.\\

 Un  autre problème proposé par M. Duflo dans
\cite{duflo-japon} est résolu en section \S~\ref{sectionopdiff} :
l'écriture, en coordonnées exponentielles, des opérateurs
 différentiels invariants (Théorème~\ref{theoecritureexp}).
 Notre réponse s'exprime en termes de graphes de Kontsevich.\\

En fin de Section~\ref{sectionopdiff}, on définit une déformation le long des axes de la formule de Campbell-Hausdorff pour les paires symétriques, dans l'esprit de la conjecture de Kashiwara-Vergne. On montre (Théorème~\ref{theoKVsym} et  Proposition~\ref{propKVsym}) que dans le cas des algèbres de Lie quadratiques, considérées comme des paires très symétriques quadratiques, notre déformation implique la conjecture de Kashiwara-Vergne. Plus généralement on conjecture que notre fonction $E$ vaut $1$ dans le cas des algèbres de Lie, considérées comme des paires symétriques. Cette conjecture implique alors la conjecture de Kashiwara-Vergne.\\

 Dans la Section~\ref{sectionHomomorphismeHC}, on s'intéresse à
l'homomorphisme d'Harish-Chandra pour les paires symétriques. En
effet il existe deux choix de supplémentaires de $\k^\perp$,
essentiellement celui donné par la décomposition de Cartan et
l'autre donné par la décomposition d'Iwasawa. On montre que ces deux
choix condui\-sent à l'homomorphisme d'Harish-Chandra pour les
paires symétriques géné\-rales. Dans ce langage, l'homomorphisme
d'Harish-Chan\-dra consiste en la restriction à la petite paire
symétrique dans la décom\-position d'Iwasawa. La première
décomposition et l'entrelacement donneront alors une formule pour
l'homomorphisme d'Harish-Chandra en termes de graphes. On en
déduit au vu de l'expression que l'homomorphisme d'Harish-Chandra
généralisé est invariant par l'action du groupe de Weyl généralisé.
On espère que cette formule permettra de résoudre la
conjecture polynomiale pour les paires symétriques.\\

 Enfin dans la Section~\ref{sectionconstructioncaractere}
  de cet article, on appliquera le principe de bi-quantification \cite{CF1}
au cas des triplets $f+ \b^\perp, \g^*, \k^\perp$ o\`u $\b$ est une
polarisation pour $f\in \k^\perp$. Ces constructions fournissent des
caractères pour les algèbres d'opérateurs différentiels invariants
(Proposition~\ref{propconstructioncaractere}). C'est une nouvelle
méthode, que l'on espère prometteuse dans d'autres situations\footnote{Ces méthodes s'appliquent dans certains cas d'espaces homogènes.}.\\

Dans le cas o\`u les polarisations sont en position d'intersections
normales, on montre,  par une méthode d'homotopie sur les
coefficients (faisant intervenir une forme à $8$ couleurs), que
les caractères sont indépendants du choix des polarisations
(Proposition~\ref{propindependance}). On retrouve ainsi des
résultats classiques de la méthode des orbites dans le cas des algèbres de Lie.\\

On en déduit, pour les paires symétriques qui admettent des
polarisations $\sigma$-stables,  que l'isomorphisme de Rouvière
calcule les caractères de la méthode des orbites (Théorème~\ref{theosigmastable}).\\

On peut voir ces nouvelles méthodes comme un substitut à la
méthode des orbites.
\tableofcontents
\vspace{0,3cm}

\textbf{Bibliographie} \hfill \textbf{76}
\newpage
\section{Rappels sur la construction de
Cattaneo-Felder}\label{sectionrappels}

 Soit $\g$ une algèbre de Lie de dimension finie sur $\R$.
L'espace dual $\g^*$ est alors muni d'une structure de Poisson
linéaire. On note
 $\pi$ le bi-vecteur de Poisson associé.\\

Supposons donnée sous-algèbre $\h$  de $\g$.  Son orthogonal
$\h^\perp$ est alors une sous-variété co-isotrope de $\g^*$. Dans
\cite{CF1}, Cattaneo-Felder décrivent une construction pour une
quantification de $(S(\g)/S(\g)\cdot\h)^\h$ ou d'une sous-algèbre de cette dernière.\\
\subsection{Quantification}
Les constructions ne sont pas intrinsèques et dépendent du choix
d'un sup\-plémentaire de $\h$ dans $\g$. Notons $\q$ un tel
supplémentaire. On peut alors identifier $\h^*$ avec $\q^\perp$.\\

 Deux constructions sont données par Cattaneo et Felder, l'une en
termes de série de type Feynman avec  diagrammes colorés et l'autre
en terme de transformée de Fourier partielle. Ces constructions sont
locales mais  on peut les globaliser. On va rappeler ces dernières
dans le cadre qui nous intéresse à savoir le cas des sous-algèbres,
mais il suffira de remplacer $\h^\perp$ par une sous-variété $C$
pour obtenir la construction plus générale.

\paragraph{Rappel des constructions :}

La variété qui intervient dans cette construction est une
super-variété  : \begin{equation}M :=\h^\perp\oplus \Pi \h,\end{equation}où $\Pi$ désigne le
foncteur de  changement de parité. L'algèbre des fonctions est donc
\begin{equation}\mathcal{A} := \mathcal{C}^\infty(\h^\perp)\otimes
\bigwedge(\g^*/\h^\perp)\simeq \mathcal{C}^\infty(\h^\perp)\otimes
\bigwedge \h^*.\end{equation}
Si $C=\h^\perp$, la fibre du fibré normal $N_C$
vaut alors $ \g^*/\h^\perp$,
celle du fibré cotangent (conormal) $T_C^{\perp}=N_C^*$ vaut $\h$.\\

\subsection{Construction en terme de transformée de Fourier}

La dg-algèbre des poly-champs sur $M$ est l'algèbre symétrique
décalée de l'algèbre des dérivations de $\mathcal{A}$
\begin{equation}\mathcal{T}(A): = S_\mathcal{A}(Der(\mathcal{A})[-1])[1].\end{equation}

Par transformée de Fourier dans la fibre impaire (\textit{cf.} plus loin,)
 cette dg-algèbre est isomorphe à l'algèbre des poly-champs de
 vecteurs formels le long de $\h^\perp$, c'est-à-dire

\begin{equation}\mathcal{T}(\mathcal{B}): = S_\mathcal{B}(Der(\mathcal{B})[-1])[1],\end{equation} où $\mathcal{B}=\mathcal{C}^\infty(\h^\perp)\otimes
S((\g^*/\h^\perp)^*)= \mathcal{C}^\infty(\h^\perp)\otimes S(\h)$
est l'algèbre des fonctions sur \begin{equation}\widehat{M}:=\h^\perp\oplus \g^*/\h^\perp\end{equation}
polynomiales dans la fibre.

Une fois que l'on a fixé un voisinage tubulaire
de $\h^\perp$, on peut  identifier ce voisinage  avec
le fibré normal et on peut identifier $\g^*/\h^\perp$ avec un supplémentaire de
$\h^\perp$. On a besoin ici de choisir un supplémentaire de $\h$. \\

Les fonctions au voisinage de $\h^\perp$ s'identifient avec leur
développement de Taylor partiel dans la direction normale, c'est-à-dire des éléments de $\mathcal{B}$. L'algèbre
$\mathcal{T}(\mathcal{B})$ s'identifie alors aux poly-champs de
vecteurs formels\footnote{On entend par poly-champs de vecteurs dans
un voisinage formel de $\h^\perp$ une complétion des
$\mathcal{V}(\g^*)/I_{\h^\perp},$ où $\mathcal{V}(\g^*)$ désigne les
poly-champs de vecteurs sur $\g^*$ et  $I_{\h^\perp}$ désigne les
poly-champs avec développement de Taylor
nul sur $\h^\perp$.} le long de $\h^\perp$.\\

Le théorème de Formalité \cite{Kont}, nous dit que la dg-algèbre
$\mathcal{T}(\mathcal{A})$ est $L_\infty$ quasi-isomorphe à la
dg-algèbre $\mathcal{D}_{poly}(\mathcal{A})$. En combinant la
transformée de Fourier et le théorème de Formalité dans le cas des
super-espaces, Cattaneo-Felder démontrent le théorème suivant.

\paragraph{Théorème~(Cattaneo-Felder \cite{CF2})} \textit{La dg-algèbre des
poly-champs de vecteurs formels le long de $\h^\perp$ est $L_\infty$
quasi-isomorphe à la dg-algèbre
$\mathcal{D}_{poly}(\mathcal{A})$.}\\

Si on dispose d'une solution de l'équation de Maurer-Cartan sur
$\g^*$, alors on disposera d'une solution de Maurer-Cartan sur $M$
(par transformée de Fourier) puis d'une solution de Maurer-Cartan dans
$\mathcal{D}_{poly}(\mathcal{A})$, c'est-à-dire en général (\textit{cf.} plus
loin) d'une $A_\infty$-structure sur $\mathcal{A}$.

\paragraph{Cas des paires symétriques :} Dans le cas des paires
symétriques le  supplémen\-taire $\p$ étant canonique, on
dispose des équations suivantes.\\

Notons $(K_i)_i$ une base de $\k$, $(P_j)_j$ une base de $\p$. On
note $(K_i^*)_i, (P_j^*)_j$ les bases duales. On identifie $\k^\perp$
et $\p^*$. On identifie de même $\g^*/\k^\perp=\k^*$ et
$\p^\perp$. La super-variété
étant
\[ M :=\k^\perp\oplus \Pi \k=\p^*\oplus \Pi \k\]
on aura \[\widehat{M} :=\k^\perp \oplus \g^*/\k^\perp= \p^*\oplus
\k^*.\]  On  notera $\theta_i := \Pi K_i^*$ les fonctions de
coordonnées sur $\Pi \k$.  La dérivée dans la direction $K_i^*$ sera
notée $
\partial_{K_i^*}$; c'est un champ de vecteurs constant sur $\widehat{M}$. De
même la dérivée dans la direction  $\Pi K_i$ sera  notée
$\partial_{\Pi K_i}$; c'est un  champ de vecteurs constant sur
$M$.\\

La transformée de Fourier (avec changement de parité) change la dérivée
$\partial_{K_i^*}$ (c'est un champ de degré impair de
$\mathcal{T}(\mathcal{B})$) en la fonction $\theta_i= \Pi K_i^*$
(c'est une variable impaire de $\mathcal{A}$). De même $K_i$
fonction de coordonnée sur $\widehat{M}$  (variable paire de
$\mathcal{B}$) est changée en $\partial_{\Pi K_i}$ (variable paire de $\mathcal{T}(\mathcal{A})$).\\

\noindent  Le bi-vecteur de Poisson sur $\g^*$ associé à la
structure de
Poisson, s'écrit concrètement\footnote{On utilise la convention  $\partial_1\wedge \ldots \wedge \partial_n=\frac1{n!}\sum_{\sigma\in S_n} \partial_{\sigma(1)}\otimes \ldots \otimes\partial_{\sigma(n)}$. Les indices répétés sont sommés.}:

\begin{equation}\pi= [K_i, K_j] \partial_{K_i^*}\wedge \partial_{K_j^*} + [P_i,
P_j]
\partial_{P_i^*}\wedge \partial_{P_j^*} + 2[K_i, P_j] \partial_{K_i^*}\wedge
\partial_{P_j^*}.\end{equation}

\noindent Sa transformée de Fourier partielle, notée
$\widehat{\pi}$, est somme d'un champ de vecteurs et d'un
$3$-vecteur sur la variété $\k^\perp
\oplus \Pi \k$ :

\begin{equation}\overset{\wedge}{\pi}= \theta_i \theta_j\partial_{\Pi[K_i, K_j]} +
2[K_i, P_j] \theta_i
\partial_{P_j^*}+ \partial_{P_i^*}\wedge \partial_{P_j^*}\wedge \partial_{\Pi[P_i, P_j]}.\end{equation}

\noindent Comme on le constate, la
variété $\k^\perp \oplus \Pi \k$ ne porte pas de structure de
Poisson (car $\overset{\wedge}{\pi}$ n'est pas un $2$-vecteur),
mais une structure vérifiant l'équation de Maurer-Cartan et qui est homogène si
on tient compte de tous les degrés impairs dans
$\mathcal{T}(\mathcal{A})$,
c'est-à-dire des variables $\theta_i$ et $\partial_{P_i^*}$.

\paragraph{Remarque 1 : }La partie $1$-champ de $\overset{\wedge}{\pi}$ , est clairement
associée à la paire symétrique dégénérée (abélianisée, \textit{ie.} on a
$[\p, \p]=0$) produit semi-direct de $\k$ et $\p$:  $\k \semi
\p$. Ce $1$-champ est de carré nul. Compte tenu de la graduation, on
peut voir la paire symétrique comme une déformation de la partie
abélianisée. Si on considérait la structure déformée  formelle
$[x,y]=t^2 [x,y]$ pour $x\in \p$ et $y \in \p$ (les autres crochets
restant inchangés) on trouverait
\[\overset{\wedge}{\pi_t}= \theta_i \theta_j\partial_{\Pi[K_i, K_j]} +
2[K_i, P_j] \theta_i
\partial_{P_j^*}+ t^2\partial_{P_i^*}\wedge \partial_{P_j^*}\wedge \partial_{\Pi[P_i, P_j]}.\]\\
Pour $t=0$ on retrouve la paire symétrique abélianisée.

\paragraph{Remarque 2 :} L'utilisation de la transformée de Fourier impaire est bien connue pour les
algèbres Lie; on  considère $\Pi\g$ muni du $1$-champ impair
quadratique
\[Q=\xi_i\xi_j\partial_{\Pi[e_i,e_j]}.\] Il vérifie $Q^2=\frac
12[Q,Q]=0$. On peut alors appliquer la quantification de Kontsevich
dans ce contexte. Les formules de Kontsevich sont plus simples et on
peut par exemple décrire de manière plus naturelle le $L_\infty$
quasi-isomorphisme tangent qui, par l'argument d'homotopie,
réalisera dans le cas des algèbres de Lie, l'isomorphisme de Duflo
généralisé \[H(\g, S(\g)) \underset{\mathrm{alg\grave{e}bre}}\sim
H(\g, U(\g)).\]Cette méthode a été
expérimentée par Shoikhet \cite{Shoi} (voir aussi \cite{pev-toro}).

\paragraph{Remarque 3 :} Le $L_\infty$ morphisme, agit aussi
sur les poly-champs dont la restriction à $\h^\perp$ est nulle sur
$\bigwedge \h$, c'est-à-dire les poly-champs de vecteurs dont la
restriction est nulle sur la puissance extérieure du fibré conormal
; ce sont les poly-vecteurs dits relatifs; par exemple le bi-vecteur
de Poisson $\pi$ vérifie cette
propriété (voir \cite{CF1}).\\

\subsection{Construction en termes de diagrammes de Feynman}
La formule proposée est semblable à celle de Kontsevich \cite{Kont}
pour $\R^n$. Elle fait intervenir des graphes numérotés, des
coefficients obtenus par intégration de formes différentielles sur
des variétés de configurations dans le demi-plan de Poincaré et des
opérateurs poly-différentiels associés à ces graphes \footnote{On
pourra consulter les références \cite{AMM}, \cite{CKT}
pour une description détaillée de la construction de Kontsevich.}. Pour simplifier la lecture de cet article on rappelle brièvement les ingrédients de la construction générale.\\

\subsubsection{Variétés de configurations}
On note $C_{n,m}$ l'espace des configurations de $n$ points distincts
dans $\mathcal{H}$ le
demi-plan de Poincaré (points de première espèce ou points aériens) et $m$ points  distincts sur la droite réelle
(points de seconde espèce ou points terrestres), modulo l'action
du groupe $az+b$ (pour $a \in \R^{+*}, b\in \R$).  Dans son article \cite{Kont}, Kontsevich construit des
compactifications de ces variétés notées $\overline{C}_{n,m}$. Ce
sont des variétés à coins de dimension $2n-2+m$. Ces variétés ne sont pas connexes pour $m\geq 2$.
On notera par  $\overline{C}^{+}_{n,m}$ la composante  qui contient les configurations où les
 points terrestres sont ordonnés dans l'ordre croissant (\textit{ie.} on  a $\ov{1}< \ov{2}<\cdots< \ov{m}$).

\subsubsection{Fonctions
d'angle à deux couleurs}\label{defifonctionangle2couleurs}

Les graphes de Kontsevich vont être colorés en fonction de la
variable de dérivation  associée dans l'opérateur, on aura donc besoin d'une fonction d'angle qui dépend de deux couleurs notées $+$ et $-$.\\

\begin{defi}\label{defiangle} On définit les deux fonctions d'angles de $C_{2,0}$ dans  $\mathbb{S}^1$
\begin{equation}
\left\{\begin{array}{cccc}
 \phi_{+}(p,q) = & \underset{p}\bullet\longrightarrow \underset{q}\bullet=& \overset{\longrightarrow}{\phi}(p,q):=  & \arg(p-q) + \arg(p-\overline{q}) \\
  \phi_{-}(p,q)= & \underset{p}\bullet\dashrightarrow\underset{q}\bullet= & \overset{\dashrightarrow}{\phi}(p,q):=  & \arg(p-q) - \arg(p-\overline{q}).
\end{array}
\right.
\end{equation}
\end{defi}
Ces fonctions  d'angle s'étendent à la compactification $\overline{C}_{2,0}$. La fonction d'angle $ \phi_{+}$ sera associée aux variables tangentes (\textit{ie.} dans $ \h^\perp$ ) tandis que la fonction d'angle $ \phi_{-}$ sera associée aux variables normales (\textit{ie.} dans $\g^*/\h^\perp=\h^*$)\footnote{On a besoin ici de faire un choix d'un sup\-plémentaire de $\h$.}. Remarquons que l'on a
\begin{equation}\mathrm{d} \phi_{-}(p,q)= \mathrm{d} \phi_+(q,p).\end{equation}
Le fonction d'angle $\phi_+$ est celle définie par Kontsevich.

\subsubsection{Graphes et opérateurs
différentiels associés}\label{grapheetoperateur}

Les graphes qui vont intervenir dans le définition du $L_\infty$ quasi-isomorphisme sont
analogues à ceux de Kontsevich \cite{Kont} avec la différence
essentielle suivante, les arêtes sont colorées par nos deux couleurs $+$ (dans les dessins $\longrightarrow$) et $-$ (dans les dessins $\dashrightarrow$). \\

Soit  $\G$ est un graphe (quiver) avec $n$ sommets de première espèce (aériens) numérotés $1, 2, \ldots, n$ et $m$ sommets de seconde espèce (terrestres) numérotés $\ov{1}, \ldots, \ov{m}$.

\begin{figure}[h!]
\begin{center}
\includegraphics[width=6cm]{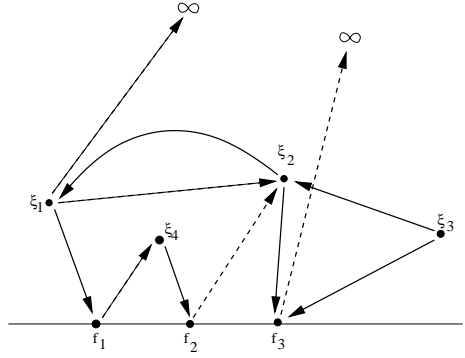}
\caption{\footnotesize Graphe type intervenant dans le calcul de
$\mathcal{U}_4$}\label{grapheU4}
\end{center}
\end{figure}

Par construction $\G$ n'a pas de boucles, ni d'arêtes doubles (même source, même but, même couleur). Les arêtes, issues des sommets de seconde espèce, sont colorées par la couleur $-$ et ne reçoivent que des arêtes colorées par la couleur $+$. Par ailleurs un certain nombre d'arêtes, portant la couleur $-$,  n'ont pas de but. On dira qu'elles vont à l'infini (\textit{cf.} Fig.~\ref{grapheU4}). \\

Notons $k_1, \ldots k_n$ le nombre d'arêtes sortant des sommets de première espèce et $k_{\ov{1}}, \ldots, k_{\ov{m}}$ le nombre d'arêtes sortant des sommets de seconde espèce. Soient $\xi_1, \ldots, \xi_n$ des poly-vecteurs,
avec $k_i$ le degré de $\xi_i$ et $f_{\ov{1}}, \ldots, f_{\ov{m}}$ des fonctions  avec $f_{\ov{j}} \in \mathcal{C}^\infty(\h^\perp)\otimes \bigwedge^{k_{\ov{j}}} \h^*$. En plaçant au sommet aérien $i$ le poly-vecteur $\xi_i$ et au sommet terrestre $\ov{j}$ la fonction $f_{\ov{j}}$, on définit après restriction à $\h^\perp$, une fonction (poly-champ)

$$B_\G(\xi_1, \ldots,
\xi_n)\big(f_{\ov{1}}, \cdots, f_{\ov{m}}) \in \mathcal{C}^\infty(\h^\perp)\otimes \bigwedge^{p}\h^*,$$
où $p$ désigne le nombre d'arêtes qui partent à l'infini. La règle complémentaire est la suivante : si l'arête porte la couleur $+$, on ne dérive que selon les variables tangentes (\textit{ie.} dans $\h^\perp$) et si la couleur est $-$, on ne dérive que selon les variables normales (\textit{ie.} dans $\h^*$). On définit ainsi un opérateur poly-différentiel $B_\G(\xi_1, \ldots,
\xi_n)$. Remarquons que les arêtes, qui partent à l'infini, contribuent à la définition de cet opérateur.

\subsubsection{Coefficients}\label{coefficient}

Soit $\G$ un graphe coloré avec $n$ sommets de première espèce et $m$ sommets de seconde espèce. On dessine le graphe dans $C_{n, m}$.  Toute arête colorée $e$, qui ne part pas à l'infini, définit par restriction
une fonction d'angle notée $\phi_{e}$  sur la variété $\overline{C}^+_{n,m}$. Si la couleur est $\epsilon\in \{+,- \}$ on choisit la fonction d'angle $\phi_\epsilon$. On note $E_{\Gamma}$
l'ensemble des arêtes du graphe $\Gamma$ qui ne partent pas à l'infini. Le produit ordonné $\Omega_{\Gamma}=\bigwedge_{e \in E_{\Gamma}}\mathrm{d}\phi_{e}$
est donc une $\sharp E_\G$-forme régulière sur  $\overline{C}^+_{n,m}$, variété compacte de dimension $2n+m-2$.
\begin{defi} Le poids, associé à un graphe coloré $\G$, est défini par l'intégrale \begin{equation}
w_{\Gamma}=\frac{1}{(2\pi)^{\sharp E_\G}}\int_{\overline{C}^+_{n,m}} \Omega_{\Gamma}.
\end{equation}
\end{defi}
Ce coefficient est nul si $\sharp E_\G \neq  2n+m-2$. Remarquons que les arêtes qui partent (resp. qui arrivent) de l'axe réel porte la couleur $-$ (resp. $+$), par conséquent les différentielles des  fonctions d'angle associées ne sont pas nulles.

\subsubsection{Construction du $L_\infty$quasi-isomorphisme}
Pour $n\geq 0$ on note
\begin{equation}\mathcal{U}_n(\xi_1, \ldots, \xi_n):=\sum_{\G} w_\G
B_\G(\xi_1, \ldots, \xi_n),\end{equation}
où la somme porte sur tous les graphes avec
$n$ sommets aériens  sur lesquels on a placés les polyvecteurs $\xi_1,
\ldots, \xi_n$ et un nombre quelconque de sommets terrestres.\\

\noindent L'opérateur $\mathcal{U}_n(\xi_1, \ldots, \xi_n)$ n'est pas
homogène pour la graduation de Hochschild, car la diminution du
nombre d'arguments peut être compensée par l'augmentation du degré
en~$\h^*$. Pour $n=0$, nécessairement $m$ vaut $2$ et on retrouve la multiplication dans
$\mathcal{C}^\infty(\h^\perp)\otimes
\bigwedge \h^*$. \\

\paragraph{Théorème~: (Cattaneo-Felder)}

\textit{La somme  $\mathcal{U}= \sum_{n\geq1} \frac 1{n!}\mathcal{U}_n$
définit un $L_\infty$ quasi-isomorphisme de la  $dg$-algèbre des
poly-champs de vecteurs formels dans un voisinage de $\h^\perp$ dans
$\mathcal{D}_{poly}\big(\mathcal{C}^\infty(\h^\perp)\otimes
\bigwedge
\h^*\big)$.}\\

 Comme dans le cas classique c'est la formule de Stokes qui
fournit  l'équation du $L_\infty$
quasi-isomorphisme.  \\

\subsection[Cas des bi-vecteurs de Poisson et cas linéaire]{Cas des bi-vecteurs de Poisson et  cas des bi-vecteurs de Poisson
linéaires}\label{quantification}

\subsubsection{Cas des bi-vecteurs de Poisson}\label{quantificationpoisson} Lorsqu'on applique
la construction précédente  dans le cas où les poly-vecteurs $\xi_i$
sont égaux à un $2$-vecteur de Poisson $\pi$, on trouve un opérateur
polydifférentiel formel non homogène du complexe de Hochschild. Mais
cet opérateur est homogène de degré $1$ si l'on
tient compte des degrés impairs.\\

Dans  cas des variétés de Poisson, étudié dans \cite{Kont},
l'opérateur est de degré~$1$ dans le complexe de Hochschild, c'est
donc  un produit associatif. Ici la situation est plus
compliquée. \\

En d'autres termes, comme la graduation tient compte du degré dans
les variables impaires, la structure obtenue est en fait une
$A_\infty$-structure sur l'espace
$\mathcal{A}=\mathcal{C}^\infty(\h^\perp)\otimes \bigwedge \h^*$
avec premier terme non nul \textit{a priori}, c'est-à-dire une
structure\footnote{On note $[\; ,\;]_{GH}$ le crochet de
Gerstenhaber.}
\[\mu=\mu_{-1}+ \mu_0+ \mu_1+ \mu_2+ \ldots \]vérifiant $\frac 12[\mu,
\mu]_{GH}=0$ ou en terme de bar-construction $ \overline{\mu}\circ
\overline{\mu}=0$, o\`u $\overline{\mu}$ désigne la codérivation de
l'algèbre tensorielle associée
à $\mu$.\\

\textit{A priori}, il existe  donc des composantes en tout degré (de
Hochschild mais ces composantes sont homogènes de degré $1$ si l'on
tient compte de
la graduation~$\bigwedge \h^*$); on a donc en général :  \\

\begin{itemize}

\item Une composante $\mu_{-1}$, que l'on note aussi $F_\pi$ (degré
$-1$ dans le complexe de Hochschild, c'est-à-dire un élément de
$\mathcal{A}$) qui ne prend pas d'argument: c'est une sorte de
courbure.

\item Une composante $\mu_0$  qui ne prend qu'un argument que l'on note
aussi $ A_\pi$ (degré $0$ dans le complexe de Hochschild, c'est-à-dire un opérateur différentiel formel sur $\mathcal{A}$) : c'est
presque une différentielle.

\item Une composante $\mu_1$  qui prend deux arguments que l'on note
aussi $ B_\pi$ (degré $1$ dans le complexe de Hochschild, c'est-à-dire un opérateur bi-différentiel formel sur $\mathcal{A}$) : c'est
presque un produit (associatif).

\item etc...

\end{itemize}

On~écrit alors $\mu=F_\pi + A_\pi +B_\pi + B_2 + \ldots .$

\paragraph{Faisons l'inventaire des arêtes: } Un graphe avec $n$
sommets aériens  et $m$ points terrestres produit (dans sa partie aérienne) $2n$ arêtes pour
une dimension de variété de configurations $2n+m-2$. Il faut donc
disposer, sur l'axe réel, $m$ poly-vecteurs fournissant au moins $m-2$
arêtes.  Si au total on a plus de $2n+m-2$ arêtes, alors il faut
faire sortir à l'infini le nombre adéquat d'arêtes, colorées par
$\h^*$.\\

Dans le cas $m=0$ il faut donc $2$ arêtes sortantes, donc $F_\pi$
est un $2$-vecteur. On remarquera que $F_\pi$ est alors de
degré $1=-1+2$. C'est une sorte de courbure.\\

Toutefois si $F_\pi=0$, alors  $A_\pi$ sera une différentielle
vérifiant $[A_\pi, B_\pi]_{GH}=0 $ et on aura l'équation\footnote{On
fera attention que l'on n'a pas $[B_2, A_\pi]_{GH}=-[A_\pi,
B_2]_{GH}$, car il faut tenir compte de la graduation des
coefficients, on a donc $[B_2, A_\pi]_{GH}=[A_\pi, B_2]_{GH}$.}
\[[B_\pi, B_\pi]_{GH}+ 2[B_2, A_\pi]_{GH}=0,\] où $B_2$ désigne la composante de
degré $2$ (qui prend $3$ arguments). En conséquence $B_\pi$
sera un produit associatif dans l'espace de cohomologie défini par $A_\pi$.\\

\subsubsection{Cas linéaire}
 On se place  dans le cas où $\pi$, le bi-vecteur de Poisson
placé aux sommets aériens, est égal à la moitié\footnote{Les graphes "géométriques" apparaissent de nombreuses fois à cause de la numérotation.} du  bivecteur de Poisson
linéaire associé à l'algèbre de Lie $\g$. On suppose toujours que
$\h$ est une sous-algèbre. Toutes les constructions se restreignent
à l'algèbre des  fonctions polynomiales notée
$C_{poly}(\h^\perp)\otimes \bigwedge \h^*$.

\paragraph{Nullité de la courbure : }
\begin{lem}\label{courburenulle}Dans le cas linéaire, pour toute sous-algèbre $\h$, on a
$F_\pi=0$.\end{lem}

\noindent \textit{\bf Preuve :} On doit avoir $2n-2$ arêtes dans le
graphe mais on ne dispose que de $n$ sommets à dériver. On doit
donc avoir par linéarité  de $\pi$, $2n-2\leq n$
c'est-à-dire  $n=1$ ou $n=2$.\\

-- Si $n=1$  la restriction de  $\pi$  à $\bigwedge^2 \h$ est nul si
et seulement si $\h$ est une sous-algèbre (on retrouve la condition de
co-isotropie de la variété $C=\h^\perp$).\\

-- Si $n=2$ les graphes qui interviennent sont comme dans
(\ref{roue2.eps}) (dessin de droite) car le $2$-vecteur est linéaire. Dans ce cas
le coefficient est nul, car on a une arête double ou deux arêtes qui
se suivent\footnote{Remarquons que pour $n=1$ ou $n=2$ les
contributions sont toujours
 nulles dans $F_\pi$, même si $\pi$ n'est pas linéaire.}.  \fin\\
\begin{equation}\label{roue2.eps}
\begin{array}{ccc}
  \includegraphics[width=6cm]{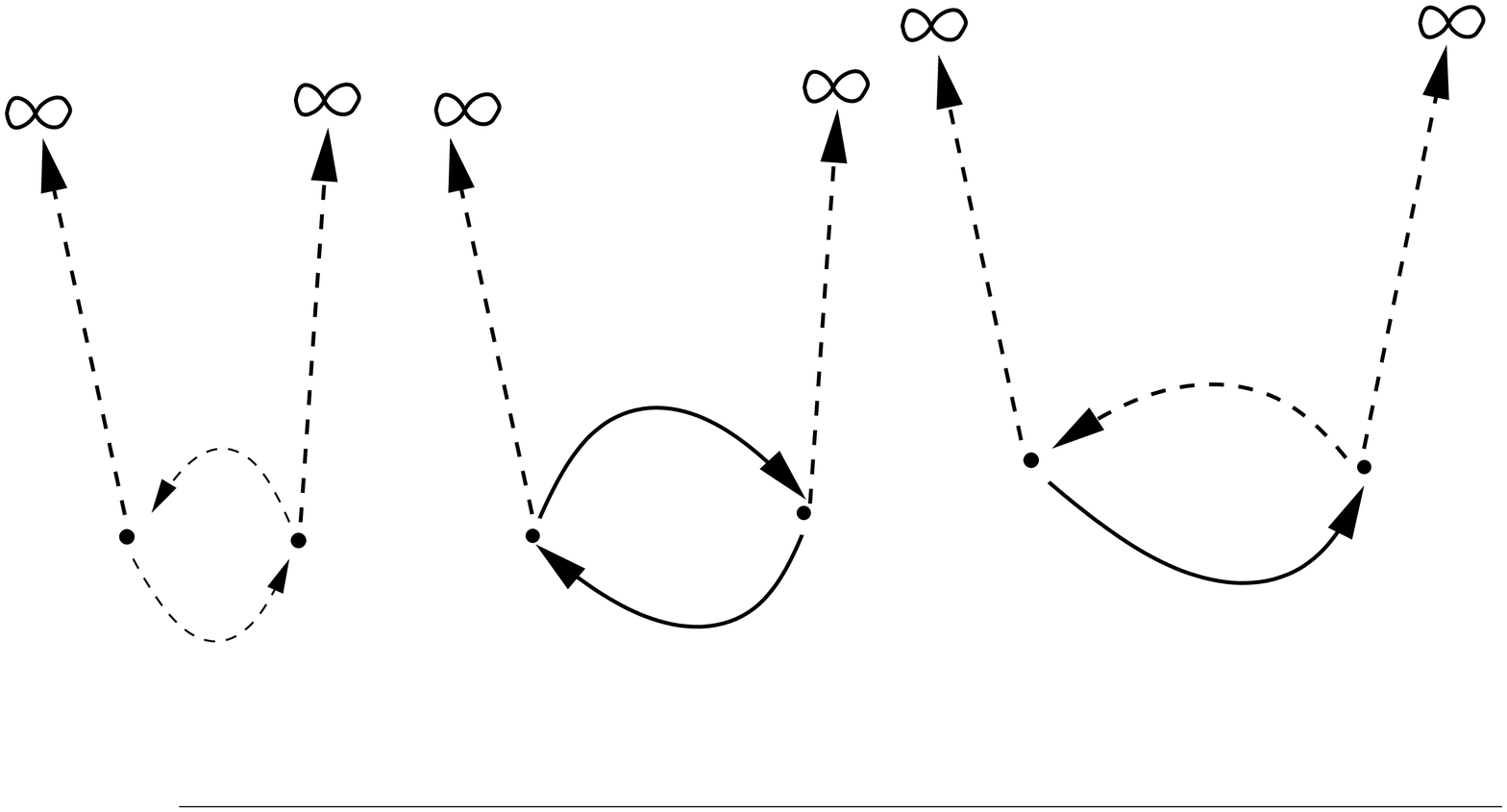} &\quad  &\includegraphics[width=6cm]{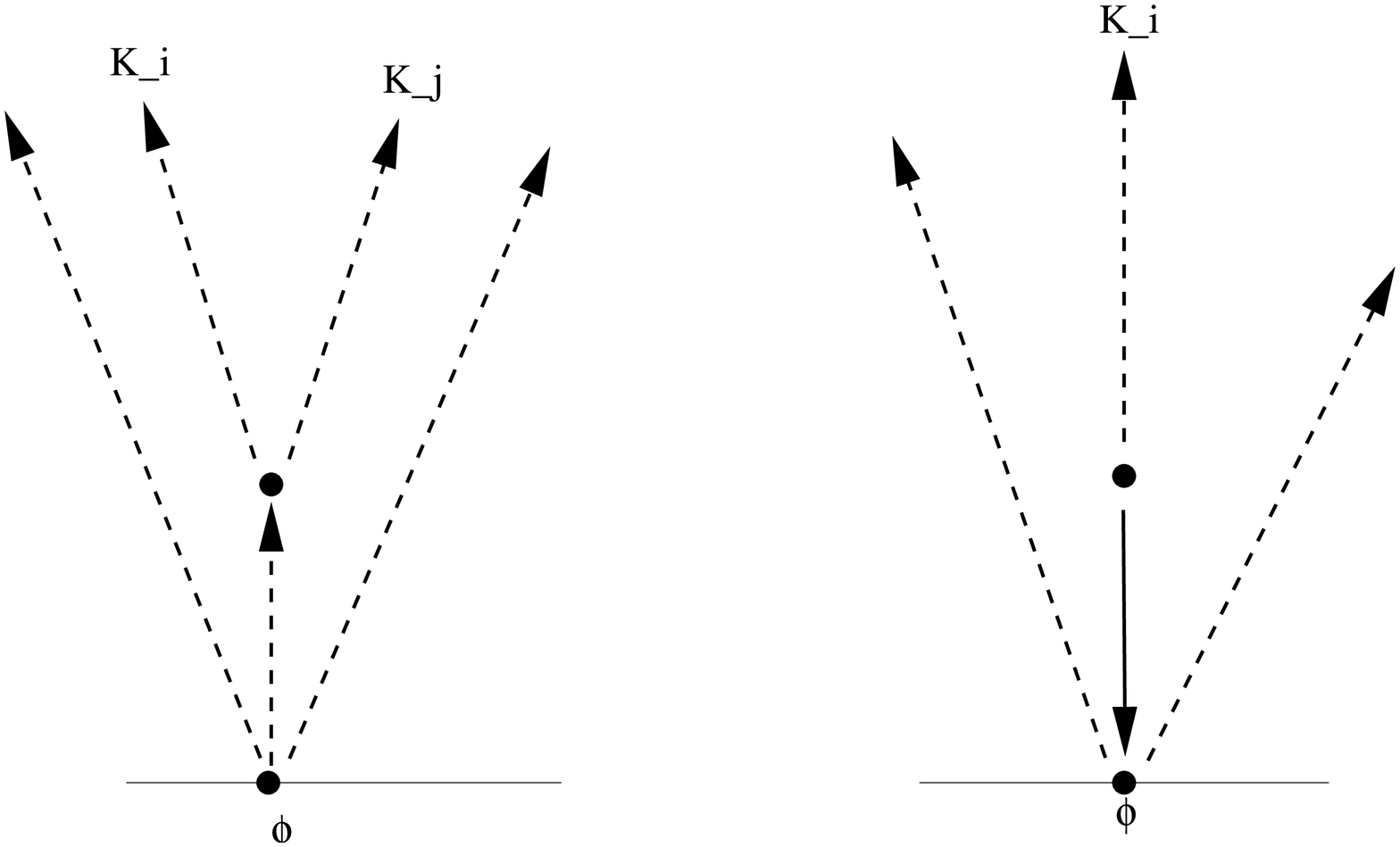} \\
  Graphes \;  pour \;
F_\pi & \quad & Graphes \; pour \; A_\pi
\end{array}
\end{equation}
\paragraph{Déformation de la différentielle de Cartan-Eilenberg :}

 Il existe une graduation liée au nombre de sommets aériens dans
les graphes de Kontsevich. On note $\epsilon$ le paramètre de
graduation; cela revient à changer $\pi$ en $\epsilon \pi$ dans les
formules. \\

\begin{lem} Le terme $A_\pi$ est  une différentielle.
On a $A_\pi =\epsilon \, d_{CE} +o(\epsilon),$ où $d_{CE}$ désigne
 la différentielle de Cartan-Eilenberg de $\h$ agissant dans
 le $\h$-module $S(\g)/S(\g)\cdot \h=S(\g/\h)=\mathcal{C}_{poly}(\h^\perp)$.

\end{lem}

\noindent \textit{\bf Preuve :} Pour $n=1$ les seuls graphes qui
interviennent  sont comme dans (\ref{roue2.eps}) (les deux dessins de droite). Pour $\phi \in C_{poly}(\h^\perp)\otimes \bigwedge \h^*$ de degré
$q-1$ les graphes dessinés contribuent comme

\[\sum\limits_{i<j} (-1)^{i+j} \phi( [K_i, K_j], K_1, \ldots,
\hat{K_i}, \ldots, \hat{K_j}, \ldots, K_q)\]et
\[\sum\limits_i(-1)^{i+1} \ad K_i
\big(\phi( K_1, \ldots, \hat{K_i}, \ldots, K_q)\big).\]
On reconnaît la différentielle de Cartan-Eilenberg. \fin\\

\begin{defi}
On note $H^\bullet_\epsilon (\h^\perp)$ l'espace de cohomologie de
réduction formelle pour la différentielle $A_\pi$ agissant dans
$\mathcal{C}_{poly}(\h^\perp)\otimes \bigwedge \h^*$. L'espace
$H^0_\epsilon (\h^\perp)$ sera  appelé l'espace de réduction.
\end{defi}

Pour $\phi\in H^\bullet_\epsilon (\h^\perp)$ on a $\phi=\phi_o+\epsilon \phi_1 +
\epsilon^2 \phi_2 \ldots.$ Le premier terme $\phi_o$ est dans l'espace de
cohomologie relative
\[H^\bullet_{CE}\big(\h, \; S(\g/  \h)\big),\] car on a
$A_\pi=\epsilon\; d_{CE} +o(\epsilon)$. Toutefois, on ne peut pas affirmer que
$H^\bullet_\epsilon(\h^\perp)$ est une déformation de
$H^\bullet\big(\h, \; S(\g/  \h)\big),$ car
l'application $\phi\mapsto \phi_o$ n'est pas \textit{a priori} surjective.\\

En particulier pour $f\in H^{0}_\epsilon(\h^\perp)$ c'est-à-dire une
fonction $f=f_0+\epsilon f_1 + \epsilon^2 f_2\ldots $ dans
$\mathcal{C}_{poly}(\h^\perp)[[\epsilon]]$
  solution de  $A_\pi(f)=0$, on aura $f_0\in S(\g/\h)^\h$.
L'application $f\mapsto f_0$ n'est pas surjective \textit{a priori}
et on ne pourra pas considérer $H^{0}_\epsilon(\h^\perp)$ comme une
déformation de l'espace $S(\g/\h)^\h$. C'est
 toutefois  le cas pour
les paires  symétriques et plus généralement lorsque $\h$ admet un
supplémentaire
$\h$-invariant (\textit{cf.} Proposition~\ref{propreductionsym}).\\

\paragraph{Structure associative en cohomologie :}
La courbure étant nulle dans notre situation linéaire, l'équation
générale
\[[B_\pi, B_\pi]_{GH}+2[B_2, A_\pi]_{GH}=0\]montre que $B_\pi$ définit un produit
associatif en cohomologie. Cette formule résulte de la formule de
Stokes avec $3$ points sur l'axe réel, il faut tenir compte des
strates de bord qui contiennent tout le chemin : cette strate
fournit un co-bord pour la différentiel $A_\pi$. Donc ce n'est qu'en
cohomologie que ces formules fournissent un produit associatif. On a
donc la proposition suivante\footnote{Si $\pi$ est un bi-vecteur
quelconque cette proposition est encore vraie dès que l'on a
$F_\pi=0$ (\textit{cf.} \cite{CF1})}.

\begin{prop}\label{propproduitasso}
Dans le cas linéaire, $B_\pi$ définit un produit associatif dans
$H^\bullet_\epsilon(\h^\perp)$.
\end{prop}

\subsection{Dépendance par rapport au choix du
sup\-plé\-men\-taire}\label{dependance} On examine, dans le cas linéaire,  la dépendance du
produit associatif de la Proposition~\ref{propproduitasso}, par rapport au choix du supplémentaire.\\

\subsubsection{Action du
champ de vecteurs de déformation}\label{actionduchamp}

On fixe une base $(K_i)_i$ de $\h$. Fixons un supplémentaire $\q_o$
de $\h$. On fixe une base $(P_a)_a$ de $\q_o$, ce qui permet de
considérer $(P_a^*)_a$ base de
$\h^\perp$.\\

Faisons choix d'un autre supplémentaire $\q_1$ dont on fixe une base
$(Q_a)_a$, ce qui permet de considérer $(Q_a^*)_a$.
On suppose que l'on a $P_a^*=Q_a^*$, \textit{ie.}~on a $P_a-Q_a \in \h$.\\

Notons $\mu_{\q_o}$ et $\mu_{\q_1}$ les structures $A_\infty$
construites pour ces choix de supplé\-men\-taires. On notera dorénavant
$\mathcal{A}$ l'algèbre des fonctions polynomiales sur la
super-variété $\h^\perp\oplus \Pi \h$ :

\[\mathcal{A} := \mathcal{C}_{poly}(\h^\perp)\otimes
\bigwedge(\g^*/\h^\perp)\simeq \mathcal{C}_{poly}(\h^\perp)\otimes
\bigwedge \h^*.\]

Changer de supplémentaire  revient à faire un changement de
coordonnées linéaires sur $\g$. La matrice de changement de bases de
$(K_i, P_a)$ vers $(K_i, Q_a)$ a  la forme suivante :

\begin{equation}\nonumber
\mathbb{M}=\left(\begin{array}{cc}
 Id & \mathbb{D}\\
  0 & Id
\end{array}\right)
\end{equation}

Notons $\mathbb{D}=[V_1, \ldots, V_ p]$ les
colonnes de la matrice $\mathbb{D}$ avec  $V_i\in \h$.\\

Par transformée de Fourier partielle, $\widehat{\partial_{K_i^*}}$
correspond à la fonction sur $\Pi \h$ notée $\theta_{K_i^*}$,
$\widehat{K_i}$ correspond à dérivée $\partial_{\Pi K_i}$ sur $\Pi
\h$. Les
variables $P_a$ et $\partial_{P_a^*}$ ne sont pas changées.\\

La variété $\h^\perp \oplus \Pi \h$ est intrinsèque et ne dépend pas
du choix de la décomposition. Les fonctions $P_a$ et $Q_a$ sont
égales comme fonctions sur $\h^\perp$. Les fonctions
$\theta_{K_i^*}$ (pour les deux décompositions) sont
égales sur $\Pi \h$.\\

Les  bi-vecteurs de Poisson pour  les deux décompositions s'écrivent
respectivement\footnote{On convient, comme d'habitude,  que les
indices répétés sont sommés.} :
\begin{equation}
\begin{array}{ccc}
 \pi& = & [K_i, K_j]\partial_{K_i^*}\wedge\partial_{K_j^*}+ 2[K_i,
P_a]\partial_{K_i^*}\wedge\partial_{P_a^*}+[P_a,
P_b]\partial_{P_a^*}\wedge\partial_{P_b^*}\\
  \pi_1 & = & [K_i, K_j]\partial_{K_i^*}\wedge\partial_{K_j^*}+ 2[K_i,
Q_a]\partial_{K_i^*}\wedge\partial_{Q_a^*}+[Q_a,
Q_b]\partial_{Q_a^*}\wedge\partial_{Q_b^*}.
\end{array}
\end{equation}

Par  transformée de Fourier partielle notée $\mathcal{F}$ on a :

\begin{equation}
\begin{array}{ccc}
  \widehat{\pi} &= & \partial_{\Pi[K_i, K_j]} \theta_{K_i^*}\theta_{K_j^*} + 2
\mathcal{F}\big([K_i, P_a]\big)\theta_{K_i^*}\partial_{P_a^*} +
\mathcal{F}\big([P_a,
P_b]\big)\partial_{P_a^*}\wedge\partial_{P_b^*}\\
  \widehat{\pi_1}& = & \partial_{\Pi[K_i, K_j]} \theta_{K_i^*}\theta_{K_j^*} + 2
\mathcal{F}\big([K_i, Q_a]\big)\theta_{K_i^*}\partial_{Q_a^*} +
\mathcal{F}\big([Q_a,
Q_b]\big)\partial_{Q_a^*}\wedge\partial_{Q_b^*}.
\end{array}
\end{equation}
Comme on a
$\partial_{Q_a^*}=\partial_{P_a^*}$, le changement de supplémentaire
n'opère finalement que   sur les coefficients. Soit $(e_i)_i$ une base de $\g$ adaptée à la première décomposition et
$(e_i^*)_i$ sa base duale.  Considérons le bivecteur :
\begin{equation}\pi_\mathbb{M} :=\mathbb{M}^{-1}[\mathbb{M}e_i, \mathbb{M}e_j]\partial_{e_i^*}\wedge
\partial_{e_j^*}.\end{equation}
Alors on aura \[\widehat{\pi_\mathbb{M}}=\widehat{\pi_1}.\]

\begin{lem} Notons $v= -V_a
\partial_{P_a^*}$ le champ de vecteurs sur $\g^*$. On a $[v,v]_{SN}=0$
\footnote{On note $[\;, \;]_{SN}$ le crochet de
Schouten-Nijenhuis.}
 et
\[\pi_\mathbb{M}= e^{\ad v} \cdot \pi= \pi + [v, \pi]_{SN}+\frac {1}{2}[v,[v,\pi]]_{SN}.\]
L'action du champ $-v$ sur le bivecteur $\pi$ correspond au
changement de supplémentaire.
\end{lem}

\noindent \textit{\bf Preuve :}   Vérifions ce fait dans le cas des
paires symétriques (mais le résultat reste vrai en général). On a,
avec les notations introduites pour la matrice~$\mathbb{M}$, les équations suivantes :
\begin{eqnarray}\mathbb{M}^{-1}[\mathbb{M}K_i, \mathbb{M}K_j]=[K_i, K_j]\\\label{eqPIM1}
2\mathbb{M}^{-1}[\mathbb{M}K_i, \mathbb{M}P_a]=2\mathbb{M}^{-1}[K_i, P_a+ \underset{\in
\k}{\underbrace{\mathbb{D}P_a}}]=2[K_i,
P_a]-\underset{(1)}{\underbrace{2\mathbb{D}[K_i,
P_a]}}+\underset{(2)}{\underbrace{2[K_i, \mathbb{D}P_a]}}
\end{eqnarray}

\begin{multline}\label{equationPIM}
\mathbb{M}^{-1}[\mathbb{M}P_a, \mathbb{M}P_b]=
\mathbb{M}^{-1}[P_a+\mathbb{D}P_a,
P_b+\mathbb{D}P_b]=\\\underset{\in \k}{\underbrace{[P_a, P_b]}}+
\underset{(3)}{\underbrace{\underset{\in
\k}{\underbrace{[\mathbb{D}P_a,
\mathbb{D}P_b]}}}}+\underset{(4)}{\underbrace{\underset{\in\p}{\underbrace{\big([\mathbb{D}P_a,
P_b]+[P_a, \mathbb{D}P_b]\big)}} }}
-\underset{(5)}{\underbrace{{{\mathbb{D}\big(\underset{\in\p}{\underbrace{[\mathbb{D}P_a,
P_b]}}+[P_a, \mathbb{D}P_b]\big)}} }}
\end{multline}

\noindent Calculons $[v, \pi]_{SN}=v\bullet\pi+ \pi\bullet v$. On a
les formules suivantes :
\begin{equation}
\begin{array}{cccll}
 -v\bullet \pi &= & 2 V_a \langle P_a^*, [K_i, P_b]\rangle
\partial_{K_i^*}\wedge\partial_{P_b^*}&=&
2\underset{\in\k}{\underbrace{\mathbb{D}[K_i,
P_b]}}\partial_{K_i^*}\wedge\partial_{P_b^*} \\
-\pi\bullet v&= &2\underset{\in\k}{\underbrace{[V_a, K_i]}}\partial_{K_i^*}\wedge\partial_{P_a^*}+
2\underset{\in\p}{\underbrace{[V_a,
P_b]}}\partial_{P_b^*}\wedge\partial_{P_a^*}&&\\
v\bullet(v\bullet \pi)&= &0&&\\
v\bullet(\pi\bullet v )&= &2 V_c \langle P_c^*,[V_a,
P_b]\rangle \partial_{P_b^*}\wedge\partial_{P_a^*}&=&2
\underset{\in\k}{\underbrace{\mathbb{D}[V_a,
P_b]}}\partial_{P_b^*}\wedge\partial_{P_a^*}\\
(v\bullet \pi)\bullet v&= &2\underset{\in\k}{\underbrace{\mathbb{D}[K_i, P_b]}}\langle K_i^*, V_c\rangle
\partial_{P_b^*}\wedge\partial_{P_c^*}&=&
2\underset{\in\k}{\underbrace{\mathbb{D}[V_c,
P_b]}}\partial_{P_b^*}\wedge\partial_{P_c^*}\\
(\pi\bullet v)\bullet v&= &
2\underset{\in\k}{\underbrace{[V_a, K_i]}}\langle K_i^*,
V_c\rangle\partial_{P_a^*}\wedge\partial_{P_c^*}&=&2\underset{\in\k}{\underbrace{[V_a,
V_c]}}\partial_{P_a^*}\wedge\partial_{P_c^*}.
\end{array}
\end{equation}

On a donc
\[[v, \pi]_{SN}=\underset{(1)}{\underbrace{-2\underset{\in\k}{\underbrace{\mathbb{D}[K_i,
P_a]}}\partial_{K_i^*}\wedge\partial_{P_a^*}}}+\underset{(2)}{\underbrace{2\underset{\in\k}{\underbrace{[K_i,
\mathbb{D}P_a]}}\partial_{K_i^*}\wedge\partial_{P_a^*}}}+
\underset{(4)}{\underbrace{2\underset{\in\p}{\underbrace{[\mathbb{D}P_a,
P_b]}}\partial_{P_a^*}\wedge\partial_{P_b^*}}}\]

et

\[\frac 12[v, [v, \pi]]_{SN}=- \underset{(5)}{\underbrace{2
\underset{\in\k}{\underbrace{\mathbb{D}[\mathbb{D}P_a,
P_b]}}\partial_{P_a^*}\wedge\partial_{P_b^*}}}-
\underset{(3)}{\underbrace{
\underset{\in\k}{\underbrace{[\mathbb{D}P_a,
\mathbb{D}P_b]}}\partial_{P_a^*}\wedge\partial_{P_b^*}}}.\] On
retrouve bien les termes de $\pi_\mathbb{M}$ des équations
(\ref{eqPIM1}) et (\ref{equationPIM}). \fin

\subsubsection{Dépendance par rapport au supplémentaire }
Considérons maintenant la transformée de Fourier partielle. L'action
infinité\-simale est donnée alors par $[\widehat{\,v},
\widehat{\,\pi}]_{SN}$. L'expression
\[\widehat{\,\pi}+ [\widehat{\,v}, \widehat{\,\pi}]_{SN}+ \frac 12
{[\widehat{\,v},[\widehat{\,v},\widehat{\,\pi}]]}_{SN}\] représente
la transformée de Fourier du bi-vecteur $\pi$ pour le deuxième choix
de supplémentaire, exprimée dans les coordonnées
 intrinsèques de $\h^\perp \oplus \Pi \h$. \\

 La transformée de Fourier partielle $\widehat{\,v}$ du  champ
$v$, s'écrit $-\partial_{\Pi V_a}\partial_{P_a^*}$.  C'est donc un
$2$-vecteur à coefficients constant,  qui vérifie
aussi l'équation $[\widehat{\,v}, \widehat{\,v}]_{SN}=0$.\\

Appliquons   le $L_\infty$ quasi-isomorphisme de Kontsevich noté
$U$. On note  $m$ la multiplication dans
$\mathcal{C}_{poly}(\h^\perp)\otimes \bigwedge \h^*$. On définit
$\pi_t$ par\footnote{On oublie l'indice $SN$ s'il n'y  a pas de
confusion possible.}

\begin{equation}\pi_t := e^{t\,  \ad v}\cdot \pi =\pi + t[v, \pi] + \frac {t^2}2 [v,[v,\pi]]].\end{equation}
On a $\pi_{t=0}=\pi$ et $\pi_{t=1}=\pi_1$. On désigne par $\mu_t$
la structure $A_\infty$ correspondante :
\begin{equation} \mu_t :=U(
e^{\widehat{\,\pi_t}})=m+\sum_{n\geq 1} \frac {\epsilon^n}{n!}U_n
\left(\widehat{\,\pi_t}, \ldots, \widehat{\,\pi_t}\right).\end{equation}Il
vient en conséquence  :
\begin{multline}
\mu_{\q_1}=\mu_{t=1}=U(e^{\widehat{\,\pi}+ [\widehat{\,v},
\widehat{\,\pi}]+\frac
12 [\widehat{\,v}, [\widehat{\,v}, \widehat{\,\pi}]]})= \\
m+\sum_{n\geq 1} \frac {\epsilon^n}{n!}U_n \left(\widehat{\,\pi}+
[\widehat{\,v}, \widehat{\,\pi}]+\frac 12 [\widehat{\,v},
[\widehat{\,v}, \widehat{\,\pi}]], \ldots, \widehat{\,\pi}+
[\widehat{\,v}, \widehat{\,\pi}]+\frac 12 [\widehat{\,v},
[\widehat{\,v}, \widehat{\,\pi}]]\right).\end{multline} La question
est de savoir comment écrire $\mu_{t=1}$ en fonction de $\mu_{t=0}$. La dérivée $DU_{\widehat{\,\pi_t}}$ au point
$\widehat{\,\pi_t}$ est un morphisme de complexes.  En dérivant en  $t$, on obtient l'équation différentielle :

\begin{equation}\frac{\partial \mu_t}{\partial t}= \epsilon
\,DU_{\widehat{\,\pi_t}} \left({[\widehat{\,v},
\widehat{\,\pi_t}]}_{SN}\right)={[DU_{\widehat{\,\pi_t}}
\left(\widehat{\,v}\right), \mu_t]}_{GH}.\end{equation}

\begin{lem}

L'opérateur $DU_{\widehat{\,\pi}} \left(\widehat{\,v}\right)$ est de
degré $0$ (si l'on tient compte de tous les degrés impairs) et n'a
pas de composantes de degré $-1$ dans le complexe de Hochschild.
\end{lem}

Dans notre situation linéaire, c'est le même raisonnement que dans
le Lemme~\ref{courburenulle}. En effet pour des raisons de dimension
et de degré, la composante de degré $-1$ ne pourrait apparaître que
pour le graphe réduit à $v$. Sa contribution est nulle car $v$
est nul sur $\h^\perp$. \fin

\paragraph{Remarque 4: }Notons toutefois que l'opérateur
$DU_{\widehat{\,\pi}} \left(\widehat{\,v}\right)$ possède  \textit{a
priori} des composantes en tout degré
positif pour le complexe de Hochschild.\\

\noindent Notons $Y$ le champ sur la variété de Maurer-Cartan défini  au point $\nu=U(e^{\widehat{\pi}})$ par :
\[Y(\nu)= [DU_{\widehat{\,\pi}} \left(\widehat{\,v}\right),
\nu]_{GH}.
\]
C'est plus généralement un champ de vecteurs sur la variété formelle image
$U(T_{poly}(\mathcal{A})_1)$ (l'image par $U$ des poly-vecteurs de
degré total $1$ ), c'est à dire une co-dérivation de la cogèbre
image $U(S^{+}(T_{poly}(\mathcal{A})_1))$. \footnote{Si $V$ est un
espace gradué on note $S^{+}(V)$  les éléments non constants dans
$S(V[-1])$.}

\noindent L'équation précédente montre que la courbe intégrale du
champ $Y$ relie $\mu_{t=0}$ et $\mu_{t=1}$. On en déduit que les
deux structures $\mu_{t=0}$ et $\mu_{t=1}$ sont conjuguées par un
élément du groupe formel des difféomorphismes de la variété formelle
pointée $U(T_{poly}(\mathcal{A})_1)$ (le groupe à un paramètre qui
intègre l'action
du champ $Y$). \\

Toutefois (dans le cas non gradué) d'après \cite{Kont} \S~4.5.2, les
deux structures  $\mu_{t=0}$ et $\mu_{t=1}$ sont équivalentes par le
groupe de jauge (des star-produits  formels). En effet l'équation
d'évolution s'écrit aussi

\[\frac{\partial \mu_t}{\partial t}=
[DU_{\widehat{\,\pi_t}} \left(\widehat{\,v}\right), Q]_{\mu_t},\]
où $Q$ désigne le champ impair sur la variété formelle
correspondant à la structure d'algèbre de Lie différentielle
graduée. C'est bien l'équation de \cite{Kont} \S~4.5.2. Par
conséquent les structures sont bien conjuguées (comme série formelle
en $t$) par le groupe de jauge. L'élément du groupe de jauge se
construit de proche en proche par rapport au degré formel de la
variable\footnote{Voir plus loin la construction perturbative.}~$t$.
Cet élément dépend évidemment du point $\mu_{t=0}$, tandis que
l'élément du groupe des transformations formelles (qui est un
groupe plus gros) ne dépend que du champ $\widehat{v}$ \cite{Man}.\\

Pour avoir les idées claires, tout élément du  groupe de jauge (des
star-produits formels ou des $A_\infty$-structures) s'écrit sous la
forme
\[g=\exp\left(\epsilon \Delta_1 + \epsilon^2 \Delta_2 + \ldots\right),\] où les  $\Delta_i$ sont des éléments de degré (total) $0$  dans $D_{poly}(\mathcal{A})$. Le  produit se fait grâce à la
formule de Campbell-Hausdorff formelle pour l'algèbre de Lie
graduée de $D_{poly}(\mathcal{A})$. Pour $g=\exp(D)$ on a par définition \[g\cdot \mu=e^{\ad D} \mu.\]

\begin{prop}\label{propdependance}
Lorsque l'on change de supplémentaire, les espaces de réduc\-tion
sont isomorphes. L'isomorphisme (comme série formelle en $t$) est
donné par l'action d'un élément du groupe de jauge de  la forme
\[\exp\left(\sum_{n\geq 1} t^n E_n\right)\] o\`u les opérateurs formels (en $\epsilon$)
 $E_n$ sont de degré total $0$
(on tient compte de tous les degrés impairs). En particulier les
structures $A_\infty$ données par $\mu_{t=0}$ et $\mu_{t=1}$ sont
$A_\infty$ quasi-isomorphes (formellement en $t$).
\end{prop}
Le problème délicat est que la composition pour les opérateurs
polydifféren\-tiels n'est pas associative. Par ailleurs il faut
faire intervenir la notion d'éléments de type super-groupe dans les
cogèbres cocommutatives colibres (voir \cite{Man} et \cite{AMM}). Grâce à la bar- construction, on peut
composer les morphismes de cogèbres. Si $D=(D_i)_{i\geq 0}$ est un
élément plat (c'est-à-dire sans terme de degré $-1$) de
$D_{poly}(\mathcal{A})$ on considère $\overline{D}$ la co-dérivation
de
cogèbre (sans co-unité) dont les coefficients de Taylor sont les $D_i$.\\

Alors $e^{\overline{D}}$ a un sens et c'est un morphisme de
cogèbres. On vérifie que l'on a pour  $\overline{\mu}$ (la
co-dérivation de carré nul associée à $\mu$) :
\[e^{\overline{D}} \overline{\mu} e^{-\overline{D}}=\overline{e^{\ad D}
\mu}\]ce qui montre que dans notre situation les structures
$A_\infty$ sont bien $A_\infty$ quasi-isomorphes.

\begin{cor} Les
structures $A_\infty$ correspondant à deux choix de
supplémen\-taires sont $A_\infty$-quasi-isomorphes.
\end{cor}
\subsubsection{Entrelacement des espaces de
cohomologie}\label{soussectionentrelacement}

Fixons deux supplémentaires et notons dans cette sous-section
$\mu_1 $ et $\mu_2$ les structures $A_\infty$ correspondantes.\\

D'après la partie précédente il existe  une série formelle en $t$,  $D=t \,E_1+ t^2 E_2
\ldots$ de degré total $0$  tel que l'on ait $e^{\overline{D}} \overline{\mu} e^{-\overline{D}}=\overline{e^{\ad D}
\mu}$. On écrit $D=D_0+ D_1 +\ldots$ pour la décomposition dans le
complexe de Hochschild décalé. L'opérateur $D_0$ ne prend qu'un
seul argument, $D_1$ prend $2$ arguments mais est de degré total $0$
et par conséquent $D_1(f,g)=0$ si $f, g$ n'ont pas de composantes
dans $\bigwedge
\h^*$.\\

 On aura alors
\[e^{\ad D} \mu_1=\mu_2.\] Si on note
$\mu_1=\mu_1^{(0)}+\mu_1^{(1)}+ \ldots$ les composantes homogènes
(pour la graduation de Hochschild décalée) de $\mu_1$, on aura alors
\[[D, \mu_1]_{GH}=[D_0,\mu_1^{(0)}]_{GH}+ \ldots,\]par conséquent la
différentielle $\mu_1^{(0)}$ se transforme selon le champ de la
variété formelle $\xi \mapsto [D_0, \xi]_{GH}$. On en déduit que les
deux différentielles vérifient
\[e^{\ad D_0} \mu_1^{(0)}=\mu_2^{(0)}.\]Les
différentielles  $ \mu_1^{(0)}$ et $\mu_2^{(0)}$ sont donc
conjuguées par l'élément du groupe de jauge pour le complexe de
Hochschild
\[\phi: f \mapsto e^{D_0} f=f + D_0(f)+ \frac 12
D_0^2(f) +\ldots, \] c'est-à-dire que l'on a $\phi \circ
\mu_1^{(0)}=\mu_2^{(0)}\circ \phi.$ En particulier $\phi$ est un
isomorphisme de $H^{\bullet}\left(\mathcal{A}, \mu_1^{(0)}\right)$ sur $H^{\bullet}\left(\mathcal{A},
\mu_2^{(0)}\right)$.\\

 Calculons à l'ordre $1$ dans le complexe
de Hochschild (c'est-à-dire que l'on ne garde que les opérateurs au
plus $2$-différentiels): on utilise la formule à l'ordre $1$ en $Y$
\[e^{X+Y}=e^X\left(1+ \frac{1-e^{-\ad X}}{\ad X} Y\right).\]
On calcul à l'ordre $1$ :

\begin{multline}\nonumber
\mu_2^{(0)}+ \mu_2^{(1)}= e^{\ad D_0 + \ad D_1 }(\mu_1^{(0)}+
\mu_1^{(1)})=\\e^{\ad D_0}\left(1+ \ad\left(\frac{1-e^{-\ad
D_0}}{\ad D_0} D_1\right)
\right)(\mu_1^{(0)}+ \mu_1^{(1)})=\\
e^{\ad D_0} \mu_1^{(0)} + e^{\ad D_0}
\left(\mu_1^{(1)}+\ad\left(\frac{1-e^{-\ad D_0}}{\ad D_0} D_1\right)
\mu_1^{(0)}\right).
\end{multline}
On en déduit que l'on a

\begin{equation} \mu_2^{(1)}= e^{\ad D_0}
\left(\mu_1^{(1)}+\ad\left(\frac{1-e^{-\ad D_0}}{\ad D_0} D_1\right)
\mu_1^{(0)}\right)= \\
e^{\ad D_0} \mu_1^{(1)}+\underset{cobord\,\, pour \,\,
\mu_2^{(0)}}{\underbrace{[\frac{e^{\ad D_0}-1}{\ad D_0} D_1,
\mu_2^{(0)}]}}.\end{equation}

  En conclusion
les deux produits $e^{\ad D_0} \mu_1^{(2)}$ et $\mu_2^{(2)}$
définissent les mêmes produits dans l'espace de cohomologie défini
par $\mu_2^{(0)}$ et $\phi$ réalise un isomorphisme d'algèbres
de $H^{\bullet}\left(\mathcal{A}, \mu_1^{(0)}\right)$ muni du
produit $\mu_1^{(1)}$ sur $H^{\bullet}\left(\mathcal{A},
\mu_2^{(0)}\right)$ muni du produit $\mu_2^{(1)}$.\\

En particulier lorsque $f, g $ sont deux fonctions dans l'espace de
réduction
 (\textit{ie.} sans composantes
dans $\bigwedge \h^*$) on aura\footnote{Rappelons que l'on a
$\tilde{D}_1(f,g)=0$ lorsque le degré total de $\tilde{D}_1$ est
nul.} :
\begin{equation}\mu_2^{(1)}(e^{D_0}f,e^{D_0}g)=
\left(e^{\ad D_0}
\mu_1^{(1)}\right)(e^{D_0}f,e^{D_0}g)=e^{D_0}\left(\mu_1^{(1)}(f,
g)\right),\end{equation}ce qui montre que l'application $\phi= e^{D_0}$ réalise
un isomorphisme d'algèbres sur les espaces de réduction.\\

On résume cette section en énonçant le théorème suivant :

\begin{theo}\label{theodependance}
Lorsque l'on change de supplémentaire les espaces de réduction sont
isomorphes. L'isomorphisme est donné par l'exponentielle d'un
opérateur différentiel (série formelle en $t$) de degré $0$. C'est
aussi un isomorphisme d'algèbres de l'espace de cohomologie $\Big(
H^{\bullet}\left(\mathcal{A}, \mu_1^{(0)}\right),  \mu_1^{(1)}\Big)$
sur l'espace de cohomologie $\Big(H^{\bullet}\left(\mathcal{A},
\mu_2^{(0)}\right), \mu_2^{(1)}\Big)$.

\end{theo}

\subsection{Bi-quantification}\label{soussectionbiquantification}
C'est une idée essentielle de l'article de Cattaneo-Felder \cite{CF2}.\\

Dans leur  article, les auteurs définissent une fonction d'angle
dépendant de $4$ couleurs. En considérant deux sous-variétés
co-isotropes de $\g^*$, par exemple $\h_1^\perp\;$ et  $\; \h_2^\perp$, cette fonction d'angle à $4$ couleurs permet de
définir sur un espace de déformation modelé sur $(\h_1+ \h_2)^\perp=
\h_1^\perp \cap \h_2^\perp$ une action  de l'espace de réduction $
H^\bullet_\epsilon(\h_1^\perp)$ à droite et une action de l'espace
de réduction $H^\bullet_\epsilon(\h_2^\perp)$ à gauche.

\subsubsection{Définition de la fonction d'angle
à $4$ couleurs}\label{defi4couleurs}

Soit $\epsilon_1, \epsilon_2$ dans $\{-1, 1\}$. Les couples $(\epsilon_1, \epsilon_2)$ sont les
couleurs. Pour $p, q$ dans le premier quadrant\footnote{Les nombres complexes
vérifiant $Im(p)>0, Imp(q)>0, Re(p)>0$ et $ Re(q)>0$.}, on définit
la fonction d'angle à $4$ couleurs :

\begin{equation}\phi_{\epsilon_1, \epsilon_2}(p,q)=\arg (p-q)+\epsilon_1 \arg
(p-\overline{q}) + \epsilon_2 \arg (p+ \overline{q})
+\epsilon_1\epsilon_2 \arg(p+q).\end{equation}

\begin{figure}[h!]
\begin{center}
\includegraphics[width=8cm]{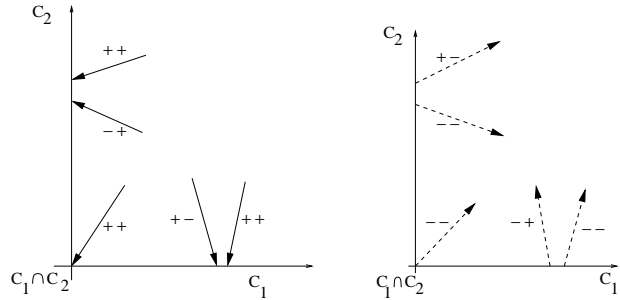}
\caption{\footnotesize Fonction d'angle à
$4$ couleurs}\label{4-couleurs2.eps}
\end{center}
\end{figure}

\noindent Cette fonction d'angle vérifie les propriétés de nullité
résumées dans Fig.~\ref{4-couleurs2.eps} o\`u on a représenté les
fonctions d'angle non nulles \textit{a priori}.\footnote{Bien que cela soit redondant par rapport à la couleur  $(\epsilon_1, \epsilon_2)$, dans  la figure ci-dessous les arêtes en trait plein "dérivent" au bord, tandis que les arêtes en pointillé "sortent" des bords. Cette convention est compatible avec les fonctions d'angle à 2 couleurs.} \\

On va considérer les intégrales des formes d'angle associés aux
graphes de Kontsevich colorés. Les  variétés, sur lesquelles on
intègre  ces formes, sont les variétés de configurations de
points dans le premier quadrant, modulo
l'action des dilatations. On utilise les mêmes compactifications que
celles décrites dans Kontsevich mais adaptées à notre situation.

\paragraph{Propriétés de bord :}

\begin{enumerate}
\item Lorsque $p ,\,q$ se concentrent sur l'axe horizontal,
 les formes d'angle $\mathrm{d} \phi_{\epsilon_1, \epsilon_2}(p,q)$
  tendent vers la $1$-forme d'angle \[\mathrm{d}\phi_{\epsilon_1}(p,q)= \mathrm{d}\arg (p-q)+\epsilon_1 \mathrm{d}\arg
(p-\overline{q}) \]

\item Lorsque $p ,\,q$ se concentrent sur l'axe vertical,
 les formes d'angle $\mathrm{d}\phi_{\epsilon_1, \epsilon_2}(p,q)$
  tendent vers la $1$-forme d'angle\[\mathrm{d}\phi_{\epsilon_2}(p,q)= \mathrm{d}\arg (p-q)+\epsilon_2 \mathrm{d}\arg
(p+\overline{q}). \]

\end{enumerate}

En conséquence  dans la formule de Stokes,  les variétés qui
apparaissent dans les concentrations près de l'axe horizontal et de
l'axe vertical,  correspondent aux variétés de configurations dans le
demi-plan de Poincaré avec formes d'angle à $2$ couleurs décrites
dans \S~\ref{defifonctionangle2couleurs}. On retrouve la situation
avec une seule
variété co-isotrope $C_1$ ou $C_2$.

\subsubsection{Graphes à $4$ couleurs}
On considère un graphe $\G$ numéroté et coloré par nos $4$ couleurs. Ce graphe est dessiné dans le premier quadrant. On place les points de première espèce à l'intérieur du quadrant et  les points de seconde espèce sur les axes ou à l'origine. Les arêtes issues  des points de seconde espèce (ainsi que les arêtes arrivant sur les points de seconde espèce) sont colorées comme dans Fig. \ref{4-couleurs2.eps}. Certaines arêtes, de couleur $(-,-)$, peuvent ne pas avoir de but, on dira qu'elles partent à l'infini.\\

Ces données permettent de définir un coefficient $w_\G$ comme dans \S~\ref{coefficient}.

\subsubsection{Définition
de la structure de bi-module}

Soit $C_1$ (par exemple $\h_1^\perp$) et $C_2$ (par exemple $\h_2^\perp$) deux sous-variétés co-isotropes. On fixe un
système de coordonnées et on suppose que l'on a localement une
situation d'intersection normale. On peut donc identifier  le fibré
normal  de $C_i$ à un voisinage de $C_i$ ($i=1,2$).\\

On considère un graphe de Kontsevich à $4$ couleurs. Aux sommets de première espèce
 on place un bivecteur
 de Poisson\footnote{Ici $\pi$ est un bi-vecteur de Poisson général. Dans le cas des algèbres de Lie, on considérera la moitié du crochet de Lie.}.\\

Aux sommets de l'axe horizontal, on place des poly-vecteurs restreints à
 $C_1$ avec dérivées normales par rapport à  $C_1$ (couleurs $(-, +)$ ou $(-,-)$). Dans notre cas il s'agit d'éléments de $\mathcal{A}_1:=\mathcal{C}_{poly}(\h_1^\perp)\otimes \bigwedge\h_1^*$.

Aux sommets de l'axe vertical,  on place des poly-vecteurs restreints à
 $C_2$ avec dérivées normales par rapport à  $C_2$ (couleurs $(+, -)$ ou $(-,-)$). Dans notre cas il s'agit d'éléments de $\mathcal{A}_2:=\mathcal{C}_{poly}(\h_2^\perp)\otimes \bigwedge\h_2^*$.

A l'origine,  on place un poly-vecteur restreint à $C_1\cap C_2$ avec
 dérivées normales à $C_1$ et $C_2$ (couleur $- \, -$). Dans notre cas il s'agit d'un élément  de $\mathcal{A}_{1,2}:=\mathcal{C}_{poly}(\h_1^\perp\cap\h_2^\perp)\otimes \bigwedge(\h_1^*\cap\h_2^*)$.\\

Comme en \S~\ref{grapheetoperateur} on utilise la règle complémentaire de dérivation suivante:  si la variable~est
\begin{enumerate}

\item dans $C_1$ et  dans $C_2$, la couleur sera $(\epsilon_1, \epsilon_2)=(+ \;, +)$
\item dans $C_1 $ mais pas dans $ C_2$, la couleur sera $(\epsilon_1, \epsilon_2)=(+ \;, -)$
\item dans $C_2$ mais pas dans $  C_1$, la
couleur sera~$(\epsilon_1, \epsilon_2)=(- \;, +)$
\item ni dans $C_1$ ni  dans
$C_2$, la couleur sera $(\epsilon_1, \epsilon_2)=(- \;, -)$.
\end{enumerate}

Après restriction à $C_1\cap C_2$, on définit un élément de $\mathcal{A}_{1,2}$ et un opérateur
$$B_\G :\mathcal{A}_2\otimes \mathcal{A}_{1,2}\otimes \mathcal{A}_1  \longrightarrow \mathcal{A}_{1,2}.$$

On considère alors  $\sum_{\G} w_\G B_\G$, où la somme porte sur  tous les graphes colorés à $4$ couleurs. C'est  une
opération tri-linéaire de
$\mathcal{A}_2\otimes \mathcal{A}_{1,2}\otimes \mathcal{A}_1$ dans $\mathcal{A}_{1,2}$\\

On note $\underset{1}\star$ le produit de droite et
$\underset{2}\star$ le produit de  gauche.

\subsubsection{Définition de l'espace de réduction  de deuxième
espèce de l'intersection $C_1\cap C_2$}
\paragraph{Hypothèse  1: } On fait l'hypothèse de  \textbf{courbure
nulle}, c'est à dire  $F_{\pi, C_i}=0$ pour $i=1$
et $i=2$
(procédure de quantification du  \S~\ref{quantificationpoisson}). \\

 On considère l'opérateur associé  aux contributions de tous les
graphes avec un sommet terrestre placé à l'origine. C'est un
opérateur qui agit sur les poly-vecteurs sur $ C_1\cap C_2$ (avec
dérivées transverse à $C_1\cup C_2$). Cet opérateur est encore,
grâce à la formule de Stokes et l'hypothèse $1$ ci-dessus, une
différentielle notée $\mathcal{A}_{\pi, C_1, C_2}$. L'espace de
réduction pour $C_1\cap C_2$, correspond à la cohomologie associée à
cet
opérateur.\\

\noindent \textbf{Hypothèse 2:} On suppose que la fonction constante
égale à~$1$ est dans l'espace de cohomologie, c'est-à-dire
$\mathcal{A}_{\pi, C_1, C_2}(1)=0$.\footnote{Cette hypothèse est
analogue à l'hypothèse sur la courbure pour les sous-variétés $C_i$.
En d'autres termes on suppose que la contribution des graphes
avec aucun point sur les axes ni à l'origine, est nulle. Elle assure que notre espace de réduction n'est pas nul.}\\

\paragraph{Remarque 5 : }Ces hypothèses sont automatiquement
vérifiées dans le cas des bi-vecteurs de Poisson linéaires.

\subsubsection{Compatibilité en cohomologie}

Expliquons par exemple, la compatibilité à  droite de l'action
$\underset{1}\star$ en cohomologie  : il suffit d'examiner les strates terrestres qui vont
intervenir dans les
différents  bords.\\

\noindent Soit $f$ dans l'espace de cohomologie pour la
différentielle $\mathcal{A}_{\pi, C_1, C_2}$ et soient $g$ et $h$
dans l'espace de cohomologie pour la différentielle
$\mathcal{A}_{\pi, C_1}$.\\

\noindent On place $f$ à l'origine, $g$ et $h$ sur l'axe horizontal.\\

\begin{figure}[h!]
\begin{center}
\includegraphics[width=7cm]{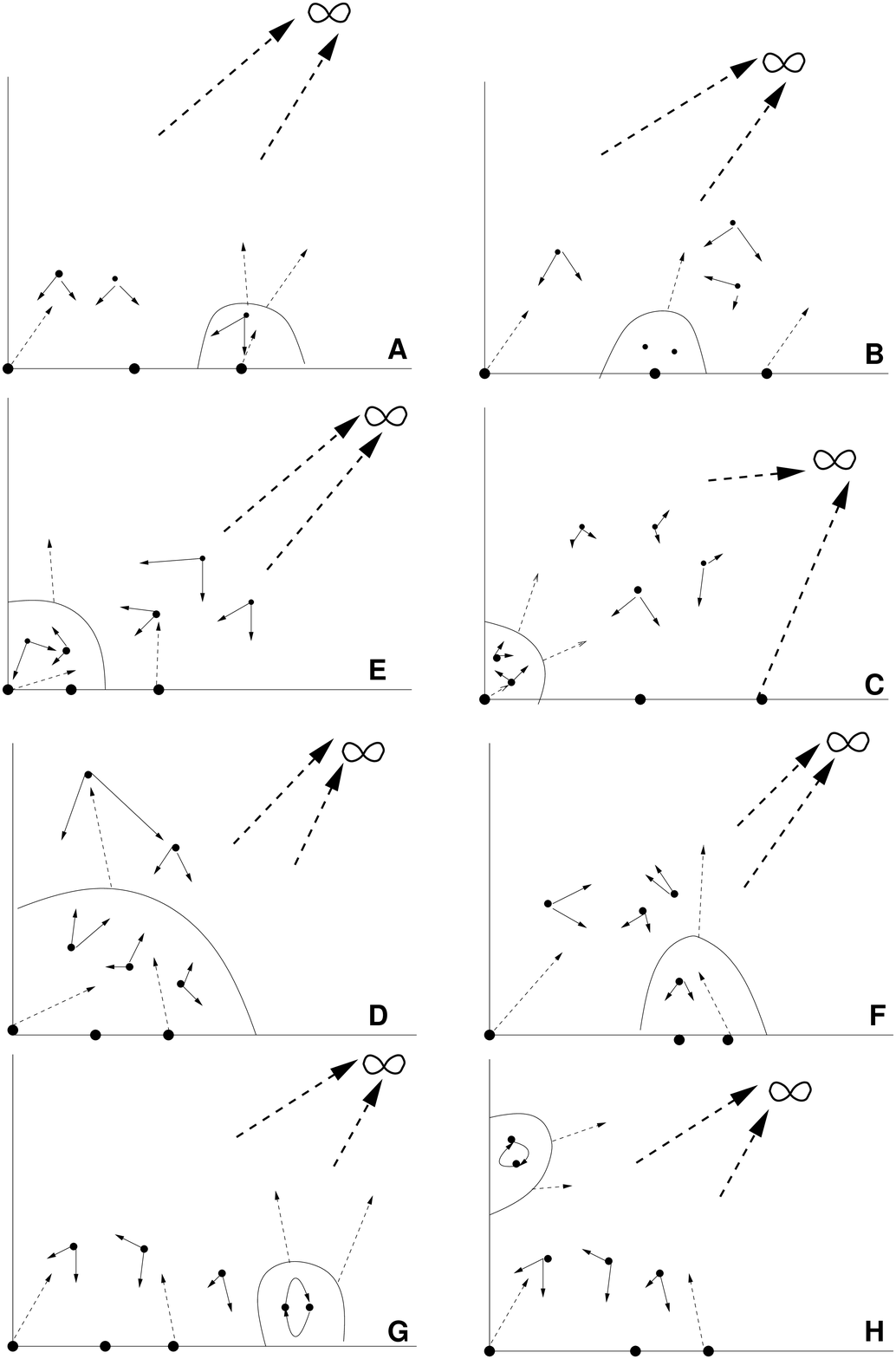}
\caption{\footnotesize Diverses strates pour la compatibilité de
l'action à droite en
cohomologie}\label{compatibilitecohomologie.eps}
\end{center}
\end{figure}
\noindent La formule de Stokes, appliquée dans ce contexte, va donner
une
compatibilité en cohomologie.\\

\begin{itemize}
\item
Les concentrations sur les axes sont nulles  (
Fig.~\ref{compatibilitecohomologie.eps} \textsc{cas} G,
Fig.~\ref{compatibilitecohomologie.eps} \textsc{cas} H) grâce à l'Hypothèse~1 de courbure nulle  ci-dessus.

 \item

Il est presque clair qu'avec ces conventions, les concentrations
horizontales ( Fig.~\ref{compatibilitecohomologie.eps}  \textsc{cas}
A , Fig.~\ref{compatibilitecohomologie.eps} \textsc{cas} B )
redonnent la différentielle $\mathcal{A}_{\pi, C_1}$ d'une seule
variété co-isotrope agissant sur $g$ ou $h$ .

\item Les concentrations (Fig.~\ref{compatibilitecohomologie.eps}  \textsc{cas} C )
 donnent
la différentielle $\mathcal{A}_{\pi, C_1, C_2}$ de la variété
$C_1\cap C_2$ agissant sur $f$.

\item La concentration  (Fig.~\ref{compatibilitecohomologie.eps}  \textsc{cas} E ) donne un terme de la forme
\[(f\underset{1}\star g)\underset{1}\star h\]

\item Les concentrations  (Fig.~\ref{compatibilitecohomologie.eps} \textsc{cas} F ) donnent un terme de la
forme \[f\underset{1}\star (g\underset{CF}\star h).\]

\item Il reste les concentrations
 (Fig.~\ref{compatibilitecohomologie.eps} \textsc{cas}  D ), qui produisent un cobord
pour $C_1\cap C_2$.
\end{itemize}

La compatibilité des actions à gauche et à droite résulte de la
formule de Stokes comme toujours.  Les strates de bord de type
Fig.~\ref{compatibilitecohomologie.eps} \textsc{cas} D fournissent
un cobord de l'espace de réduction de $C_1\cap C_2$, donc la
compatibilité des actions à droite  n'est vraie qu'en cohomologie.

\fin

 En considérant deux points terrestres (un à l'origine et
l'autre sur l'axe des abscisses) on en déduirait le lemme suivant.

\begin{lem}
L'action des cobords à droite (ou à gauche) sur les cocycles fournit
un cobord dans l'espace de $C_1\cap C_2$.
\end{lem}
On peut donc conclure cette section par la
compatibilité en cohomologie.
\begin{prop}

La description en termes de diagrammes de Feynman, fournit une
action en cohomologie  des espaces de réduction  de $C_1$ à droite
et de $C_2$ à gauche sur l'espace de cohomologie  de réduction de
\textbf{deuxième espèce} associé à $C_1\cap C_2$.

\end{prop}
En degré $0$ on retrouve le résultat de Cattaneo-Felder \cite{CF1}.
 \begin{cor}

 Au niveau des fonctions, on a une action à droite et à gauche
 des espaces de réductions de $C_1$ et de $C_2$ sur l'espace
 de réduction de deuxième espèce de $C_1\cap C_2$.

 \end{cor}

\paragraph{Remarque 6: } On fera attention au fait que la définition de  l'espace de réduction de
l'intersection, appelé de deuxième espèce, est différent de l'espace
de réduction de première espèce (pour les variétés co-isotropes).
D'ailleurs $C_1\cap C_2$ n'est pas une sous-variété co-isotrope,
 on ne peut donc pas appliquer la procédure de quantification du
\S~\ref{quantificationpoisson}.

\section{Exemples d'espaces de réduction}\label{ExempleReduction}
Prenons $\G$ un graphe coloré. Enon\c cons un lemme utile dans la
simplification des calculs.

\begin{lem}\label{couleurssuivent}
Si $\G$ a un sommet  de la forme $\overset{x}{\bullet} \leftarrow
\overset{z}{\bullet}\dashrightarrow \overset{y}{\bullet} $, alors
$w_\G=0$. \end{lem}

\noindent \textit{\bf Preuve :} Pour  un tel graphe, la forme
$\Omega_\Gamma$ (que l'on intègre) présente à ce sommet la
configuration suivante :

$$\overset{x}{\bullet} \leftarrow \overset{z}{\bullet}\dashrightarrow
\overset{y}{\bullet} \quad  = \quad \overset{x}{\bullet} \leftarrow
\overset{z}{\bullet}\leftarrow \overset{y}{\bullet} .$$

\noindent Les dérivées en $z$ dans $\Omega_\G$ ne proviennent que de
cette partie du diagramme, donc  le calcul se fait à $x, y$ fixés.
On a donc $w_\G=0$ d'après  Lemme~\S 7.3.3.1 de~\cite{Kont}.

\fin
\subsection[Cas des paires symétriques]{Cas des paires symétriques (ou cas d'un supplémentaire
stable)}\label{sectioncasdespairessymetriques}

On se place dans le cas des paires symétriques. On considère
 un caractère $\lambda $ de $\k$,  c'est-à-dire une forme linéaire sur $\k$
telle que $\lambda([\k, \k])=0$.  Alors $\lambda +\k^\perp$ est
l'espace affine des formes linéaires dont la restriction à $\k$ vaut
$\lambda$. C'est une sous-variété
co-isotrope.\\

\begin{prop}\label{propreductionsym} Dans le cas des paires symétriques (et plus généralement dans le cas où
$\h$ admet un supplémentaire stable sous l'action adjointe de $\h$)
on a $A_\pi= \epsilon \; d_{CE}$. Par conséquent l'espace de
réduction pour la sous-variété $\lambda+ \k^\perp$ est
$$H^0_\epsilon(\k^\perp)=\mathcal{C}_{poly}(\k^\perp)^{\k}[[\epsilon]]
=S(\p)^{\k}[[\epsilon]].$$

\end{prop}

\noindent \textit{\bf Preuve :} La démonstration se base sur
l'inventaire des graphes qui peuvent apparaître dans $A_\pi$.
Pour une paire symétrique (et plus généralement si $\h$ admet un
supplémentaire stable), les crochets vérifient $[\k, \p]\subset
\p$, ce qui imposent certaines restrictions sur les couleurs des
graphes. On place au point terrestre un élément de
$\mathcal{C}_{poly}(\lambda + \k^\perp) \otimes \bigwedge \k^*$ de
degré $q$. Il sort nécessairement $q+1$ arêtes à l'infini, donc au
moins une arête sort d'un point aérien du graphe $\G$. Pour
$n\geq 2$ les graphes qui
interviennent sont de deux types. \\

-i- Si dans le graphe aérien on a une arête colorée par $\k^*$ qui
sort à l'infini, alors au sommet de sortie on a dans le graphe $\G$:

-soit deux arêtes de même couleur qui se suivent  (car on a
$[\k, \p]\subset \p$ et $[\k, \k]\subset \k$)   et le coefficient
$w_\G$ vaut $0$ d'après un Lemme~\ref{couleurssuivent},

-soit une seule arête (colorée nécessairement par $\p^*$ car
$\lambda[\k,\k]=0$), auquel cas le coefficient $w_\G$ est encore nul
car $n\geq 2$ et ce sommet ne contribue qu'une fois\footnote{Ce
sommet se déplace sur $2$ dimensions, tandis que cette variable ne
contribue qu'une fois dans $\Omega_\G$.} dans la forme
$\Omega_\G$.\\

-ii- Si deux arêtes colorées par $\k^*$ sortent d'un même
sommet, alors ce sommet ne contribue qu'au plus une fois dans la
forme différentielle. Comme $n\geq 2$, le coefficient sera nul. \\

On conclut que les arêtes qui sortent à l'infini proviennent de
l'axe réel. Comme au moins une arête sort d'un sommet aérien,
toutes les contributions pour $n\geq 2$ sont nulles.\\

 En conséquence seuls les
graphes pour $n=1$ interviennent c'est-à-dire que l'on retrouve
comme espace de réduction, l'espace $\mathrm{d}_{CE} (\phi)=0$ pour
la différentielle de Cartan-Eilenberg.

\fin
\begin{cor}\label{cohomologiecassymetrique}Les constructions de Cattaneo-Felder
munissent $H^\bullet(\k, S(\p))$ d'une structure d'algèbre associative.
\end{cor}

\subsection{Cas linéaire pour les fonctions}\label{reductionlineaire}

Pla\c cons nous  dans le cas linéaire et prenons  $\h^\perp$ comme variété co-isotrope. On ne
suppose pas que $\h$ admet un supplémentaire stable.

\paragraph{Graphes intervenant dans le calcul de $A_\pi$ pour les
fonctions : }

Lorsqu'on place au sommet terrestre une fonction, les graphes qui
interviennent dans le calcul de la différentielle $A_\pi$ sont de
trois types\footnote{On pourrait regrouper les deux derniers types.}:\\

\begin{equation}\label{graphedifferentielle}\begin{array}{lcc}
               \includegraphics[width=5cm]{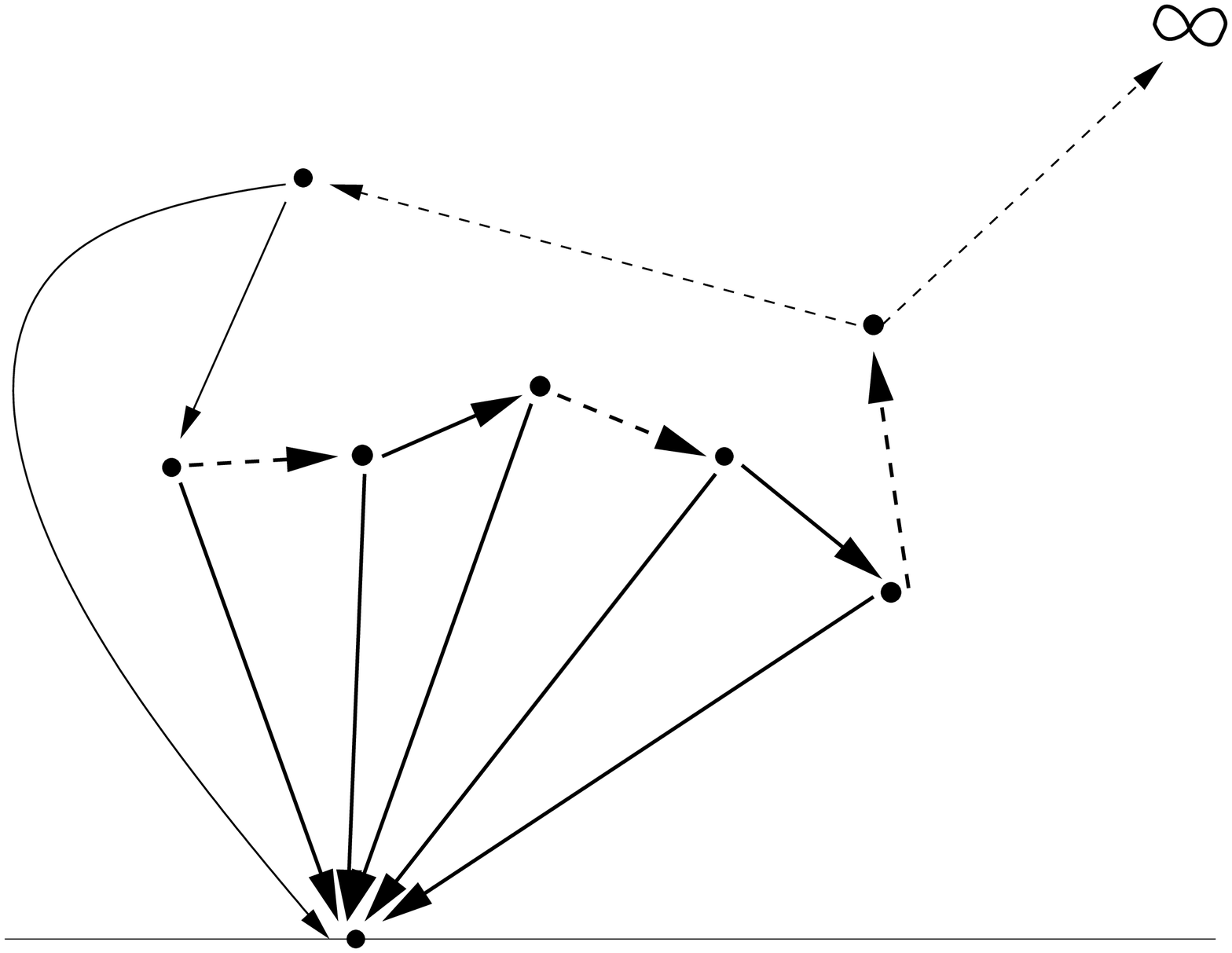}  & \includegraphics[width=5cm]{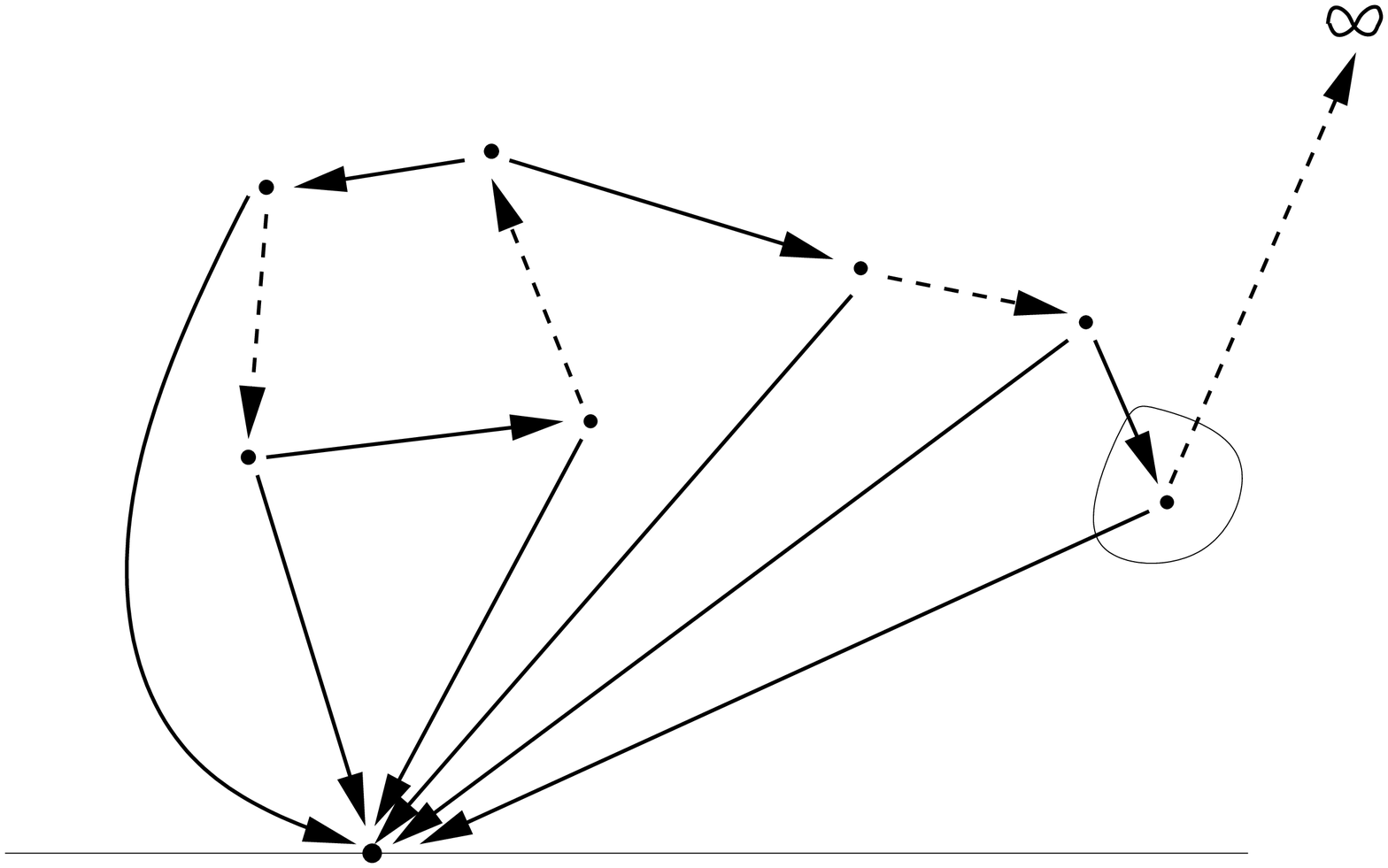} &  \includegraphics[width=5cm]{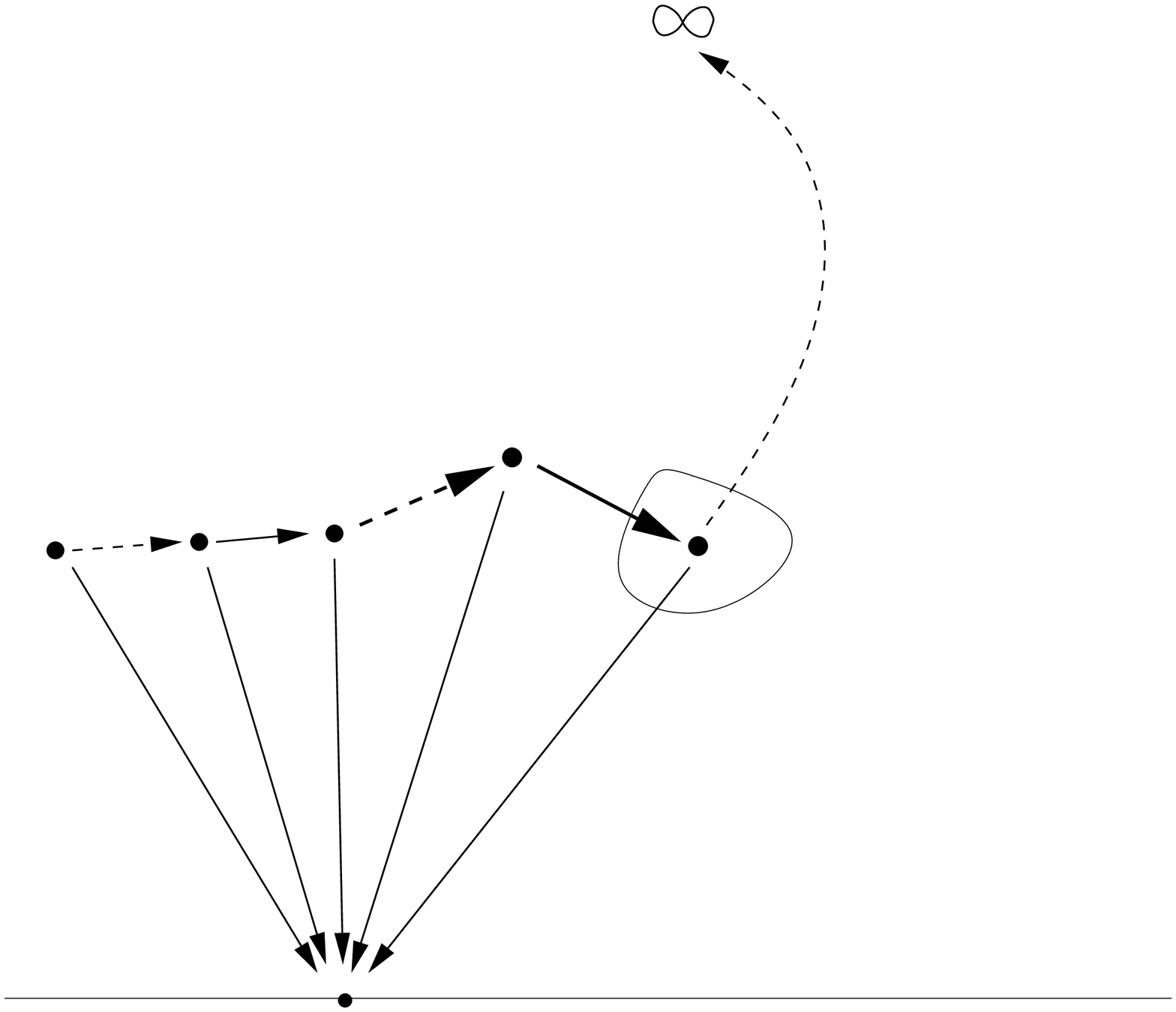}\\
                  Roue \; pure   & Roue-Bernoulli & Bernoulli
                \end{array}
\end{equation}

-i-  les graphes de type Bernoulli, avec la dernière arête partant
à
l'infini (\textit{cf.} (\ref{graphedifferentielle}) dessin de droite) \\

-ii- les graphes de type Roue, avec des rayons attachés directement
à l'axe réel sauf pour l'un d'entre eux qui est attaché à un
graphe de type Bernoulli dont la dernière arête part à l'infini
 (\textit{cf.} (\ref{graphedifferentielle}) dessin du milieu)\\

-iii- les  graphes de type Roue avec des rayons attachés
directement à l'axe réel sauf pour l'un d'entre eux qui part à
l'infini
 (\textit{cf.} (\ref{graphedifferentielle}) dessin de gauche)\\

\noindent  En effet ce sont les mêmes graphes que ceux que l'on
rencontre dans la formule de Kontsevich $f\underset{Kont} \star g$
pour les algèbres de Lie, lorsqu'on  prend pour $g$
une forme linéaire. La forme linéaire correspond ici à la direction qui part à l'infini.\\

\paragraph{Parité et homogénéité :}

\begin{lem}\label{lemmereductionlineaire} Les coefficients associés aux graphes intervenant dans le calcul de $A_\pi(f)$ sont nuls
si le nombre de sommets aériens est pair.
\end{lem}
\noindent \textit{\bf Preuve :} La variété de configurations admet
une symétrie par rapport à l'axe vertical. S'il y a $n$ sommets aériens, il
y a $2n-1$ arêtes (car l'une d'entre elles part à l'infini). En
faisant agir la symétrie, on trouve un facteur $(-1)^{n+1}$, ce qui
montre que $n$ doit être impair pour que le coefficient  soit non
nul. \fin

\begin{lem}\label{lemmeparite} Si $f=f_0+\epsilon f_1 + \epsilon^2 f_2 + \ldots $ est
dans l'espace de réduction, alors $F=f_0+ \epsilon^2 f_2 +
\epsilon^4 f_4 + \ldots $ l'est aussi et l'on peut prendre $F$
homogène si l'on tient compte  du degré de $\epsilon$ (qui vaut par
convention $1$).

\end{lem}

\noindent \textit{\bf Preuve :} Les graphes qui interviennent
dans les équations de réduction sont de trois types :

-i- Bernoulli $B_n$ : l'opérateur dérive $n$ fois mais ajoute un
degré (à cause de la racine), le degré total (y compris $\epsilon$)
est donc $1$.

-ii- Roue (attachée ou non à un Bernoulli ) : l'opérateur dérive
$n-1$ fois, le degré total  (y compris $\epsilon$) est donc $1$.

Les graphes avec un nombre pair de sommets n'interviennent pas,
donc les termes modulo $2$ vérifient encore le système de réduction.

Les équations s'écrivent (où les $D_n$ désignent les contributions des graphes
avec $n$ sommets) :

\[D_1(f_0)=0, \quad  D_1(f_2)+D_3(f_0)=0 \quad \ldots\]
donc les composantes homogènes (en tenant compte de $\epsilon$)
vérifient encore le système. \fin

\subsection{Cas linéaire pour les poly-vecteurs }

Décrivons les graphes qui vont intervenir dans le calcul de la
différentielle $A_\pi$ lorsqu'on l'applique sur un élément de $
\mathcal{A} = \mathcal{C}^\infty(\h^\perp)\otimes \bigwedge \h^*$, où $\h$ est encore une sous-algèbre quelconque de $\g$.

\begin{lem}Les graphes intervenant dans le calcul de la
différentielle des éléments de $\mathcal{A}$, sont de deux types:
les graphes déjà rencontrés pour les fonctions (la composante dans
$\bigwedge \h^*$ n'intervient pas) et les graphes qui sont le
crochet de deux Bernoulli, avec une dérivation unique de la racine
par une arête issue de l'axe réel (voir (\ref{differentielleforme}) dessin de droite).

\end{lem}

\noindent \textit{\bf Preuve :} Deux cas se présentent.\\

i- Toutes les arêtes issues de l'axe réel vont à l'infini; dans ce
cas les graphes qui interviennent sont ceux que l'on a détaillés pour les
fonctions.\\

ii- Une ou plusieurs arêtes, issues de l'axe réel, vont sur un sommet
aérien. Ces sommets sont des racines de graphes de type Lie (dont
certaines arêtes pourraient partir à l'infini).

Il n'y a qu'une arête issue de l'axe réel qui va vers un sommet
aérien sinon la forme à intégrer serait à variables séparées et on
aurait un problème de dimension (à cause du groupe de dilatation).

De ce graphe il part $2$ arêtes à l'infini. Raisonnons sur le
symbole $\G(X, Y)$  de ce graphe. Le mot $\G(X,Y)$ doit être de type
Lie et de degré deux en $Y$ (représentant les deux arêtes qui
partent à l'infini). Le mot s'écrit sous la forme $[A, B]$. Comme
l'arête issue de l'axe réel dérive le sommet du graphe, les mots $A$
et $B$ doivent contenir $Y$, sinon $A=X$ ou $B=X$ et on aurait une
arête double comme dans (\ref{differentielleforme}) (dessin de gauche\footnote{N'oublions pas que
la couleur inverse le sens de l'arête.}). En conséquence les mots $A$ et $B$
sont des mots de type Bernoulli c'est-à-dire de la forme
\[[\underset{n\,  \mathrm{fois} }{\underbrace{X,[X, \ldots , [X}}, Y]]],\]
l'arête partant à l'infini étant placée en dernière position. Ces
graphes sont comme dans (\ref{differentielleforme}) (dessin de droite).
 Lorsque les mots sont réduits à $Y$ on retrouve le graphe
 intervenant dans la différentielle de Hochschild.
 \fin

\begin{equation}\label{differentielleforme}\begin{array}{lr}
                  \includegraphics[width=6cm]{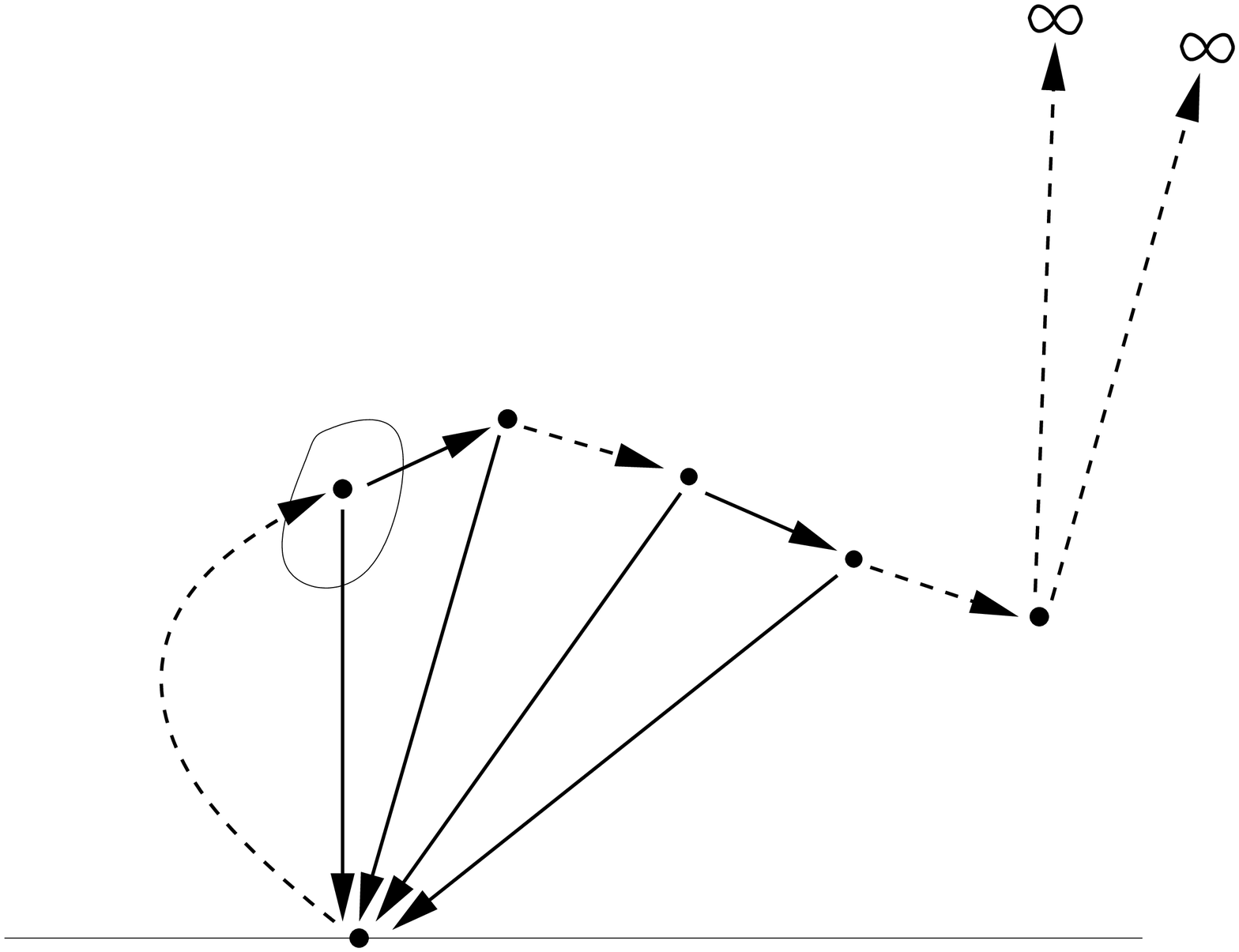} & \includegraphics[width=7cm]{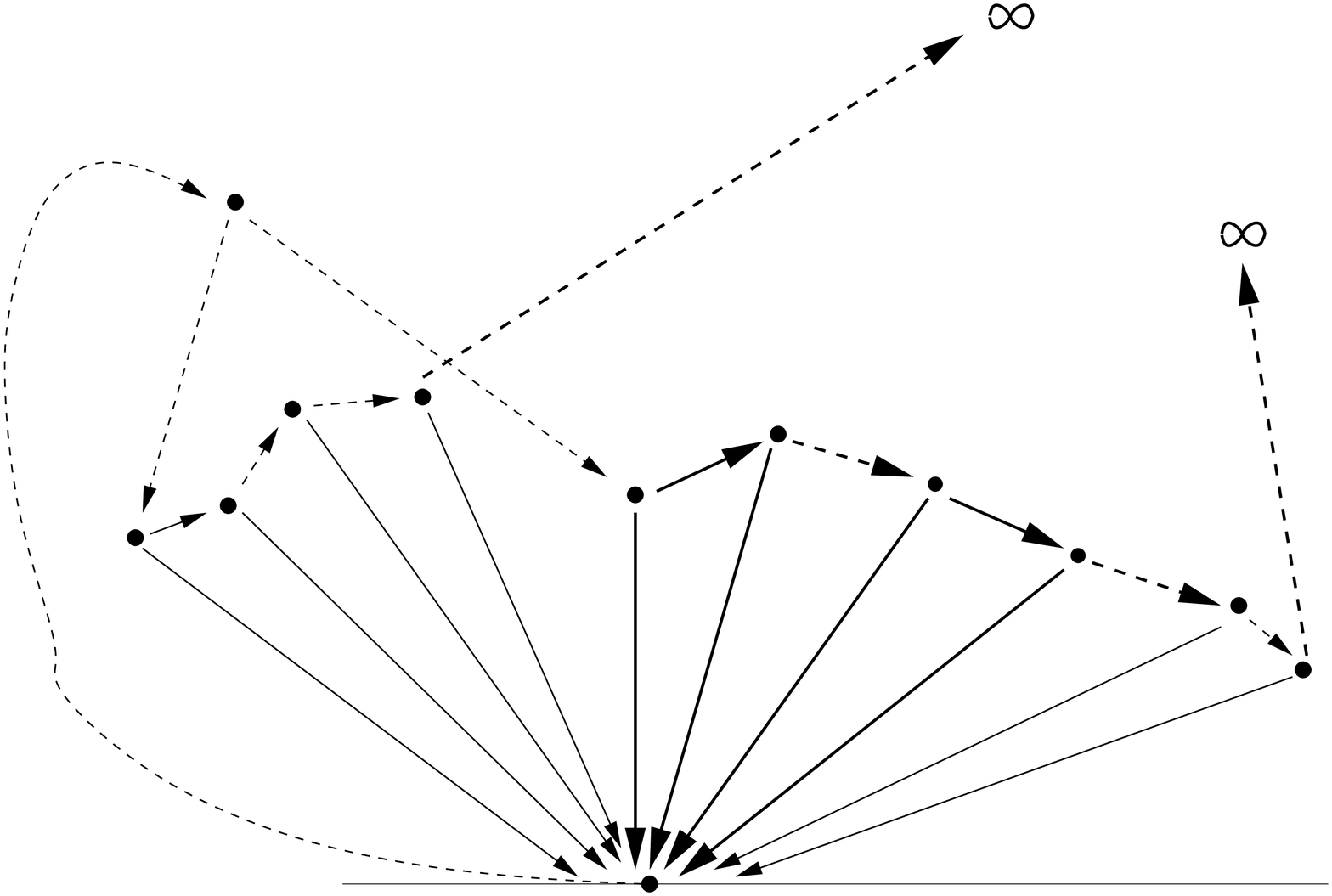}
                \end{array}
\end{equation}

\subsection{Cas Iwasawa}\label{sectioniwasawa}

Considérons une paire symétrique. Nous avons défini dans \cite{To3},
une notion de décomposition d'Iwasawa, généralisant la décom\-po\-si\-tion
d'Iwasawa pour les paires symétriques réductives.  Pour simplifier
la lecture de cet article, rappelons brièvement cette construction
qui intervient dans la définition de l'homomorphisme
d'Harish-Chandra généralisé.

\subsubsection[Décomposition d'Iwasawa généralisée]{Rappels  sur la notion de décomposition d'Iwasawa
généralisée }\label{iwasawa}

 Pour $\xi\in \k^\perp$
 on note $\g(\xi)$ le centralisateur de $\xi$ pour
l'action co-adjointe. C'est une sous-paire symétrique algébrique;
elle admet une décomposition de  Cartan notée
 $\k(\xi)\oplus \p(\xi)$.\\

Lorsque la dimension de $\g(\xi)$ est minimale pour  $\xi$ parcourant
$\k^\perp$, on dira que $\xi$ est régulier.  Les éléments réguliers
forment un ouvert de Zariski. On a alors $[\k(\xi), \p(\xi)]=0$, ce
qui permet de montrer que la sous-algèbre engendrée par $\p(\xi)$
est nilpotente. On note $\s_\xi$ sa partie semi-simple, c'est donc
une sous-algèbre abélienne formée d'éléments semi-simples (un tore)
et qui est dans $\p$. On dira que $\xi$ est générique lorsque la
dimension de ce tore est maximale. Ce tore généralise la notion de
\textit{\textbf{sous-espace de
Cartan}}.\\

On peut alors décomposer $\g$ sous l'action adjointe de $\s_\xi$
lorsque ce tore est diagonalisable (déployé\footnote{C'est le cas
si le corps de base est $\mathbb{C}$.}). C'est ce que nous
supposerons\footnote{Cette restriction n'est pas très sévère, car il
suffit de tout complexifier et de considérer la paire symétrique
complexe obtenue comme une paire symétrique réelle.}.  On note
$\g_o$ le centralisateur de $\s_\xi$. C'est une sous-paire
symétrique, appelée \textbf{\textit{petite paire symétrique.}}
On note $\g_o=\k_o\oplus \p_o$ sa décomposition de Cartan. Comme
$\s_\xi\subset \p$, si $\alpha$ est une racine alors $-\alpha $
aussi;  ceci permet d'exhiber, en séparant les racines et leurs
opposées, des décompositions triangulaires : $\g=\n_-\oplus \g_o\oplus \n_+$  et des décompositions d'Iwasawa $\g=\k\oplus \p_o \oplus \n_+$.

\subsubsection{Espace de réduction}\label{sectionespacereductioniwasawa} On fixe une décomposition
d'Iwasawa généralisée dans le cas des paires symétriques
\[\g=\k\oplus \p_o \oplus \n_+,\] et on note $\g_o=\k_o\oplus \p_o$.

\begin{prop}\label{propreductioniwasawa}

L'espace de réduction pour la sous-variété co-isotrope $(\k_o\oplus
\n_+)^\perp$ est
$$H^0_\epsilon((\k_o\oplus
\n^+)^\perp)=\mathcal{C}_{poly}(\k_o^\perp)^{\k_o}[[\epsilon]]
=S(\p_o)^{\k_o}[[\epsilon]].$$
\end{prop}

\noindent \textit{\bf Preuve :} On fixe un supplémentaire de $\k_o$
dans $\k$ noté $\mathfrak{r}$.
Les couleurs des arêtes sont donc :\\

-- les arêtes colorées par $ \dashrightarrow$ correspondent aux
variables dans $\n_+^*, \k_o^*$\\

-- les arêtes colorées par $\longrightarrow$ correspondent aux
variables $\p_o^*,
\mathfrak{r}^*$ que l'on peut identifier à $\k^*/\k_o^*$.\\

Soit $f$ vérifiant l'équation $A_\pi(f)=0$. D'après le Lemme~\ref{lemmeparite} on peut écrire
$f=f_0+\epsilon f_2+ \epsilon^4 f_4\ldots$ avec $D_1(f_0)=0$.  Donc $f_0$ est $\k_0\oplus  \n_+$-invariant. D'après \cite{To3} un
élément de $S(\p_0)^{\k_0}$.\\

Les graphes qui interviennent dans $D_n$ sont  de type Bernoulli ou de type Roue attachée à un Bernoulli.

 Examinons l'action de ces  opérateurs  sur $f_0 \in
S(\p_0)^{\k_0}$. Comme on ne peut dériver $f_0$ qu'en les directions
$\p_o^*$, toutes les arêtes qui arrivent sur le sommet terrestre sont
colorées par $\longrightarrow$. L'arête sortant à l'infini est soit
dans $\k_o^*$, soit dans
$\n_+^*$.\\

\begin{equation}\label{figureiwasawa} \begin{array}{lr}
                  \includegraphics[width=6cm]{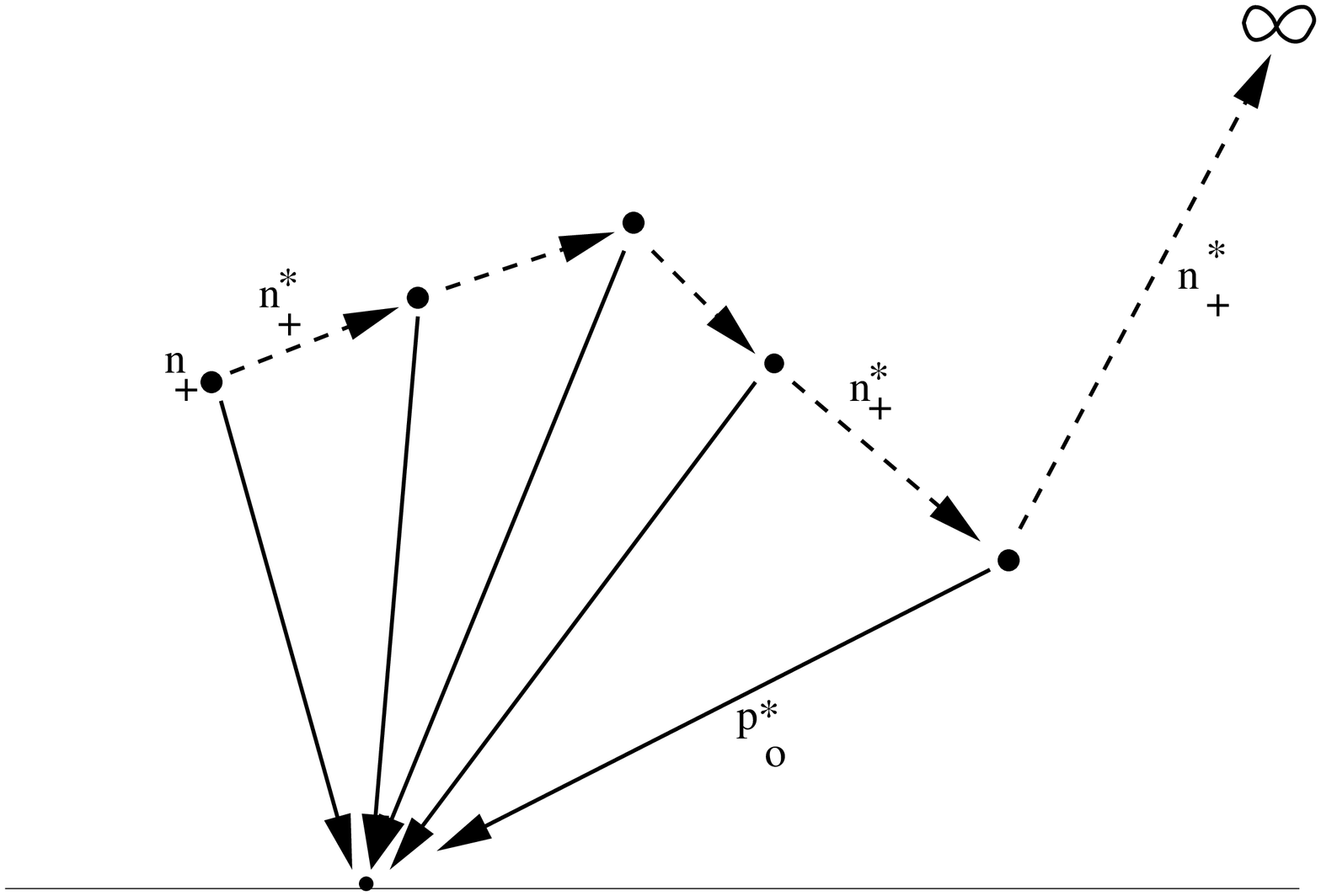} & \includegraphics[width=6cm]{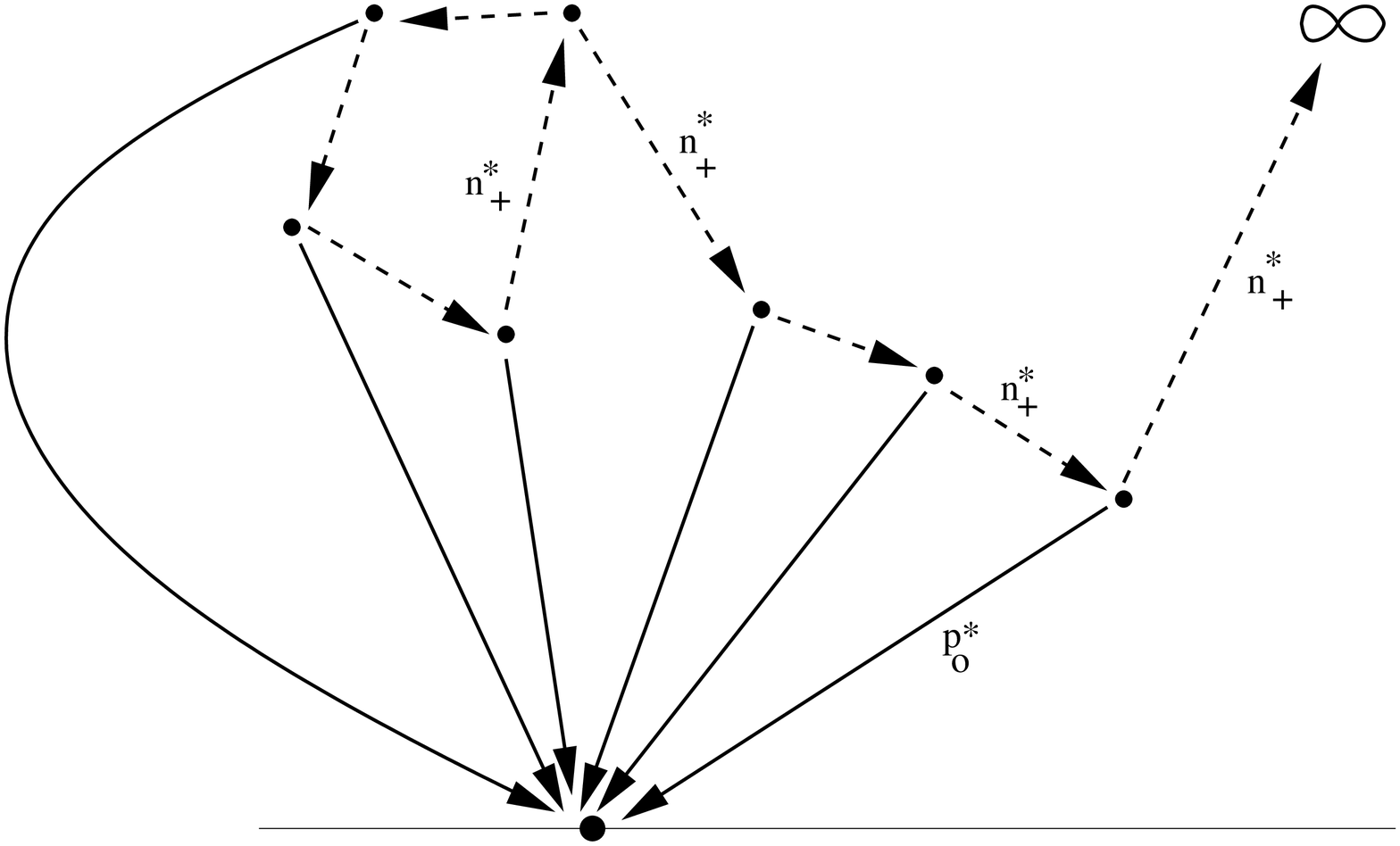} \\
                  Cas \; 1  & Cas \; 2
                \end{array}
\end{equation}
\paragraph{Premier cas :} C'est le type Bernoulli comme (\ref{figureiwasawa}) (dessin de gauche).\\

Le coefficient sera nul, si on a deux arêtes de type
$\longrightarrow$ qui se suivent au sommet o\`u l'arête part à
l'infini. Par conséquent l'arête qui part à l'infini est colorée par
$\n_+^*$. Par ailleurs on a  $[\p_0,
\n_+]\subset \n_+$, donc  la couleur $\dashrightarrow$ se propage dans le graphe. La racine du graphe représente  ,pour
l'opérateur différentiel, un coefficient dans $\n_+$ qui est donc nul par restriction à $(\k_o\oplus \n_+)^\perp$.\\

\paragraph{Deuxième cas :} Comme (\ref{figureiwasawa}) (dessin de droite),  le graphe à considérer de
manière générale est une roue attachée à l'axe réel et à un
sous-graphe de Bernoulli (ce graphe pouvant être vide auquel cas
l'arête part directement à l'infini et est colorée par
$\n^*_+$). Comme précédemment, pour que le coefficient soit non nul,
l'arête qui part à l'infini est colorée par $\n_+^*$. Ceci implique
que l'attache sur la roue se fait dans des dérivations selon
$\n_+^*$ (voir (\ref{figureiwasawa}) dessin de droite). La dérivation dans la roue
se fait alors soit de manière homogène selon les dérivées $\n_+^*$
pour la couleur $\dashrightarrow$, soit selon des dérivées dans la
direction $\k^*/\k_o^*$ avec couleur $\dashrightarrow$. Dans les
deux cas l'opérateur fera apparaître un terme de la forme

\[\tr_{\n_+}( \ad P_1 \ldots \ad P_l \ad N_1)=0\]ou un terme de
la forme

\[\tr_{\g/(\g_o\oplus \n_+)}( \ad P_1 \ldots  \ad P_l \ad N_1)=0.\]

L'opérateur différentiel sera alors nul quand on l'applique à des
fonctions qui sont dans $S(\p_o)^{\k_o}$. Par conséquent le système
d'équations est tout simplement $D_1(f_p)=0$. \fin

\subsection{Cas des polarisations}\label{sectioncasdespolarisations}

On prend $\xi\in \g^*$ une forme linéaire.  La forme bilinéaire
alternée $B_\xi$,  définie par  $B_\xi(X,Y)=\xi[X,Y]$ pour $X, Y \in
\g$, admet pour noyau  $\g(\xi)$. Une polarisation en $\xi$ dans
$\g$ est une sous-algèbre isotrope pour $B_\xi$ et de dimension
maximale. C'est automatiquement une sous-algèbre algébrique. On dit
que $\b$  vérifie la condition de Pukanszky si on a $\b=\g(\xi) +
\b_u$, où  $\b_u$ désigne le
radical unipotent de $\b$.   \\

Alors $\xi+\b^\perp$ est une sous-variété co-isotrope. On note $B$
un sous-groupe algébrique connexe de $G$ d'algèbre de Lie $\b$.

\begin{prop}\label{propreductionpolarisation}
L'espace de réduction pour la sous-variété $\xi+\b^\perp$ est le
corps $\R$, c'est-à-dire que l'on a
$$H^0_\epsilon(\xi+\b^\perp)=\R.$$

\end{prop}

\noindent \textit{\bf Preuve :} Considérons une fonction $f=f_0+
\epsilon^2 f_2+\epsilon^4 f_4 \ldots$ dans l'espace de réduction. La
première équation fournit des fonctions $B$-invariantes sur
$\xi+\b^\perp$, c'est-à-dire des fonctions constantes, car $B\cdot
f$ est ouvert dans $\xi+\b^\perp$. Comme précédemment les autres
équations se résument par récurrence aux
équations $D_1(f_n)=0$ pour tout $n$.
\fin

\section{Fonction $E(X,Y)$ pour les paires
symé\-tri\-ques}\label{sectionfonctionE}

On se place dans le cas des paires symétriques $\g=\k\oplus \p$ et
on considère le cas de la sous-variété $\lambda +\k^\perp$ où
$\lambda$ désigne un caractère de $\k$ (\textit{ie.} $\lambda \in (\k/[\k,\k])^*$). On
considère le produit $B_\pi$ (qui prend deux arguments) introduit au
\S~\ref{quantificationpoisson} pour $\pi$ le bi-vecteur de Poisson
associé à la moitié du  crochet de Lie.

\subsection{Définition de la fonction $E(X,Y)$ pour les paires
symé\-tri\-ques}\label{defiE}

\noindent  Pour $X, Y\in \p$ on considère les fonctions
exponentielles $e^X$ et $e^Y$. On considère alors $B_\pi(e^X, e^Y)$
que l'on restreint à $\lambda +\k^\perp$. On pourra consulter~\cite{AST} pour une description précise des graphes dans le cas
linéaire.

\begin{lem} Dans le calcul de $B_\pi(e^X, e^Y)$ les graphes de type
Lie avec racine dans $\p$ n'interviennent pas.
\end{lem}
\noindent \textit{\bf Preuve :} D'après le Lemme~\ref{couleurssuivent}, si la racine d'un tel graphe est dans $\p$,
alors le  poids associé est nul car le graphe présente, à la racine,
un sommet de la forme $\overset{x}{\bullet} \leftarrow
\overset{z}{\bullet}\dashrightarrow \overset{y}{\bullet} $.
\fin

Dans le cas des paires symétriques, seuls  les graphes avec racine
 dans $\k$ ou de type Roue interviennent, donc par
restriction on trouve l'action du caractère. Comme d'habitude, les
graphes sont des superpositions de graphes simples de type Lie ou Roue. Le résultat
est donc une exponentielle. \\

La composante logarithmique  de $B_\pi(e^X, e^Y)$  dans $\p$ est
donc $X+Y$ et celle dans $\k$ s'écrit $H(X, Y)$, avant restriction à
$\lambda + \k^\perp$.\\

Lorsque $\lambda=0$, seuls les graphes de type Roue vont intervenir de manière
significative.

\begin{defi}
 On note $E(X, Y)$ les contributions de tous les graphes de
type Roue. C'est une fonction analytique en $X, Y$ près de $0$.

\end{defi}

 On a donc  : \begin{equation}B_\pi(e^X, e^Y)=E_\lambda(X, Y)e^{X+Y},\end{equation}
où on a noté  $E_\lambda(X, Y) : =e^{\lambda(H(X, Y))}E(X, Y)$.\\

 La fonction $E_\lambda(X, Y)$ est le symbole formel du
star-produit de Cattaneo-Felder associé à la sous-variété
co-isotrope $\lambda + \k^\perp$.

\begin{lem}\label{symetrieE}

On a les symétries suivantes $E_\lambda(X, Y)=E_{-\lambda}(Y, X)$ et $E_\lambda(-X, -Y)=E_\lambda(X, Y)$
\end{lem}

\noindent \textit{\bf Preuve : } Dans le cycle d'une roue, lorsqu'une sortie est colorée par $\p^*$
il se produit un changement de couleur dans le cycle. En définitive,
le nombre de changements de couleur est pair, c'est-à-dire qu'un
tel cycle est attaché à un nombre pair de sorties dans $\p$, disons
$2p$ et à un nombre $q$ de sorties dans $\k$. Par ailleurs un
sous-graphe (de type Lie) dont la racine est dans $\p$ a un nombre pair
de sommets, et chaque graphe dont la  racine est dans $\k$ a un nombre
impair de sommets (car les "pattes" de ces graphes sont attachées
sur $e^X, e^Y$, c'est-à-dire que les pattes sont colorées par
$\p^*$). Il y a donc $2p+q$ sommets dans le cycle et $q \pmod 2$
sommets en dehors du cycle. Au total on a un nombre pair de
sommets dans
$\G$ et donc $E(-X, -Y)=E(X, Y)$.\\

Faisons une symétrie axiale. On a
\begin{equation}\label{symcoeff}
\G^{\vee}(X, Y)=\G(Y,X) \quad \mathrm{ et } \quad  w_{\G^{\vee}}=(-1)^{\sharp\{\mathrm{E_\G}\}} w_{\G}= w_\G,
\end{equation}ce qui
montre la symétrie recherchée.

Pour les graphes de type Lie, qui contribuent pour la composante
dans $\k$, on obtient $w_{\G^{\vee}}=- w_\G$, car on a un nombre
impair de sommets aériens. Les mots de Lie qui sont dans $\k$ sont nécessairement pairs.

\fin\\

\begin{defi}

On définit le produit de Cattaneo-Felder $\underset{CF,\,
\lambda}\star$ sur $S(\p)$ par la relation pour $f, g \in
S(\p)=\mathcal{C}_{poly}(\k^\perp)$ et $\xi \in \k^\perp=\p^*$:
\[(f\underset{CF,\, \lambda}\star g)(\xi)=B_\pi(f,g)=E_\lambda(\partial_\nu,
\partial_\eta)(f\otimes g)|_{\nu=\eta=\xi}.\]
\end{defi}

Pour $\lambda=0$, on notera ce produit $\underset{CF}\star$. D'après la Proposition~\ref{propreductionsym} et (\ref{symcoeff}) on déduit le lemme suivant.

\begin{lem}
Le produit $\underset{CF,\, \lambda}\star$ est
associatif sur $S(\p)^\k$ et vérifie  $f\underset{CF,\,
\lambda}\star g=g\underset{CF,-\lambda}\star f$.
\end{lem}
Le corollaire suivant précise le
Corollaire~\ref{cohomologiecassymetrique}
 (\S\ref{sectioncasdespairessymetriques}) et décrit ce qui se passe
en cohomologie.
\begin{cor}
 Le produit associatif de Cattaneo-Felder
 dans $H^\bullet(\k, S(\p))$ est  donné par l'action de
 l'opérateur $E$ sur les coefficients.  Le
 produit en cohomologie est donc l'extension du produit sur $S(\p)^\k$.
 \end{cor}

\noindent \textit{\textbf{Preuve :}}  Si une arête issue de l'axe
réel intervient dans un
 graphe, il existerait une arête colorée par
 $\k^*$ qui sortirait à l'infini. Or au sommet d'où part une telle arête, on
 aurait  deux arêtes de même couleur qui se suivent et donc le
 coefficient serait nul par le  Lemme~\ref{couleurssuivent}.

 \fin
\subsection{Contributions dans
$E(X,Y)$}\label{sectioncontributionsdansE}

Dans cette section, nous exploitons une symétrie pour les graphes qui interviennent
dans $E(X, Y)$. Pour $X, Y \in \p$ on note $x=\ad X$ et  $y=\ad Y$
les opérateurs adjoints.

\begin{prop}\label{alternance} Dans le  logarithme de $E(X,Y)$,  interviennent des termes de la forme
$$w_\G\Big( \tr_\p(x_1\cdots x_{n})+ (-1)^{n-1}\tr_\k( x_{n}\cdots
x_2 x_1)\Big),$$

\noindent où $w_\G$ est le coefficient associé à un graphe coloré de
type Roue, $x_i=\ad X_i$ et $X_i$ est un mot de type Lie en
$X,Y$.
\end{prop}

\noindent \textit{\textbf{Preuve :}} On change l'orientation et la coloration des arêtes dans le cycle pour conserver, au signe près, la forme $\Omega_\G$.\\

\begin{figure}[h!]
\begin{center}
\includegraphics[width=10cm]{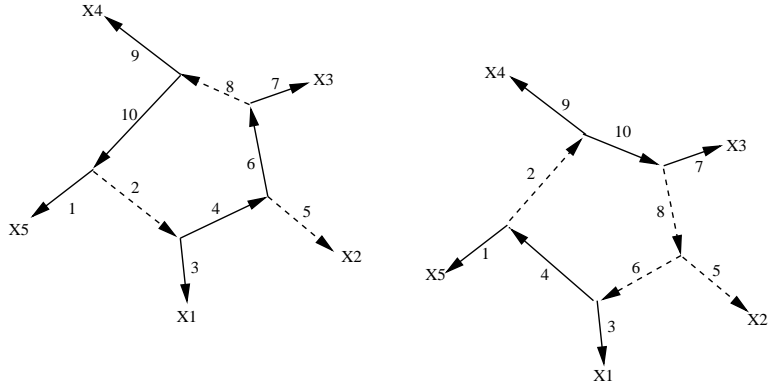}
\caption{\footnotesize Roue dans un sens et dans
l'autre}\label{doublesens.eps}
\end{center}
\end{figure}

Par exemple sur Fig.~\ref{doublesens.eps} les arêtes $2, 4, 6, 8,
10$ se retrouvent dans la forme $\Omega_\G$, en
position $10, 8, 6, 4, 2$. Comme on a inversé le sens des arêtes, la
forme est la même au signe près, donné par la signature de l'inversion de
numérotation dans le cycle; sur l'exemple Fig.~\ref{doublesens.eps}
le signe est donc  $(-1)^4=1$. On trouve dans $\ln(E(X,Y))$, une contribution de la forme
\[w_\G\Big(\tr_\p(x_5x_1x_2x_3x_4)+\tr_\p(x_5x_4x_3x_2x_1)\Big).\]

\noindent En général si  l'un des $X_i$ est dans $\p$, par exemple $X_1$ recevant la
couleur $\p^*$ dans la roue, alors on trouve des contributions de la
forme :
\[\tr_\p(x_1\cdots x_{n})\] et
\[(-1)^{n-1}\tr_\p(x_1\cdot x_{n}\cdot x_{n-1}\cdots x_2)=
(-1)^{n-1}\tr_\k( x_{n}\cdots x_2 x_1).\]Au total on a une contribution de la forme :
\[\tr_\p(x_1\cdots x_{n})+ (-1)^{n-1}\tr_\p(x_1 x_{n}\cdots x_2)=
\tr_\p(x_1\cdots x_{n})+ (-1)^{n-1}\tr_\k( x_{n}\cdots x_2 x_1).\]
Si tous les $X_i $ sont dans $\k$ on trouve
\begin{equation}\nonumber\tr_\p(x_1\cdots x_{n})+ (-1)^{n-1}\tr_\k(x_1 x_{n}\cdots x_2)=
\tr_\p(x_1\cdots x_{n})+ (-1)^{n-1}\tr_\k( x_{n}\cdots x_2
x_1).\end{equation}
 \hfill$\blacksquare$

\subsection{Cas résoluble}
\begin{prop}\label{propE=1resoluble} Quand $(\g, \sigma)$ est une paire symétrique
 résoluble on a $E=1$. Plus généralement quand  $X$ est  dans un idéal résoluble on a $E(X, Y)=1$ pour tout $Y\in \p$.
\end{prop}

 \begin{cor}

Pour toute paire symétrique  $(\g, \sigma)$, on a  $f\underset{CF}\star g = fg$ pour
 $f,g \in S(\p)^\k$ dès que  $f\in S(J)^\k$ où $J$
 est  idéal $\sigma$-stable et
résoluble de $\g$.
\end{cor}

\noindent \textit{\textbf{Preuve :}} En effet tous les  crochets en $X, Y$ sont unipotents. Donc il n'y a
que des roues pures, nécessairement  de taille $2n$. On a donc des termes de la
forme (avec  $x_i=\ad X $ ou $\ad Y$)

\[w_\G \big(\tr_\p(x_1\ldots x_{2n})-\tr_\k(x_1\ldots
x_{2n})\big). \] Cette expression est nulle  car on peut permuter les termes deux à deux (ça
fait apparaître un élément unipotent), et on peut changer la trace
sur $\k$ en une trace sur $\p$ pour la même raison.

\fin

On peut appliquer la proposition au cas ou $X=\pm Y$ et on obtient
\begin{equation} E(X, X)=E(X, -X)=1.\end{equation}

On montrera au Théorème~\ref{calculCF} (\S\ref{sectionRou=CF}) que le
produit $\underset{CF} \star$ et le produit de Rouvière coïncident
sur $S(\p)^\k$. On déduit de la proposition précédente une  autre
preuve du résultat fondamental de Rouvière \cite{Rou86} :
\begin{theo}[\cite{Rou86}, Théorème~5.1]\label{theoResoluble} La formule de Rouvière est un
isomorphisme pour la convolution des germes de distributions
$K$-invariantes sur l'espace symétrique résoluble~$G/K$.
\end{theo}

\subsection{Propriétés
supplémentaires}\label{sectionpropritesupplementaires}
\begin{lem}

La fonction $E(X, Y)$ ne dépend que d'une composante réductive $\r$
de $\g$. Cette fonction est la même que pour la paire symétrique
dégénérée $\hat{\g}= \r
 \semi (\g_u)_{ab}$.
\end{lem}

\begin{cor}

Le produit $\underset{CF}\star $ ne fait pas intervenir de dérivées
d'éléments unipotents.
\end{cor}

 Une autre conséquence de la Proposition~\ref{alternance} est que
$\ln (E(X, Y))$ s'écrit $\tr_\p(A)$ o\`u $A$ est dans l'idéal
bilatère $\mathcal{I}$ engendré par $\ad[X, Y]=xy-yx$.\footnote{Ce genre de
résultat devrait s'énoncer plus correctement  dans l'algèbre libre
engendrée par $\ad X, \ad Y$.} En effet considérons une roue $\Gamma$ attachée aux mots $X_1, \ldots, X_n$. Si l'un des mots $X_i$ est de longueur plus grande que $1$ alors $A=\ad X_1 \ad X_2 \ldots \ad X_n \in \mathcal{I}$. Si tous les $X_i$ valent $X$ ou $Y$ alors, d'après la Proposition~\ref{alternance}, les contributions se regroupent sous la forme
\[w_\G \big(\tr_\p(x_1\ldots x_{2n})-\tr_\k(x_1\ldots
x_{2n})\big)= w_\G \big(\tr_\p(x_1\ldots x_{2n})-\tr_\p(x_{2n}x_1\ldots
x_{2n-1})\big).\]On écrit $\tr_\p(x_1\ldots x_{2n})=\tr_\p(x_2x_1x_3\ldots
x_{2n})+\tr_\p([x_1, x_2]x_3\ldots x_{2n})$, ce qui montre que, modulo des termes de la forme $\tr_\p(\mathcal{I})$, on a  $\tr_\p(x_1\ldots x_{2n})= \tr_\p(x_{2n}x_1\ldots
x_{2n-1})$. On retrouve ainsi le résultat de \cite{Rou86} Théorème~3.15.\\

\subsection{Calculs numériques à l'ordre $4$ pour $\ln (E_\lambda(X,Y)) $}\label{calculEordrequatre}

\subsubsection{Composante sur $\k$}

On calcule d'abord le poids des graphes qui vont intervenir à
l'ordre $4$. Ces calculs se font en utilisant la formule de Stokes. On trouve  le développement suivant :
\begin{equation}\frac12[X, Y]+ \frac{-1}{24}[X,[X,[X, Y]]]
+\frac{-1}{24}[Y,[Y,[X, Y]]] +\frac{-1}{48}[X,[Y,[X, Y]]]+
\frac{-1}{48}[Y,[X,[X, Y]]].\end{equation}

Lorsque le caractère $\lambda$ de
la sous-algèbre $\k$ vaut $\frac 12\tr_\k \ad $, on peut écrire ce
terme à l'ordre $4$  sous la forme

\begin{equation}\frac14 \tr_\k(xy-yx) + \frac1{12}\tr_\k(yx^3-x^3y+y^3x-xy^3)
+\frac{1}{24} \tr_k\big((yx)^2-(xy)^2\big),\end{equation}
où $x=\ad X$ et $y=\ad Y$. C'est bien l'expression que l'on trouve dans \cite{Rou90} page
256.

\subsubsection{Composante scalaire}

D'après le Lemme~\ref{symetrieE}, la fonction $\ln(E(X, Y)$ est paire. Les contributions à l'ordre $2$ sont nulles, car d'après la Proposition~\ref{alternance} elles sont de la forme $\tr_\p(\ad X \ad Y) - \tr_\k (\ad Y \ad X)=0$.\\

Afin de comparer notre fonction $E(X, Y)$
à celle de Rouvière $e(X, Y)$, on regarde les contributions à l'ordre $4$ en
$X,Y$, en fonction de la taille des roues. \\

Lorsque la roue est de  taille $2$, les arêtes sortantes sont
colorées par $\k^*$, sinon on aurait une arête double. On trouve, pour le
terme de degré $4$ en $X,Y$, une expression de la forme :

\begin{equation}w_\G (\tr_\p-\tr_\k)\Big(\ad[X, Y]\ad[X, Y]\Big)=w_\G b([X,Y], [X, Y]),\end{equation}
où $b(A, B)=K_\g(A, B)-2K_\k(A, B)$ est la différence des  formes de Killing sur $\g$ et
$\k$, comme dans \cite{Rou90}.

Pour les roues de taille $3$, les contributions à l'ordre $4$ font intervenir $2$ sorties dans $\p^*$.
Compte tenu des symétries, le seul cas non trivial concerne les sorties $[X, Y], X, Y$ (dessin de gauche de (\ref{Eordre4})). On trouve une contribution de la forme précédente $ C_3 b([X, Y], [X,
Y])$.

Pour les roues de tailles $4$, qui sont nécessairement des roues pures,  les seules contributions non triviales pour le symbole  font intervenir des sorties attachées sur $X, X, X, Y$ ou $Y, Y, Y, X$ (dessin de droite de (\ref{Eordre4})). Toutefois le coefficient associé est nul car on
peut intervertir les positions $2$ et $4$  ce qui échange $3$ groupes de $2$ arêtes.

\begin{equation}\label{Eordre4}
\begin{array}{lcr}
  \includegraphics[width=7cm]{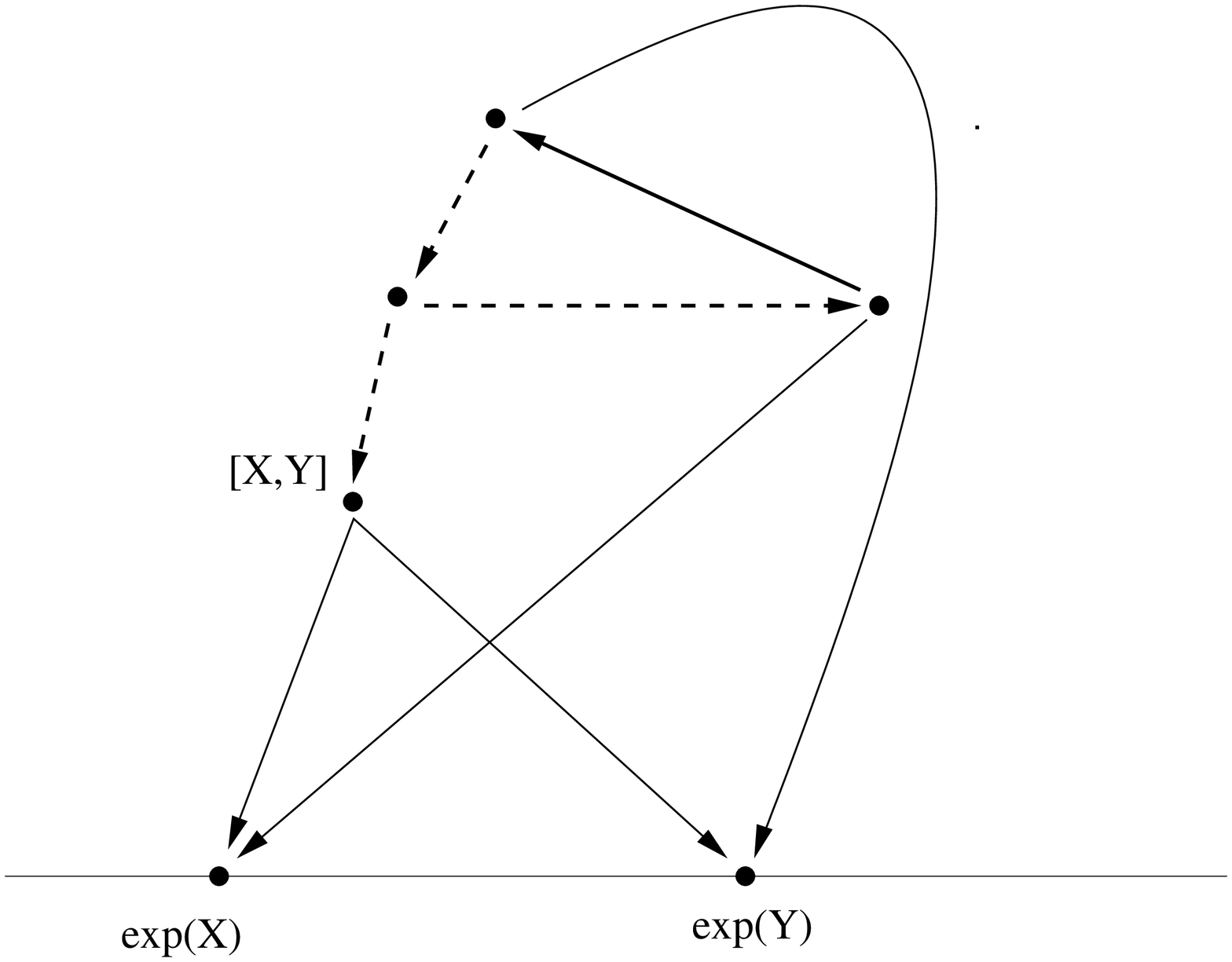} &\quad & \includegraphics[width=7cm]{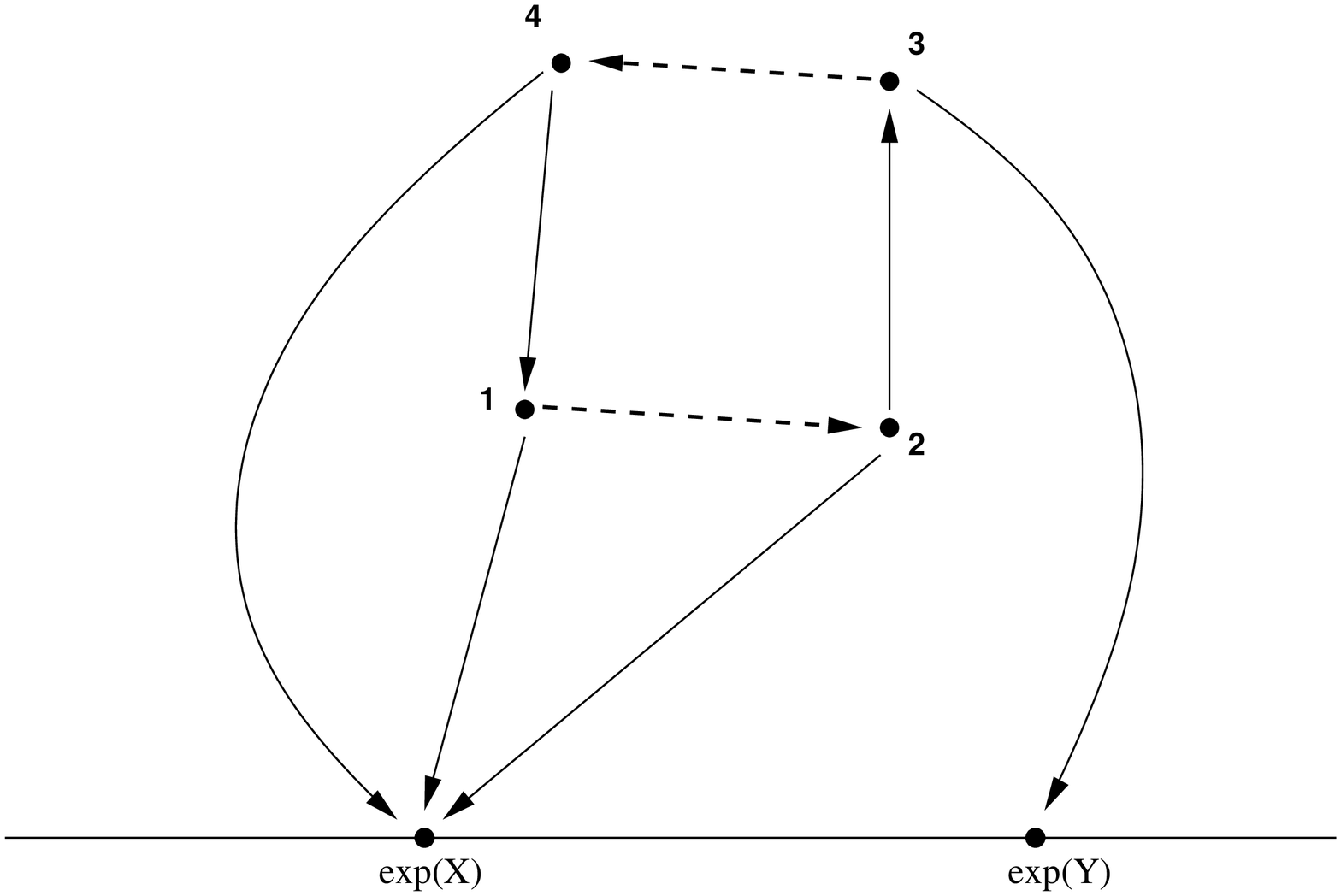} \\
   Roue\;  taille \; 3 \; sur \; [X, Y], X, Y & \quad & Roue \; pure \; taille\; 4 \; sur \; X,X,X,Y
\end{array}
\end{equation}

\noindent On résume ces calculs par la proposition suivante :

\begin{prop}\label{propordre4}
Jusqu'à l'ordre $4$ au moins, la fonction $E_{\frac 12
\tr_k}(X, Y)$ coïncide avec la fonction $e(X, Y)$ de Rouvière.
\end{prop}

\noindent \textit{\bf Preuve : } D'après les calculs précédents les
contributions scalaires à l'ordre $4$ correspondent à un terme de la
forme \[\exp\Big(c (\tr_\p -\tr_\k)(\ad[X,Y])^2 + \ldots \Big).\] On
verra au  Théorème~\ref{calculCF} (\S\ref{sectionRou=CF}) que le
produit de Cattaneo-Felder coïncide avec celui de Rouvière ce  qui
montre que la constante $c$ est uniquement déterminée par un calcul
explicite. On peut le faire dans $sl(2)$ et on trouve $c= 1/240$
résultat conforme (au signe près) à celui de Rouvière (\cite{Rou91}
page 257).

\noindent Voici le détail du calcul: $\widehat{H},\widehat{X},\widehat{Y}$
désigne la base de $sl(2)$ et $\k=\langle \widehat{X}-\widehat{Y}
\rangle$ et $\p=\langle \widehat{H},
\widehat{X}+\widehat{Y}\rangle$. On pose $\omega=\widehat{H}^2+
(\widehat{X}+\widehat{Y})^2$ et $\Omega=\beta(\omega)$. On a pour
$X\in \p$
\begin{eqnarray}
 J^{1/2}(X)=1+\frac1{12} \tr_\p(\ad X)^2
+\frac1{360} \tr_\p(\ad X)^4+\ldots.\end{eqnarray}
On trouve
$\partial_{J^{1/2}}(\omega^2)=\omega^2+\frac{16}3\omega+
\frac{128}{45}$ et $\beta(\omega^2)=\Omega^2 -\frac83 \Omega
\pmod{U(\g) \cdot \k}$. Modulo $U(\g)\cdot \k$ on a donc
\begin{eqnarray}\nonumber \partial_{J^{1/2}}(\omega^2)\cdot \partial_{J^{1/2}}(\omega^2)=\Omega^2 +\frac83
\Omega+\frac{16}3=\beta\left(\partial_{J^{1/2}}(\omega^2-\frac{16}{15})\right).
\end{eqnarray} Pour $sl(2)$ on calcule $(\tr_\p
-\tr_\k)(\ad[X,Y])^2$ agissant comme opérateur bidifférentiel sur
$\omega\otimes \omega$. On trouve $-256$, la constante vaut donc
$\frac{16}{15\times 256}=\frac1{240}$.

\fin
\subsection{Cas quadratique}
 On dit qu'une paire symétrique $(\g, \sigma)$ est quadratique si
elle est munie d'une forme
  bilinéaire invariante et non dégénérée (mais on ne suppose pas que
  c'est la forme de Killing !).
\subsubsection{Cas des paires d'Alekseev-Meinrenken}

 C'est par définition le cas des paires quadratiques  pour
lesquelles la forme bilinéaire est
\textbf{$\sigma$-anti-invariante.} On a alors
  les
identifications $\p^*\sim \k$ et $\k^*\sim \p$.   En utilisant la
Proposition~\ref{alternance} on trouve :

\begin{multline}\tr_\p(x_1\cdots x_{n})+ (-1)^{n-1}\tr_\k(x_{n}\cdots
x_1)=\tr_\p(x_1\cdots x_{n})+ \\(-1)^{n-1}\tr_{\p^*}(x_{n}\cdots
x_1)=\tr_\p(x_1\cdots x_{n}) - \tr_{\p}(x_{1}\cdots x_n)=0.
\end{multline}

On en déduit le résultat important suivant:
\begin{prop}\label{propE=1AM}
Pour une paire symétrique d'Alekseev-Meinrenken la fonction $E$
vaut~$1$.
\end{prop}

 Dans le Théorème~\ref{calculCF}, comme
déjà annoncé,   on fera le lien entre la fonction $E$ et le produit
de Rouvière, ce qui entraînera notamment que notre fonction $E$
décrit la convolution des distributions $K$-invariantes sur l'espace
symétrique $G/K$. Par conséquent on retrouve un résultat démontré
par Alekseev-Meinrenken \cite{AM} et Torossian \cite{To4}, à savoir
que la symétrisation, modifiée par la racine carrée du jacobien, est
un isomorphisme  d'algèbres dans le cas des paires symétriques
 avec forme $\sigma$-anti-invariante. Dans notre situation on
démontre une conjecture formulée dans \cite{To4} : cet isomorphisme
s'étend aux
germes de distributions invariantes car la fonction $E(X,Y)$ vaut~$1$. On résume la situation par le théorème suivant. \\

\begin{theo}\label{theoAM}
Pour une paire symétrique d'Alekseev-Meinrenken, la formule de
Rouvière est  encore un isomorphisme au niveau des germes de
distributions $K$-invariantes.

\end{theo}

\subsubsection{Cas quadratique $\sigma$-invariant}

  C'est par définition le cas des paires quadratiques  pour
lesquelles la forme bilinéaire est \textbf{$\sigma$-invariante}. On
a alors les identifications
 $\p\sim \p^*$ et $\k \sim
\k^*$  et l'endomorphisme transposé vérifie $^{t}\ad X= -\ad X$. On
a donc :

\begin{multline}\tr_\p(x_1\cdots x_{n})+ (-1)^{n-1}\tr_\k(x_{n}\cdots
x_1)=\\\tr_\p(x_1\cdots x_{n})+ (-1)^{n-1}\tr_{\k^*}(x_{n}\cdots
x_1)=\tr_\p(x_1\cdots x_{n})- \tr_{\k}(x_{1}\cdots x_n).
\end{multline}
Si de plus,  la paire symétrique est très-symétrique\footnote{Une
paire symétrique est dite très-symétrique s'il existe un opérateur
$A$ sur $\g$ tel que
 $A\circ \ad X= \ad X \circ A$ avec $A: \k
\longrightarrow \p$ et $\p\longrightarrow\k$.}, alors  par un calcul
élémentaire en termes de matrices blocs $2\times 2$, on trouve :
\[\tr_\p(x_1\cdots x_{n})- \tr_{\k}(x_{1}\cdots x_n)=0.\]
On en déduit la proposition suivante.

\begin{prop}\label{propE=1tressym}
Si $(\g, \sigma)$ est une  paire très-symétrique et quadratique
avec forme $\sigma$-invariante, alors la fonction $E$ vaut~$1$.
\end{prop}

Par exemple une algèbre de Lie considérée comme un espace symétrique
est une paire très-symétrique, car on a $\sigma(x,y)=(y,x)$ et
$A(x,y)=(x, -y)$. Pour une paire symétrique complexe, on peut
prendre pour $A$ la multiplication par $i$. La proposition ci-dessus
s'applique notamment dans le cas des algèbres de Lie quadratiques
vues comme des paires symétriques $\g \times \g/\mathrm{diagonale}$.

\begin{cor}
Si $\g$ est une algèbre de Lie quadratique vue comme une paire
symétrique $\g \times \g/\mathrm{diagonale}$, alors la fonction $E$ vaut~$1$.
\end{cor}

\paragraph{Conjecture 1:} Pour les algèbres de Lie, vues comme paires symétriques
$\g \times \g/\mathrm{diagonale}$, la fonction $E$ vaut $1$.

 \section[Opérateurs différentiels invariants]{Écriture en coordonnées exponentielles des opérateurs dif\-féren\-tiels invariants sur un espace symé\-trique}\label{sectionopdiff} On considère une
paire symétrique $(\g, \sigma)$. Dans le diagramme de
Cattaneo-Felder pour la bi-quantification on considère le couple de
variétés co-isotropes $C_2= \k^\perp$, mis en position verticale, et
$C_1=0^\perp=\g^*$, mis en position horizontale (\textit{cf.}
\S~\ref{soussectionbiquantification}). Comme dans \cite{Kont} on
considère
 pour le bi-vecteur de Poisson, la moitié du crochet de Lie. Alors $(U(\g), \cdot) $ est
isomorphe à $(S(\g), \underset{DK}\star)$.
\subsection{Liens entre les produits de Rouvière et de Cattaneo-Felder}

\subsubsection{Produit de Duflo-Kontsevich}\label{produitDK}

Rappelons que le produit
de Duflo-Kontsevich sur $S(\g)$ est donné par la formule

\begin{equation}\label{produitDuflo}\beta\Big(\partial_{q^{1/2}} (f\underset{DK}\star
g)\Big)=\beta\big(\partial_{q^{1/2}} f\big)\cdot
\beta\big(\partial_{q^{1/2}} g\big),\end{equation}

\noindent où $q$ désigne  la fonction
\begin{equation}\label{foncionq}q(X)=\det_{\g}\Big(\frac{\sinh \frac{\ad X}2}{\frac{\ad
X}2}\Big).\end{equation}
\subsubsection{Roues horizontales et verticales}\label{sectionrouesHV}
 L'espace de réduction, pour la partie verticale, est l'algèbre
$S(\p)^\k$ munie du produit~$\underset{CF}\star$. L'espace de
réduction, pour la partie horizontale,
 est l'algèbre $S(\g)$ munie du produit de Duflo-Kontsevich. \\

\begin{figure}[h!]
\begin{center}
\includegraphics[width=8cm]{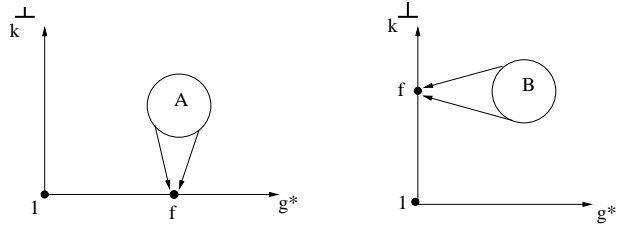}
\caption{\footnotesize Contributions des roues pures sur les axes}\label{rouesAB.eps}
\end{center}
\end{figure}

Pour $f\in S(\g)$, notons $A(f)$ les contributions sur l'axe horizontal $\g^*$  (\textit{cf.} Fig. \ref{rouesAB.eps})

\begin{equation}1\underset{1}\star f=A(f).\end{equation}
Alors $A$ est un opérateur à valeurs dans $S(\p)$, donné par des graphes de Kontsevich à
$4$ couleurs, qui n'est pas restreint aux roues pures. Cet opérateur est compliqué et n'est pas à coefficients
constants car les graphes avec arêtes
doubles colorées par $+\, + $ et $+ \, -$ ne sont pas nuls \textit{a priori}.\\

En fait, l'opérateur $A$ correspond à une modification de la
projection de $S(\g)$ sur $S(\p)$ via la symétrisation (modifiée),
c'est-à-dire la décomposition $U(\g)= \beta(\partial_{q^{1/2}}S(\p))
\oplus \k\cdot U(\g)$. En effet, lorsqu'on munit $S(\g)$ du produit
de Duflo-Kontsevich, on a $S(\g)=S(\p)\oplus \k\underset{DK}\star
S(\g)$ et donc, grâce à l'application $\beta\circ
\partial_{q^{1/2}}$, on a :

\[ U(\g)= \beta\circ
\partial_{q^{1/2}} (S(\p)) \oplus \k\cdot U(\g).\]
 En particulier on a l'équivalence \[1\underset{1}\star f=0 \quad
\Longleftrightarrow \quad f\in \k\underset{1}\star S(\g).\]

Lorsque $f$ est dans $S(\p)$, l'opérateur $A$ correspond aux
contributions des roues pures; c'est un opérateur à coefficients
constant  que l'on notera encore $A$.\\

Notons $B(f)$ les contributions sur l'axe vertical  $\k^\perp$; c'est un opérateur à coefficients constants correspondant
aux contributions des roues pures attachées sur $\k^\perp$ :
\begin{equation}f\underset{2}\star 1=B(f).\end{equation}

\noindent On notera $A(X)$ et $B(X)$\footnote{On verra dans
\S~\ref{sectionsigmastable} prop. \ref{B=1} que l'on a $B(X)=1$ pour
$X\in \p$. Mais on n'utilise pas ici ce résultat.} les symboles
associés. On a donc pour $X\in \p$,
\begin{equation}A(X)=(1\underset{1}\star e^X)e^{-X} \quad \mathrm{et} \quad
B(X)=(e^X\underset{2}\star 1)e^{-X}.\end{equation}
\noindent Ce sont des symboles inversibles.
\subsubsection{Comparaison des star-produits  $\sharp_{Rou}$ et $\star_{CF}$}\label{sectionRou=CF}

Soient  $f, g \in S(\p)^\k$. Par compatibilité des produits
$\underset{1}\star$ et $\underset{2}\star$ en cohomologie,  on a $f\underset{2}\star
1=1\underset{1}\star A^{-1}B(f)$ et par suite
\[(f\underset{CF}\star g)\underset{2}\star 1= 1\underset{1}\star
A^{-1}B(f\underset{CF}\star g),\] puis

\[
(f\underset{CF}\star g)\underset{2}\star 1=1\underset{1}\star
\Big(A^{-1}B(f)\underset{DK}\star A^{-1}B(g)\Big).
\]
On en déduit que l'on a

\begin{equation}\label{projection}
A^{-1}B(f\underset{CF}\star g)=A^{-1}B(f)\underset{DK}\star
A^{-1}B(g) \quad \mathrm{modulo} \quad \k\underset{DK}\star
S(\g).\end{equation}

Comme $\beta\circ\partial_{q^{1/2}}$ transforme le produit de Duflo
en le produit dans l'algèbre enveloppante, on en déduit que
\[\beta\circ\partial_{q^{1/2}}\circ A^{-1}B\] transforme le produit
$\underset{CF}\star $ en le produit dans l'algèbre enveloppante
modulo~\mbox{$\k\cdot U(\g)$}, c'est-à-dire que l'on a
$\mathrm{modulo} \quad \k\cdot U(\g)$ :

$$\left(\beta\circ\partial_{q^{1/2}}\circ
A^{-1}B\right)(f\underset{CF}\star g)=
\left(\beta\circ\partial_{q^{1/2}}\circ
A^{-1}B\right)(f)\left(\beta\circ\partial_{q^{1/2}}\circ
A^{-1}B\right)(g)$$
Soit $J(X)=\det\limits_{\p}\Big(\frac{\sinh(\ad X)}{\ad X}\Big)$. Remarquons que pour $X\in \p$, on a $q^{1/2}(X)=J(\frac X2)$. Le produit de Rouvière \cite{Rou86} est défini par la formule, pour $p,q \in S(\p)$ :

$$\beta\Big(\partial_{J^{\frac 12}}(p\underset{Rou}\sharp
q)\Big)=\beta\Big(\partial_{J^{\frac
12}}(p)\Big)\cdot\beta\Big(\partial_{J^{\frac 12}}(p)\Big) \quad
 \mathrm{modulo}\quad \k\cdot U(\g).$$

\noindent On va voir que les deux produits coïncident sur les éléments $\k$-invariants.

\begin{theo}\label{calculCF} Les produits de Rouvière et de
Cattaneo-Felder coïncident pour toute paire symétrique.  On a, pour
$f, g \in S(\p)^\k$,
\[\beta\Big(\partial_{J^{\frac 12}}(f\underset{CF}\star
g)\Big)=\beta\Big(\partial_{J^{\frac
12}}(f)\Big)\cdot\beta\Big(\partial_{J^{\frac 12}}(g)\Big) \quad
 \mathrm{modulo}\quad \k\cdot U(\g).\]
 \end{theo}

\noindent \textit{\textbf{Preuve du théorème:}}  La preuve se fait
comme dans \cite{Kont}. On compare deux isomorphismes.
 On démontre un lemme et une formule intermédiaire. \\

\begin{lem} Il existe des algèbres de Lie résolubles pour
lesquelles   $\tr_\p(\ad X)^{2n}$ n'agit pas comme une
dérivation sur $S(\p)^\k$.
\end{lem}
\noindent \textit{\textbf{Preuve du lemme:}} Soit $\g$ l'algèbre de
Lie résoluble engendrée par $t, x, y, z$ et les relations
\[[x,y]=z, \quad [t,x]=-x, \quad [t,y]=y.\]
 C'est une paire symétrique pour l'involution
\[\sigma(t)=-t,\quad \sigma(x)=y,\quad \sigma(z)=-z.\]Le sous-espace $\p$ est
engendré par $t, x-y, z$ et le sous-espace $\k$ est engendré par
$x+y$. On vérifie facilement que l'on a  $S(\p)^\k=\C[z, 4zt+
(x-y)^2]$. La forme $t^*$ est clairement $\k$-invariante car on a
$t^*([\k, \p])=0$. Par ailleurs l'application polynomiale $X\mapsto
\tr_\p(\ad X)^{2n}$ s'identifie à  $(t^*)^{2n}$, qui n'agit pas comme
une dérivation dans $S(\p)^\k$.

\fin

\begin{prop}\label{AB}
On a l'égalité de fonctions universelles suivante\footnote{Compte
tenu du fait que $B=1$ comme on le verra dans
\S~\ref{sectionsigmastable} prop. \ref{B=1}, cette proposition donne
une formule pour $A(X)$ lorsque $X\in \p$.} pour $X\in \p$,

\begin{equation}\label{formuleAJ=Bq}A(X)J^{\frac 12}(X)=B(X) q^{\frac 12}(X).
\end{equation}

\end{prop}

\noindent \textit{\textbf{Preuve de la proposition:}} Les fonctions
$A, B$, $q, J$ ne font intervenir que les $\tr_\p (\ad
 X)^{2n}$ et sont universelles.  Lorsque  $\g$ est résoluble,
 on sait que l'on a les égalités $E=e=1$ donc le
produit de Rouvière et le  produit $\underset{CF}\star$
correspondent au  produit standard. Raisonnons par l'absurde. On
disposerait d'un automorphisme non trivial pour le produit standard dans
$S(\p)^\k$ donné par l'exponentielle d'une série en $\tr_\p (\ad
 X)^{2n}$. On en déduirait qu'il existerait $n>0$ tel que $\tr_\p
(\ad  X)^{2n} $ agirait comme une dérivation sur $S(\p)^\k$ ce qui
n'est
pas vrai d'après le lemme ci-dessus. \fin\\

\noindent \textit{(fin de la preuve du théorème)} Grâce à la
Proposition~\ref{AB} et à l'équation (\ref{projection}) le
produit $\underset{CF}\star$ et le produit de Rouvière coïncident sur les éléments $\k$-invariants.
\fin

\paragraph{Remarque 7: } Pour $X\in \p$ on a  $q^{\frac 12}(X)=J(\frac
X2)$ et d'après la Proposition~\ref{B=1} (\S~\ref{sectionsigmastable})
on a $B(X)=1$. On en déduit la formule suivante pour $A(X)$
\begin{equation}A^{-1}(X)q^{\frac 12}(X)=A^{-1}(X)J(\frac X2 )=J^{\frac 12}(X).\end{equation}

\paragraph{Remarque 8: }  Si on avait inversé  le role de $C_1=\g^*$ et
 $C_2=\k^\perp$ on aurait trouvé de même, pour $f,g \in S(\p)^\k$ :

\[\beta\Big(\partial_{J^{\frac
12}}(f\underset{CF}\star g)\Big)=\beta\Big(\partial_{J^{\frac
12}}(f)\Big)\cdot\beta\Big(\partial_{J^{\frac 12}}(g)\Big) \quad
 \mathrm{modulo}\quad U(\g)\cdot \k.\]

\subsection{Cas à paramètre}

\subsubsection{Commutativité d'algèbres d'opérateurs
différentiels invariants sur les $z$-densités}

On met en position verticale maintenant la sous-variété co-isotrope
$\lambda + \k^\perp$ avec $\lambda$ un caractère de $\k$ et $\g^*$
en position horizontale.   Notons  $\k^\lambda$ le sous-espace
 de $S(\g)$ défini par  $\{X+\lambda(X), \; X\in \k\}$.   Comme précédemment on
aura, pour $f,g $ deux éléments $\k$-invariants :

\[ A^{-1}B(f\underset{CF , \, \lambda }\star g)=A^{-1}B(f)\underset{DK}\star
A^{-1}B(g) \quad \mathrm{modulo} \quad
\k^{-\lambda}\underset{DK}\star S(\g),\] car on a $1\underset{1}\star
K=\lambda(K)$.\\

\noindent \textbf{Action de la symétrie par rapport à la diagonale :
} Faisons agir $s$ la symétrie par rapport à la première diagonale. Les fonctions d'angle vérifient
les relations suivantes :
\begin{equation}\label{symetriesangle}
\begin{array}{ccc}
\mathrm{d} \phi_{+ +}(p,q) &= & \mathrm{d}\phi_{- -}(q,p)\\
\mathrm{d} \phi_{+ +}(s(p), s(q))&= &-\mathrm{d} \phi_{+ +}(p, q)\\
\mathrm{d}\phi_{+ -}(s(p), s(q))&= &-\mathrm{d} \phi_{- +}(p,q).
\end{array}
\end{equation}
Comme les roues pures $A$ et $B$ sont paires, l'action de la symétrie $s$ les préserve et nos notations restent cohérentes ;  $A$ est l'action des roues pures attachées sur $\g^*$ que l'on choisisse la position horizontale ou verticale.\\

L'action de la symétrie revient à  mettre en position horizontale la sous-variété $\lambda +
\k^\perp$ et en position verticale~$\g^*$.\\

Considérons deux fonctions $f, g$ positionnées sur
l'axe horizontal et regardons le produit limite quand les positions de $f$ et $g$ se rapprochent. Dans la première configuration limite ($\g^*$ en position horizontale) on retrouve le produit de Duflo-Kontsevich
$f\underset{DK}\star g$. Par symétrie  diagonale  on trouvera dans
la deuxième configuration limite ($\g^*$ en position
verticale) $g\underset{DK, \; signes}\star f$, les signes provenant de l'action de $s$ sur les fonctions d'angle. Or les coefficients avec fonction d'angle à une couleur (ou deux couleurs) ont une propriété de symétrie par rapport à l'axe vertical ($(-1)^n
w_\G=w_{\overset{\wedge} \G}$). En tenant compte de cette symétrie on constate que $g\underset{DK, \; signes}\star f$ vaut $f\underset{DK}\star g$.

 Du coté du produit de
Cattaneo-Felder pour le paramètre $\lambda$, lorsque $f$ et $g$ sont positionnées sur l'axe vertical,  on trouvera dans la
première configuration limite ($\lambda + \k^\perp$ en position verticale )  $f\underset{CF,\, \lambda} \star g$ et par
symétrie dans l'autre configuration limite ($\lambda + \k^\perp$ en position
horizontale)  $g\underset{CF,\,-
\lambda} \star
f$.\\

On en déduit les égalités, pour $f,g \in S(\p)^\k$:
\begin{eqnarray}\label{gauche-droite1} A^{-1}B(f\underset{CF,\,\lambda
}\star g)=A^{-1}B(f)\underset{DK}\star A^{-1}B(g) \quad
\mathrm{modulo} \quad \k^{-\lambda}\underset{DK}\star S(\g)\\\label{gauche-droite2}
 A^{-1}B(g\underset{CF,\,-\lambda }\star
f)=A^{-1}B(f)\underset{DK}\star A^{-1}B(g) \quad \mathrm{modulo}
\quad S(\g)\underset{DK}\star \k^{-\lambda}.
\end{eqnarray}

Or d'après  le Lemme~\ref{symetrieE} (\S~\ref{defiE}) on a
$f\underset{CF,\, \lambda }\star g=g\underset{CF,\, -\lambda }\star
f$ et d'après un résultat de Duflo~\cite{du79}, résultant de la
dualité de Poincaré, on a :
\[\k^{-\lambda}\cdot U(\g)\cap U(\g)^\k=U(\g)^\k\cap  U(\g)\cdot
\k^{-\lambda + \tr_\k }.\]
 On en déduit le théorème  suivant, qui étend le critère
 de commutativité des algèbres d'opérateurs différentiels
 sur les demi-densités\footnote{ cas $z=\frac 12$  \cite{du79}, on utilise la relation
 $\frac 12 \tr_\k=-\frac 12 \tr_{\g/\k}+\frac 12 \tr_\g$ et le fait que l'action
 du caractère $\tr_\g$ n'est pas significative.}
 et sur les fonctions.\footnote{cas $z=0$
 \cite{lich}.}

\begin{theo}\label{theoz}

 Pour $f, g$ dans $S(\p)^\k$ et pour tout $z\in \R$, on a la relation
  \[f\underset{CF,\, \lambda}\star g= f\underset{CF,\,\, \lambda
+ z \, \tr_\k }\star g.\] En particulier le produit naturel dans
$\big(U(\g)/U(\g)\cdot \k^{z\, \tr_{\k}}\big)^\k$, l'algèbre des
opérateurs différentiels invariants sur les $z$-densités, est
commutatif pour tout $z$.
\end{theo}

\noindent \textit{\textbf{Preuve :}} Les fonctions $E_{\lambda
+\tr_\k}$ et $E_\lambda$ définissent les mêmes star-produits sur
les éléments $\k$-invariants d'après la relation de Duflo et les formules (\ref{gauche-droite1}), (\ref{gauche-droite2}). Il résulte par récurrence qu'il en est de même
pour $E_{\lambda + n \tr_\k}$ et par prolongement des identités
polynomiales pour $E_{\lambda + z\,\tr_\k}$.

\fin
\subsubsection{Comparaison des fonctions $E(X, Y)$ et $e(X, Y)$}
On applique  à l'équation~(\ref{gauche-droite2}), l'opérateur
$\beta(\partial_{q^{1/2}})$ qui transforme le produit
$\underset{DK}\star$ en le produit dans l'algèbre enveloppante.  Compte tenu de (\ref{formuleAJ=Bq}), on trouve

\[ \beta(\partial_{J^{1/2}}(f\underset{CF,\, \lambda}\star
g))=\beta(\partial_{J^{1/2}}f) \cdot \beta(\partial_{J^{1/2}}g)
\quad \mathrm{modulo} \quad U(\g) \cdot \k^{-\lambda}.\]
Les
fonctions $e_{-\lambda}(X,Y)$ vérifient le même type de relation.
On en déduit que $E_\lambda(X, Y)$ et $e_{-\lambda}(X,Y)$ doivent
définir
le même star-produit. \\

Or d'après \cite{Rou94} les fonctions $e_{-\lambda}(X,Y)$ ne sont
symétriques que pour $-\lambda =\frac 12 \tr_\k$ et  les premiers
termes du développement calculé en Proposition \ref{propordre4} (\S~\ref{calculEordrequatre})
laissent à penser que l'on a l'égalité conjecturale $E_{\frac 12 \tr_\k}(X,Y)\overset{??}=e(X, Y)$.
On devrait avoir aussi  de manière conjecturale l'égalité $E_\lambda(X, Y)\overset{??}=e_{-\lambda+\frac 12 \tr_k}(X,Y)$.

\paragraph{Conjecture 2 : } \textit{On a $E_\lambda(X, Y)=e_{-\lambda+\frac 12
\tr_k}(X,Y)$.}\\

\subsection[Opérateurs différentiels en coordonnées exponentielles]{Opérateurs différentiels invariants en coordonnées
exponentielles}\label{sectionOpdiffexp}

On considère le couple de variétés co-isotropes $C_1=\k^\perp$ et
$C_2=0^\perp=\g^*$. Le bi-vecteur de Poisson vaut la moitié du crochet de Lie et $(U(\g), \cdot) $ est isomorphe à $(S(\g),
\underset{DK}\star)$.\\

Les espaces de réduction sont, pour la partie horizontale, l'algèbre
$S(\p)^\k$ munie du produit
$\underset{CF}\star=\underset{Rou}\sharp$ et pour la partie
verticale, l'algèbre $S(\g)$ munie du produit de Kontsevich. L'espace
de réduction à l'origine est $S(\p)$, c'est-à-dire les distributions
portées par l'origine
sur l'espace symétrique $G/K$.\\

On en déduit pour $S(\p)$ l'existence d'une action à gauche de
$(S(\g), \underset{Kont}\star)$ et d'une action à droite de
$(S(\p)^\k, \underset{CF}\star)$.

Soit $R\in S(\p)^\k$.  On place dans le quadrant de Cattaneo-Felder,
$e^X$ à l'origine et $R$ sur l'axe horizontal: on obtient donc
$e^X\underset{1}\star R$. C'est donc un élément de $S(\p)$ qui
dépend de $X$, que l'on voit comme un opérateur différentiel sur
$\p$ au point $X$.\footnote{Le terme $e^X$ est vu comme la
distribution de Dirac au point $X$.} On a le résultat suivant qui
résout de manière satisfaisante
un problème posé par  Duflo (\cite{duflo-japon} problème 7).\\

\begin{theo}\label{theoecritureexp}
La formule $e^X\underset{1}\star R$ est l'écriture, en coordonnées
exponentielles du conjugué par l'application $J^{1/2}(X)$ de
l'opérateur différentiel invariant sur l'espace symétrique $G/K$,
donné par  $D_{\beta(J^{1/2}(\partial)R)}$.
\end{theo}

\noindent \textit{\textbf{Preuve :}} Soit $f\in
\mathcal{C}^\infty(G/K)$ et $R\in S(\p)^k$. On a
\[D_{\beta(R)} \,f(gK)= \langle R^{(Y)}, \, f(g\exp(Y)K)\rangle\]
o\`u $R$ est vu comme une distribution de support $Y=0$. On note
$\Exp$ la fonction $X\in \p \mapsto \exp_G(X)K \in G/K$ avec  $X\in
\p$.\\

\noindent Si on tient compte du facteur $q^{1/2}$  entre
$\mathcal{C}^\infty(G/K)$ et $\mathcal{C}^\infty(\p)$, on pose pour $X\in \p$,  $$\phi(X):=q^{1/2}(X) f(\exp_G(X)K).$$
Pour $u=\beta(\partial_{q^{1/2}} R)$ et on note $D_u$ l'opérateur
différentiel invariant à gauche sur l'espace symétrique associé à
$u$. Il vient en écriture exponentielle

\begin{equation}
D^{Exp}_u (\phi)(X) :=q^{1/2}(X)D_u(f)(\exp_G(X)K)=\\\langle \,
(\partial_{q^{1/2}} R)^{(Y)}, \, q^{1/2}(X)
\left(\frac{\phi}{q^{1/2}}\right)(P(X,Y))\, \rangle,
\end{equation}
 où on a écrit $\exp_G(X)\exp_G(Y)=\exp_G(P(X,Y))\exp_G(K(X,Y))$
avec $P(X,Y) \in \p$ et $K(X,Y)\in \k$. Ces facteurs $P, K$
dépendent
de manière analytique de $X, Y$ dans un voisinage de $0$. Le facteur $P(X, Y)$ est la série de Campbell-Hausdorff pour les espaces symétriques. \\

Rappelons que le produit de Kontsevich-Duflo est donné par la
formule (\ref{produitDuflo}). Pour $X, Y$ dans $\p$, on peut donc écrire :
\[e^X\underset{DK}\star e^Y= \frac{D(X,Y)}{D(P, K)}\,
e^P\underset{DK}\star e^K,\]
avec $D(X,Y)=\frac{q^{1/2}(X)q^{1/2}(Y)}{q^{1/2}(Z(X,Y))}$ la fonction de densité sur
les algèbres de Lie.\\

On a $e^K\underset{2} \star 1 = C(K),$ où $C(K)$ désigne les contributions
des roues pures  attachées à l'axe vertical avec rayons colorés par la couleur $ - \, +$. Comme
$e^K\underset{DK}\star e^{-K}=q(K)$ on doit avoir
$e^K\underset{2}\star (e^{-K}\underset{2}\star 1)=C(K)C(-K)=q(K)$.
On a $K(X,Y)=-K(Y,X)$, donc la partie symétrique de $C(K)$ vaut
$q^{1/2}(K)$. En fait on la résultat intermédiaire suivant.

\begin{lem} La fonction $C$ est paire et vaut $q^{1/2}(K)$.
\end{lem}
\noindent \textit{\textbf{Preuve du lemme:}} On a
\[\left(e^{K/2}\underset{DK}\star e^{K/2}\right)\underset{2}\star 1=\frac{
q(K/2)}{q^{1/2}(K)}e^K\underset{2}\star 1=C(K)\frac{ q(K/2)}{q^{1/2}
(K)}.\] Par ailleurs on a aussi

\[\left(e^{K/2}\underset{DK}\star e^{K/2}\right)\underset{2}\star
1=e^{K/2}\underset{2}\star\left(e^{K/2}\underset{2}\star 1 \right)=
C(K/2)^2,\] d'o\`u l'on tire
\[\left(\frac{C(K/2)}{q^{1/2}(K/2)}\right)^2=\frac{C(K)}{q^{1/2}(K)},\]puis
par récurrence\[\left(\frac{C(\frac K{2^n})}{q^{1/2}(\frac
K{2^n})}\right)^{2^n}=\frac{C(K)}{q^{1/2}(K)}.\] Comme la fonction
$C$ et la fonction $q^{1/2}$ s'écrivent $\exp\left(\sum\limits_{n\geq 2}w_n\tr_\p(\ad
K)^n+\sum\limits_{n\geq 2} w_n'\tr_\k(\ad K)^n\right)$, on en déduit
que l'on a $\frac{C(K)}{q^{1/2}(K)}=\lim\limits_{n\mapsto \infty}\left(\frac{C(\frac K{2^n})}{q^{1/2}(\frac
K{2^n})}\right)^{2^n}=1$. \fin

\noindent \textit{(Fin de la preuve du théorème :)} On note comme en
\S\ref{sectionrouesHV}, $A$ l'opérateur vertical (correspondant à
$\g^*$) des contributions des roues pures et $B$ l'opérateur
horizontal\footnote{Les axes sont inversés par rapport à la situation \S\ref{sectionrouesHV} et  $B=1$ d'après la Proposition~\ref{B=1}
 (\S~\ref{sectionsigmastable}).}. Grâce au lemme  ci-dessus, on en
déduit que l'on a

$$(e^X\underset{DK}\star e^Y) \underset{2}\star
1=\frac{D(X,Y)}{D(P,K)}e^P \underset{2}\star (e^K\underset{2}\star
1)= q^{1/2}(X)q^{1/2}(Y)\frac{A(P)}{q^{1/2}(P)}e^P.$$

\noindent D'après la Proposition~\ref{AB}, on a pour $X \in \p$,  $A(X)J^{\frac 12}(X)=B(X) q^{\frac 12}(X)$, d'où l'on tire l'égalité

\begin{equation}\label{eqeXeY}
(e^X\underset{DK}\star e^Y) \underset{2}\star 1=
q^{1/2}(X)q^{1/2}(Y)\frac{B(P)}{J^{1/2}(P)} e^P.\end{equation}

Par ailleurs le terme de gauche de (\ref{eqeXeY}) peut s'écrire, compte tenu de la
compatibilité des produits :

\[ A(Y)B^{-1}(Y)\, e^X \underset{2}\star (1\underset{1}\star e^Y).\]

Il vient donc la formule

\begin{equation}e^X \underset{2}\star (1\underset{1}\star
e^Y)=J^{1/2}(Y)q^{1/2}(X)\frac{B(P)}{J^{1/2}(P)} e^P.\end{equation}

En différentiant par rapport à $Y$ selon le polynôme $R$, on peut
remplacer $e^Y$ par $R$. On peut alors utiliser la $\k$-invariance de $R$
 pour justifier de la compatibilité des produits
$\underset{1}\star $ et $\underset{2}\star$. Il vient alors

\begin{equation}\nonumber
e^X \underset{2}\star (1\underset{1}\star R) =
R^{(Y)}
J^{1/2}(Y)q^{1/2}(X)\frac{B(P)}{J^{1/2}(P)} e^{P(X, Y)}= (e^X
\underset{2}\star 1)\underset{1}\star R =A(X)e^X\underset{1} \star
R,
\end{equation} c'est-à-dire compte tenu de la relation
$AJ^{\frac12}=Bq^{\frac12}$ on a :

\begin{equation}e^X\underset{1} \star R =R(\partial_Y)
\frac{J^{1/2}(Y)J^{1/2}(X)}{B(X)} \frac{B(P)}{J^{1/2}(P)}
e^P|_{Y=0}.\end{equation}

Cette formule est précisément l'écriture de l'opérateur
$D_{\beta(J^{1/2}(\partial)R)}$ en coordonnées exponentielles sur
l'espace symétriques $G/K$,  modifié par la fonction
$J^{1/2}(X)/B(X)$. Compte tenu du résultat maintes fois annoncé
$B=1$, le théorème est démontré. \fin

\subsection{Déformation de la formule de Baker-Campbell-Hausdorff pour les paires symé\-triques}

L'application $\Exp$ définit un difféomorphisme local de $\p$
sur $G/K$. On en déduit alors l'existence d'une
formule de Campbell-Hausdorff pour les espaces symétriques, définie
de la manière suivante. Pour $X, Y \in \p$ proche de $0$ il existe
une série convergente $Z_{sym}(X,Y)$ à valeurs dans $\p$ telle que
l'on ait
$$\exp_G(X)\Exp(Y)=\Exp\big(Z_{sym}(X,Y)\big).$$
C'est le facteur $P(X, Y)$ introduit dans la section précédente.
En utilisant l'involution $\sigma$ on trouve facilement
\begin{equation}
\exp_G(2Z_{sym}(X,Y))=\exp_G(X)\exp_G(2Y)\exp_G(X).
\end{equation}

Dans \cite{KV} Kashiwara-Vergne ont conjecturé que la déformation par dilatation de la formule de Campbell-Hausdorff
était gouvernée par les champs adjoints. Cette conjecture a été démontrée  dans~\cite{AM2} comme conséquence de \cite{To1} et de la quantification par déformation. Dans \cite{Rou86} F. Rouvière  propose de faire de même pour la série $Z_{sym}(X,Y)$ en utilisant les champs $\k$-adjoints. Dans \cite{To2} une approche via la quantification de Kontsevich fut proposée.

Dans cette section nous montrons qu'une déformation par les champs $\k$-adjoints vient naturellement grâce au diagramme de bi-quantification.
\subsubsection{Première déformation}

Considérons le diagramme de bi-quantification de la section précédente, en mettant $\g^*$ en position verticale et $\k^\perp$ en position horizontale. Considérons $X, Y\in \p$ et mettons les fonctions $e^X$ et $e^Y$ sur l'axe vertical comme précédem\-ment.\\

La contribution des graphes colorés est indépendante de la position de $e^X$ et $e^Y$ sur l'axe vertical. En position limite quand $e^X$ et $e^Y$ se rapprochent,  on trouve  d'après l'équation (\ref{eqeXeY}) et $B=1$ (Proposition~\ref{B=1} (\S~\ref{sectionsigmastable})):

$$(e^X\star_{DK} e^Y)\star_2 1= \frac{q^{1/2}(X)q^{1/2}(Y)}{J^{1/2}(Z_{sym}(X, Y))} e^{Z_{sym}(X, Y)}.$$
Quand $Y$ tend vers l'origine, on trouve

$$A(Y) e^X\star_2 e^Y.$$
Ces deux expressions sont donc égales. Mettons maintenant $e^X$  sur l'axe vertical en position~$1$ et $e^Y$ sur l'axe horizontal en position~$s$. Quand $s$ tend vers $0$ on trouve  $$e^X\star_2 (1\star_1 e^Y)=e^X\star_2 e^Y=\frac{q^{1/2}(X)J^{1/2}(Y)}{J^{1/2}(Z_{sym}(X, Y))} e^{Z_{sym}(X, Y)}.$$ Quand $s$ tend vers l'infini on trouve $A(X) e^X \star_1 e^Y$. Il est alors opportun de considérer les contributions divisées par $A(X)$ pour faire apparaître le facteur $D_{sym}(X, Y)=\frac{J^{1/2}(X)J^{1/2}(Y)}{J^{1/2}(Z_{sym}(X, Y))}$.\\

\noindent Pour $s$ quelconque, les contributions totales divisées par $A(X)$ sont de la forme $$D_s^{(1)}(X, Y) e^{Z_{sym,s}^{(1)}(X, Y)},$$ où
$D_s^{(1)}(X, Y)$ est une fonction de densité et $Z_{sym,s}^{(1)}(X, Y)\in \p$ est une déformation de la fonction de Campbell-Hausdorff.  La déformation est alors contrôlée par les concentrations en $e^Y$. Les graphes qui interviennent vont se factoriser comme dans \cite{To1, To2}. L'arête issue de la position $e^Y$ va sur le sommet d'un graphe représentant un élément de $\k$. Notons $G_s^{(1)}(X, Y)$ les contributions de tous ces graphes; c'est une $1$-forme en $s$.\\

Il vient alors les deux équations d'évolution suivante :

\begin{equation}\nonumber
\mathrm{d}_s Z_{sym,s}^{(1)}(X, Y)= [Y, G_s(X, Y)]\cdot \partial_Y Z_{sym,s}^{(1)}(X, Y)
\end{equation}
\begin{equation}\nonumber
\mathrm{d}_s D_s^{(1)}(X, Y)= [Y, G_s(X, Y)]\cdot \partial_Y D_s^{(1)}(X, Y) + \tr_\k \Big(\partial_Y G_s(X, Y) \circ \ad Y\Big) D_s^{(1)}(X, Y),
\end{equation}
avec condition limite pour $s=0$ : $$Z_{sym,s=0}^{(1)}(X, Y)= Z_{sym}(X, Y)$$ et $$D_{s=0}^{(1)}(X, Y)= \frac{J^{1/2}(X)J^{1/2}(Y)}{J^{1/2}(Z_{sym}(X, Y))}=D_{sym}(X, Y).$$
La condition limite pour $s=\infty$ correspond à $e^X\star_1 e^Y$.
\subsubsection{Deuxième déformation}

Positionnons maintenant $e^X$ et $e^Y$ sur l'axe horizontal. Sans perte de généralité, on peut positionner $e^Y$ en $1$ et $e^X$ en $u\in ]0,1[$. Lorsque $u$ tend vers $0$ on trouve $e^X\star_1 e^Y$, la condition limite précédente.
Quand $u$ tend vers $1$ on trouve $$E(X,Y) e^{X+Y},$$ le produit pour les paires symétriques. \\

Pour $u$ quelconque, les contributions totales sont de la forme $$D_u^{(2)}(X, Y) e^{Z_{sym,u}^{(2)}(X, Y)},$$ où
$D_u^{(2)}(X, Y)$ désigne une fonction de densité et $Z_{sym,u}^{(2)}(X, Y)\in \p$ est une déformation de la fonction de Campbell-Hausdorff. La déformation est contrôlée par des champs $\k$-adjoints agissant cette fois-ci sur les deux variables.  Les concentrations sur $e^X$ et $e^Y$ vont définir des $1$-formes en $u$ représentant  deux séries à valeurs dans $\k$,   $F_u(X, Y)$ (les contributions des graphes avec arête colorée par $\k^*$ issue de $e^X$) et $G_u(X, Y)$  (les contributions des graphes avec arête colorée par $\k^*$ issue de $e^Y$). \\

\noindent Il vient donc les équations suivantes :

\begin{equation}\nonumber
\mathrm{d}_u Z_{sym, u}^{(2)}(X, Y)=\left( [X, F_u(X, Y)]\cdot \partial_X + [Y, G_u(X, Y)]\cdot \partial_Y \right) Z_{sym,u}^{(2)}(X, Y)
\end{equation}
\begin{multline}\nonumber
\mathrm{d}_u D_u^{(2)}(X, Y)=  \left( [X, F_u(X, Y)]\cdot \partial_X + [Y, G_u(X, Y)]\cdot \partial_Y \right) D_u^{(2)}(X, Y)=\\
\tr_\k \Big(  \partial_X F_u(X, Y) \circ \ad X +  \partial_Y G_u(X, Y) \circ \ad X \Big) D_u^{(2)}(X, Y).
\end{multline}
\subsubsection{Déformation de la formule de Campbell-Hausdorff pour les paires symétriques}

La jonction des deux déformations précédentes est alors régulière et on en déduit le théorème suivant généralisant les résultats de \cite{To2}.

\begin{theo}\label{theoKVsym}

Il existe une déformation régulière $Z_{sym, v}(X, Y)$  de la fonction de Campbell-Hausdorff pour les paires symétriques  et une déformation régulière $D_{sym, v}(X, Y)$  de la fonction de densité pour les paires symétriques telles que :

\begin{multline}
\mathrm{d}_v Z_{sym, v}(X, Y)= \left( [X, F_v(X, Y)]\cdot \partial_X + [Y, G_v(X, Y)]\cdot \partial_Y \right) Z_{sym, v}(X, Y)\\
\mathrm{d}_v D_{sym, v}(X, Y)= \left( [X, F_v(X, Y)]\cdot \partial_X + [Y, G_v(X, Y)]\cdot \partial_Y \right) D_{sym, v}(X, Y)\\+ \tr_\k \Big(  \partial_X F_v(X, Y) \circ \ad X +  \partial_Y G_v(X, Y) \circ \ad X \Big) D_{sym, v}(X, Y),
\end{multline}
où $F_v$ et $G_v$ sont des  $1$-formes en $v$ et des séries à valeurs dans $\k$ convergentes dans un voisinage de $(0,0)$. Les conditions limites sont pour $v=0$ :
\begin{equation}\nonumber D_{sym, v=0}(X, Y)=\frac{J^{1/2}(X)J^{1/2}(Y)}{J^{1/2}(Z_{sym}(X, Y))} \quad  \quad Z_{sym, v=0}(X, Y)= Z_{sym}(X, Y)
\end{equation}
et pour $v=\infty$ :
\begin{equation}\nonumber \nonumber D_{sym, v=\infty}(X, Y)=E(X, Y)  \quad \quad   Z_{sym, v=\infty}(X, Y)=X+Y.
\end{equation}
\end{theo}
\vspace{0,5cm}

Cette déformation, le long des axes, remplace la déformation dans le demi-plan complexe définie par Kontsevich et utilisée dans \cite{To1, To2}. En utilisant les techniques développées par A. Alekseev et E. Meinrenken dans \cite{AM2} et expliquées de manière simplifiée dans \cite{To5}, on montre facilement  que la condition $E=1$ implique que la déformation par dilatation est aussi contrôlée par les champs $\k$-adjoints (avec la condition de trace vérifiée). C'est le cas des paires symétriques résolubles ou les algèbres de Lie quadratiques considérées comme des paires symétriques.\\

On en déduit la proposition suivante.
\begin{prop}\label{propKVsym} Pour les algèbres de Lie quadratiques vues comme des paires symétriques, la déformation le long des axes fournit une solution à la conjecture de Kashiwara-Vergne. Plus généralement la conjecture $E=1$ pour les algèbres de Lie vues comme des paires symétriques implique la conjecture de Kashiwara-Vergne.
\end{prop}

Considérons les solutions des équations différentielles   dans
le groupe $K$, d'algèbre de Lie $\k$ avec condition initiale
triviale pour $s=0$ :

\begin{center}

$\mathrm{d}_v a_v(X, Y)=F_v(X, Y) a_v(X, Y) \quad \mathrm{et} \quad \mathrm{d}_v
b_v(X, Y)=G_v(X, Y) b_v(X, Y).$

\end{center}
Comme dans \cite{Rou90, Rou91, Rou94}, la fonction $v\mapsto
D_{sym, v}(a_v(X, Y)\cdot X, b_v(X, Y)\cdot X)$ vérifie  une équation différentielle du
type $y'(v)=\mu(v) y(v)$ ce qui donne une formule pour $y(v)$.
\begin{prop} On a l'expression de la fonction $E(X,Y)$ suivante

\begin{multline}\nonumber E(X,Y)\left(\frac{J^{1/2}(X)J^{1/2}(Y)}{J^{1/2}(Z_{sym}(X, Y))}\right)^{-1}=\\ \exp\Bigg(
\int_0^{\infty} \mathrm{tr}_\k\Big( \partial_X F_v(a_v\cdot X, b_v\cdot
Y) \circ \ad (a_v\cdot X) +  \partial_Y G_v(a_v\cdot X, b_v\cdot
Y) \circ \ad (b_v\cdot Y)\Big) \mathrm{d} v \Bigg).
\end{multline}
\end{prop}

\section[Homomorphisme d'Harish-Chandra]{Homomorphisme d'Harish-Chandra en termes de
graphes}\label{sectionHomomorphismeHC} On fixe une paire
symétrique $(\g, \sigma)$ et une décomposition d'Iwasawa (voir
\S~\ref{iwasawa})
\begin{equation}\g=\k\oplus\p_o\oplus \n_+.\end{equation}

\subsection{Diagramme d'Harish-Chandra et espaces de réduction}

On place dans le diagramme de Cattaneo-Felder pour la
bi-quantification (\S~\ref{soussectionbiquantification}), les
sous-variétés co-isotropes $\k^\perp$ en position horizontale et
$(\k_o\oplus
\n_+)^\perp$ en position verticale: c'est le  \textit{diagramme d'Harish-Chandra}.\\

On a la décomposition et les couleurs suivantes
\[\g=\underset{(-, +)}{\k/\k_o}\oplus \underset{(-, -)}{\k_o}\oplus
\underset {(+, +)}{\p_o}\oplus \underset{(+, -)}{\n_+},\] o\`u on a
noté $\k/\k_o$ un supplémentaire de $\k_o$ dans $\k$.

On fera attention que dans les graphes qui vont intervenir, il
peut y avoir une arête double si celle-ci est colorée par deux
couleurs différentes. \\

- (i) - L'algèbre de réduction  qui est placée verticalement est
tout simplement d'après la  Proposition~\ref{propreductioniwasawa}
 (\S~\ref{sectionespacereductioniwasawa}) :

$$H^0_\epsilon((\k_o\oplus
\n^+)^\perp)=\mathcal{C}_{poly}(\k_o^\perp)^{\k_o}[[\epsilon]]
=S(\p_o)^{\k_o}[[\epsilon]].$$

- (ii) - L'algèbre de réduction placée à l'origine est aussi $S(\p_o)^{\k_o}[[\epsilon]].$

- (iii)- L'algèbre de réduction, qui est en position horizontale,
correspond à l'espace de réduction de

$$H^0_{\epsilon, (\p_o\oplus \n_+)}(\k^\perp).$$

Celui-ci dépend du choix du supplémentaire de $\k$, (ici $\p_o\oplus
\n_+$) mais  est isomorphe  à l'espace de réduction standard
($S(\p)^\k, \underset{CF}\star)$ grâce au
Théorème~\ref{theodependance} (\S~\ref{dependance}).

\subsection{Espace de réduction vertical}
Dans la partie verticale on trouve l'espace de réduction
$S(\p_o)^{\k_o}$. Examinons le star-produit construit en Proposition \ref{propreductioniwasawa} (\S\ref{sectionespacereductioniwasawa}). Il est clairement donné
par une fonction de type $\widehat{\,E\,}(X, Y)$, correspondant aux
graphes de type Roue attachée à $e^X, e^Y$ pour $X, Y \in \p_o$.

\begin{prop}

On a $\widehat{\,E\, }(X, Y)=E_{\g_o}(X, Y),$ où  $E_{\g_o}(X, Y)$
est la fonction $E$ pour la petite paire symétrique $(\g_o, \sigma)$. Le
produit, dans la petite paire symétrique est donc le produit de
Rouvière standard.
\end{prop}

\noindent \textit{\bf Preuve :} Examinons les contributions qui
interviennent dans une  roue  attachée à des éléments de $\g_o$. Les
couleurs dans le cycle sont alors simples :\\

- (i)  soit le cycle dérive dans les directions de $\g_o^*$, et on
retrouve la fonction $E_{\g_o}$ pour la petite paire symétrique,\\

- (ii) soit le cycle dérive dans les directions de $\n_+$ avec
couleur~$\dashrightarrow$,\\

- (iii) soit le cycle  dérive dans les directions de $\g/(\g_o\oplus
\n_+)$ avec
couleur~$\longrightarrow$. \\

On a donc des contributions du genre
\begin{equation}\label{51}
\overset{\dashrightarrow}w_{\G_1\G_2\ldots \G_l}
\tr_{\n_+}\Big(\ad (\G_1)\circ \ad \G_2 \ldots \circ \ad
(\G_l)\Big)\end{equation} ou bien
\begin{multline}\label{52}
\overset{\longrightarrow}w_{\G_1\G_2\ldots \G_l} \tr_{\g/(\g_o\oplus
\n_+)}\Big(\ad (\G_1)\circ \ad \G_2 \ldots \circ \ad
(\G_l)\Big)=\\\overset{\longrightarrow}w_{\G_1\G_2\ldots \G_l}
\tr_{\n_-}\Big(\ad (\G_1)\circ \ad \G_2 \ldots \circ \ad (\G_l)\Big).
\end{multline}

 Remarquons alors que nous avons la symétrie $\overset{\dashrightarrow}w_{\G_1\G_2\ldots
\G_l}=(-1)^{l-1}\overset{\longrightarrow}w_{\G_1\G_l\ldots \G_2}$.\\

Notons $V(X, Y)$ toutes les contributions de type (ii) ou (iii). Ici
il se passe un phénomène remarquable. Comme la couleur est homogène dans le cycle, les dérivations se font soit toutes dans $\n_+$
soit toutes dans $\k/\k_o$. L'argument de déformation à $4$ points (\cite{CKT} \S 8.5.2),
montre qu'il n'y a pas de trace dans l'équation de déformation. Par
conséquent on peut faire intervenir la compensation dans $K_o$. On
obtient au final :
\begin{equation}\label{53} V(X, Y)^2 D(X+Y)=D(a\cdot X)D(b\cdot Y),
\end{equation}
  où $D(X)$ vaut $V(X, -X)$ et  $a, b$ sont dans $K_o$. Comme $V$ est clairement
invariante pour l'action de $K_o$, on obtient la factorisation
souhaitée. De plus $D$ est symétrique et donc les roues qui
interviennent sont de taille paire et attachées directement à l'axe
réel. Or, pour $X\in \p_o$ on a  $\tr_{\n_+}(\ad X)^{2n}=\tr_{\n_-}(\ad X)^{2n}$; des relations  (\ref{51}), (\ref{52}), (\ref{53}) on déduit que l'on a $D(X)=1$.
\fin

\paragraph{Remarque 9: } Pour les
petites paires symétriques, d'après \cite{To4},  le produit de
Rouvière correspond au produit standard de l'algèbre
symétrique via la symétrisation. c'est-à-dire que la fonction
$E_{\g_o}$ est sans effet sur les éléments $\k_o$-invariants.

\subsection{Espace de réduction horizontal}
\subsubsection{Equations de réduction} Les équations qui décrivent l'espace  de réduction
horizontal $H^0_{\epsilon, (\p_o\oplus \n_+)}(\k^\perp),$ sont données par des
graphes colorés avec un seul point terrestre. Ce sont les
graphes de la section \ref{reductionlineaire}. Les équations portent un degré (le nombre de sommets aériens) et les
termes de degré pair sont nuls d'après le
Lemme~\ref{lemmereductionlineaire} (\S~\ref{reductionlineaire}). L'opérateur $D_1$ correspond à l'action des champs $\k$-adjoints et l'opérateur $D_3$ contient un terme de type Bernoulli qui est non nul. Le dessin de gauche de (\ref{dessintransmutation}) illustre un graphe intervenant dans $D_5$.\\

Ce sont des équations compliquées. En effet, les graphes qui
apparaissent dans la détermination de l'espace de réduction sont de
type Bernoulli ou de type Roue avec une sortie vers un Bernoulli.
Les arêtes sortantes sont dans $\k^*$, tandis que les arêtes
 qui arrivent sur l'axe réel sont  colorés par $\p_o^*$ ou $\n_+^*$.
 Les arêtes intermédiaires peuvent prendre toutes les couleurs,
 ce qui ne permet pas des simplifications dans le calcul des
 coefficients ou de l'opérateur.\\

\paragraph {Explications :} L'espace de réduction pour ce choix de supplémentaire est
heuristiquement, l'écriture dans les coordonnées $(\n_+, \p_o)$ des
éléments $\k$-invariants de $U(\g)/U(\g)\cdot \k $ via la
décomposition
\[U(\g)/U(\g)\cdot \k \simeq U(\n_+) \otimes U(\g_o)/U(\g_o)\cdot
\k_o \simeq S(\n_+)\otimes S(\p_o),\] o\`u on a utilisé dans chaque
facteur  la symétrisation. On conçoit que l'écriture des éléments
$\k$-invariant ne soit  pas simple!\\

\begin{equation}\label{dessintransmutation}
\begin{array}{ccc}
  \includegraphics[width=6cm]{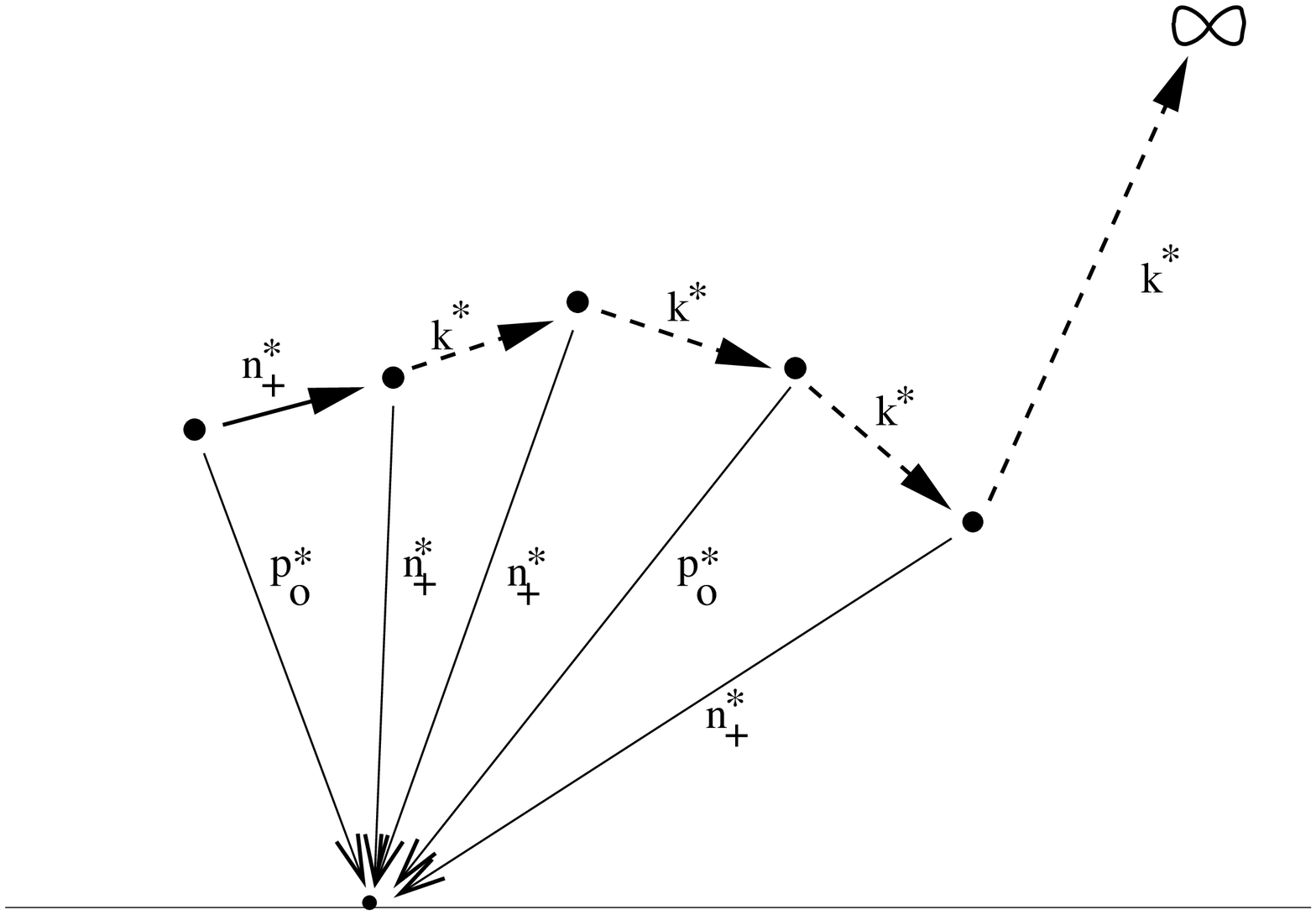} & \quad & \includegraphics[width=6cm]{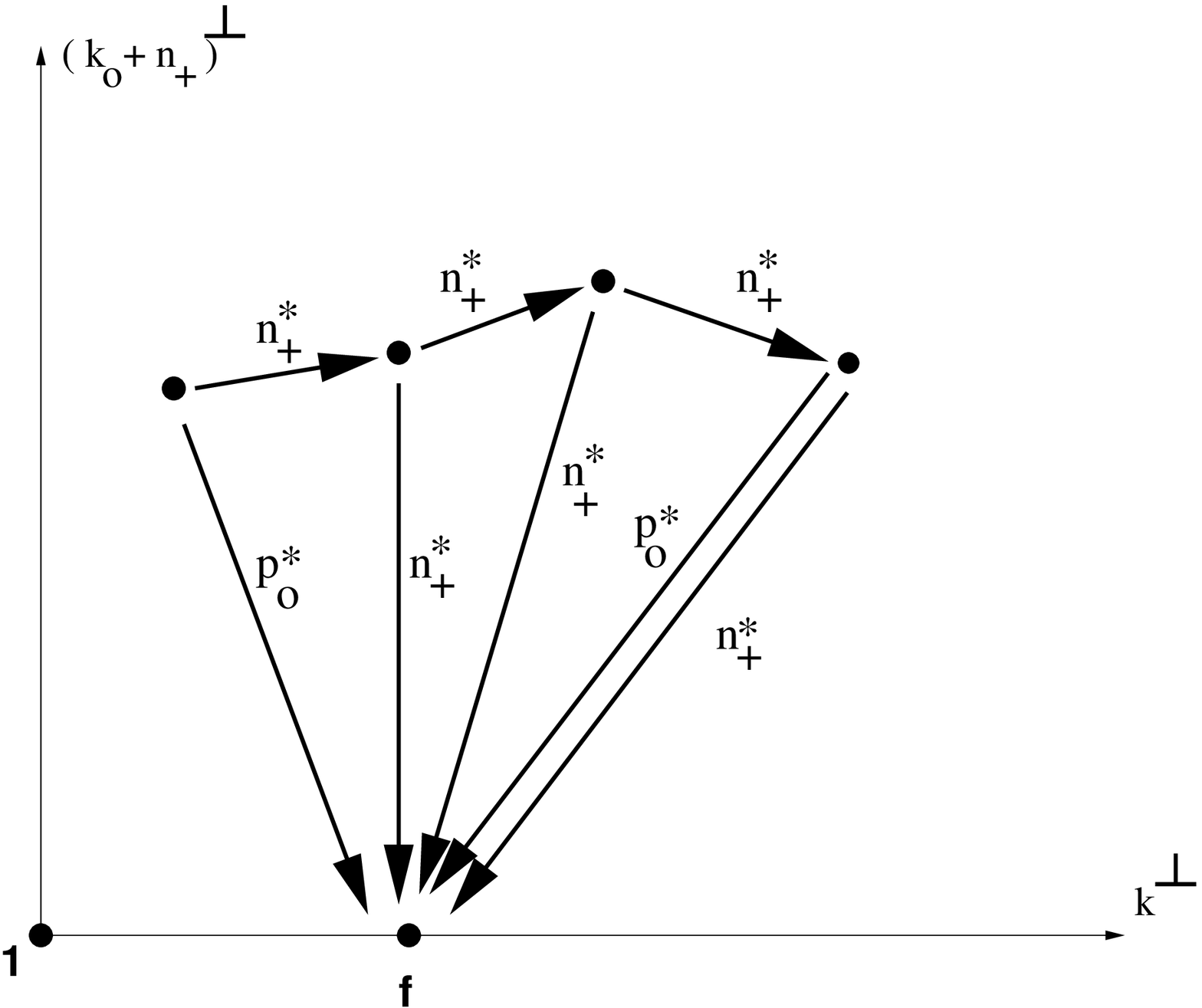} \\
  Graphe\; de \; D_5 & \quad & Transmutation \; horizontale
\end{array}
\end{equation}

\subsubsection{Description du produit $\star_{\,\p_o + \n_+}$  }
Le produit est composé soit de graphes de type Lie, soit de type Roue.\\

i- Les graphes de type Lie ont forcément une racine dans $\n_+$.
En effet, comme les couleurs sur l'axe réels sont dans $\p_o$ ou
$\n_+$ les racines sont dans $\g_o$ ou $\n_+$. Si la racine est dans
$\g_o$, alors le coefficient est nul comme dans le cas des paires
symétriques.\\

ii- Les graphes de type Roue se décomposent en deux
sous-ensembles. Le raisonnement fait pour la description de l'espace
de réduction verticale peut se recopier. On retrouve les
contributions pour le petit espace symétrique, où les cycles sont  colorés
de manière uniforme par $\n_+^*$ (= $\rightarrow$) ou  par
$\k^*/\k_o^*$ (=$\dashrightarrow$). Pour calculer ces dernières
contributions, il suffit de placer des points génériques $e^X$ et
$e^Y$ sur l'axe réel avec $X,Y \in \p_o$. L'argument à $4$ points (\cite{CKT} \S 8.5.2)
montre encore que l'on dispose d'une équation d'évolution sans
trace, ce qui permet de compenser l'évolution par un élément de
$K_o$. Pour les raisons des symétries on va trouver que ces
contributions valent encore $1$.

\subsection{Entrelacement et projection}

Rappelons la  décomposition et les couleurs  dans le diagramme de
Cattaneo-Felder associé à la décomposition d'Iwasawa :
\[\g=\underset{(-, +)}{\k/\k_o}\oplus \underset{(-, -)}{\k_o}\oplus
\underset {(+, +)}{\p_o}\oplus
\underset{(+, -)}{\n_+}.\]\\

 On définit l'opérateur de transmutation
dans le diagramme $(\k_o\oplus \n_+)^\perp$ et $\k^\perp$, c'est-à-dire que l'on va utiliser la structure de bi-module de
$S(\p_o)^{\k_o}$ pour transmuter les éléments d'un coté vers
l'autre. Pour $f\in H_{\epsilon, \p_o+\n_+}(\k^\perp)$ on va écrire

\begin{equation}\Gamma_{HC}(f)\underset{2} \star 1= 1\underset{1}\star f,\end{equation}
où $\G_{HC}(f)$ est un polynôme invariant dans $ S(\p_o)^{\k_o}$.
L'opérateur $\Gamma_{HC}$ sera appelé opérateur de transmutation. \\

\begin{lem}Les opérateurs qui interviennent dans la transmutation pour le diagramme d'Harish-Chandra
sont des roues (\textit{a priori} tentaculaires).
\end{lem}
\noindent \textit{\textbf{Preuve :}} Pour la partie horizontale, les
dérivées se font dans les directions $\p_o^*$ ou $\n_+^*$. Donc pour
les graphes de type Lie, la racine est dans $\n_+$, ce qui va
donner $0$ pour l'opérateur quand on va le restreindre à
$(\k_o\oplus \n_+)^\perp\cap \k^\perp$. Donc toutes les arêtes qui
arrivent sur l'axe réel sont dans $\p_o^*$. Mais alors il y a une
arête double de même couleur et le coefficient est nul. \\

Pour la partie verticale c'est plus simple. Comme l'espace de
réduction est $S(\p_o)^{\k_o}$,  seules les arêtes dans $\p_o^*$
peuvent dériver. On a donc une arête double dans les graphes de type
Bernoulli.

 \fin\\

Notons ici $A_{HC}$ les contributions des roues horizontales et
$B_{HC}$ les contributions des roues  verticales.

\begin{lem} La roue $A_{HC}$ n'a que des arêtes sortantes dans $\p_o^*$.
Il est de même pour la roue $B_{HC}$. Les symboles correspondants
sont donc des fonctions exponentielles sur $\p_o$. Les roues
sont pures.
\end{lem}

\noindent \textit{\textbf{Preuve :}} Supposons qu'il y ait une
sortie dans $\n_+^*$ au sommet $n$. En regardant la couleur de
l'autre arête issue de $n$,  on se
convainc que les couleurs dans le cycle sont : \\

-(i) soit toutes dans $\n_+^*$,\\

-(ii) soit toutes dans un supplémentaire de $\g_o\oplus \n_+$ pris
dans $\k$.\\

 Alors les contributions donnent selon les cas :\\
\[(i)\quad \tr_{\n_+}(\ad P_1 \cdots \ad N_1 \cdots )=0\] ou
\[(ii) \quad \tr_{\g/(\g_o\oplus \n_+)}(\ad P_1 \cdots \ad N_1 \cdots
)=0.\]Pour les roues de  $A_{HC}$, on conclut que toutes les sorties  sont
dans $\p_o^*$ et par conséquent les arêtes  finales sont colorées
par $\p_o^*$; forcément les roues sont pures. C'est bien-s\^ur
aussi le cas pour les roues de  $B_{HC}$ car l'espace de réduction est
$S(\p_o)^\k$.

\fin

\begin{prop}
On a la relation  $B_{HC}^{-1}A_{HC}=1$. L'opérateur de
transmutation correspond simplement à la restriction à $(\k\oplus
\n_+)^\perp$.

\end{prop}

\noindent \textit{\textbf{Preuve :}} Les roues, dont le cycle est coloré par
$\p_o^*$ et $\k_o^*$ (c'est-à-dire $(+, +) $ ou $(-, -)$ ),  vérifient
la même propriété que celle que nous avons déjà rencontrées à la
Proposition~\ref{alternance} (\S~\ref{sectioncontributionsdansE}). Il
y a un nombre pair de sorties, avec alternance des couleurs $(+, +)
$ et $(-, -)$, donc comme dans la Proposition~\ref{alternance},  ces
contributions s'annulent aussi bien dans $A_{HC}$ que dans
$B_{HC}$, car les roues sont pures.\\

Il reste à considérer les roues, dont le cycle est de couleur uniforme en $(+, -)  $ (dérivée en
$\n_+^*$ ) ou de couleur uniforme  en $(-, +) $ (dérivée en $\k^*/\k_o^*$).\\

Notons $w^{A}_{\G_{ + \, -}}$ le coefficient associé à la roue pure
de taille $n$  attachée sur l'axe horizontal et de couleur $(+, -)$.
On  introduit de même  $w^{A}_{\G_{- +}}$, $w^{B}_{\G_{ + \, -}}$
et
$w^{B}_{\G_{ - +}}$.\\

Compte tenu des symétries (\ref{symetriesangle})  par rapport à la
bissectrice principale on~a :

\[w^{B}_{\G_{ +\,
-}}= (-1)^n w^{A}_{\G_{- +}} \quad \mathrm{et} \quad w^{B}_{\G_{ - +}}= (-1)^n
w^{A}_{\G_{ +\,  -}}.\] Par ailleurs du coté de l'opérateur, la
dérivation ne se faisant qu'en la direction de $\p_o^*$, le symbole
correspondant peut s'écrit pour $X\in \p_o$ :

\[ \G_{(+, -)}(X)=\tr_{\n_+} (\ad X)^n \]\[
\G_{(-, +)}(X)=\tr_{\g/(\g_o\oplus\n_+)} (\ad X)^n=  \tr_{\n_-}( \ad
X)^n= (-1)^n \tr_{\n_+} (\ad X)^n. \]

En effet  on a $\tr_{\n_+}\ad X =-\tr_{\n_-} \ad X$ pour $X\in
\p_o$. Au total il y a compensation et on a les relations suivantes:

\[w^{B}_{\G_{ +\,
-}}\tr_{\n_+} (\ad X)^n =  w^{A}_{\G_{- +}} \tr_{\n_-} (\ad X)^n=
w^{A}_{\G_{- +}} \G^{A}_{(+, -)}(X)\] et
\[w^{B}_{\G_{ - +}} \G_{(-, +)}(X)= w^{A}_{\G_{ +\,  -}} \G^{A}_{(+,
-)}(X).\]

Ceci montre que l'on a $A_{HC}=B_{HC}$. Comme on a
$\G_{HC}(f)=B^{-1}_{HC}A_{HC}(f)$, on en déduit que $\G_{HC}(f)$ est
bien la restriction de $f$ à $(\k\oplus \n_+)^\perp$.

 \fin\\

\begin{defi} En théorie de Lie pour les paires symétriques, la projection sur le facteur $U(\g_o)/U(\g_o)\cdot \k_o$ dans la décomposition
\[U(\g)/U(\g)\cdot \k \simeq U(\n_+) \otimes U(\g_o)/U(\g_o)\cdot
\k_o, \] s'appelle la  projection d'Harish-Chandra. C'est un
homomorphisme d'algèbres  pour les invariants de
\[(U(\g)/U(\g)\cdot \k)^\k\quad \mathrm{sur} \quad
(U(\g_o)/U(\g_o)\cdot \k_o)^{\k_o}.\]
\end{defi}

 On peut maintenant énoncer un résultat similaire en termes de
graphes de Kontsevich.

\begin{prop} La projection d'Harish-Chandra correspond à l'opérateur de transmutation c'est-à-dire  la restriction
à $\p_o^{*}=(\k\oplus \n_+)^\perp$. Elle prend ses valeurs dans
$S(\p_o)^{\k_o}$. C'est un  homomorphisme d'algèbres de
\[\left(H^0_{\epsilon, (\p_o\oplus \n_+)}(\k^\perp)\, ,\underset {
(\p_o\oplus \n_+)}\star\right)\quad \mathrm{sur} \quad
\left(S(\p_o)^{\k_o}, \underset{Rouv=CF}\star\right).\]
\end{prop}

\noindent \textit{\textbf{Preuve :}} Le fait que la restriction à
$\p_o^*$ soit une fonction $\k_o$-invariante est clair, car toutes
les opérations que nous avons effectuées peuvent être choisies
$\k_o$-équivariantes. En effet le choix du supplémentaire $\k/\k_o$
peut être facilement choisi $\k_o$-invariant, il suffit de prendre
$(\n_{-} \oplus \n_+)\cap \k$. Les équations sont clairement
indépendantes des choix de la base (une fois choisies les
décompositions), car cela correspond à des transformations linéaires
qui préservent les couleurs des graphes. Maintenant $\k_o$
préserve la décom\-position. Donc l'espace de réduction est
$K_o$-invariant.

\fin\\

\paragraph{Commentaire : } Il est bien connu que l'homomorphisme d'Harish-Chandra est
quelque chose de très compliqué. Ici il correspond à une simple
restriction. Il faut comprendre que la difficulté est cachée dans la
description de l'algèbre de réduction. En quelque sorte il faut
maintenant écrire un isomorphisme explicite entre $S(\p)^\k$ et
$H^0_{\epsilon, (\p_o\oplus \n_+)}(\k^\perp)$, pour décrypter toute
la difficulté de l'homomorphisme d'Harish-Chandra. C'est l'objet de
la section suivante.

\subsection{L'homomorphisme d'Harish-Chandra en termes de
gra\-phes}\label{HCdiagrammes}

Le Théorème~\ref{theodependance} (\S~\ref{soussectionentrelacement})
montre qu'il existe un isomorphisme entre les deux espaces de
réduction correspondants aux décomposition de Cartan et d'Iwasawa.
On disposera alors d'une formule pour l'homomorphisme
d'Harish-Chandra en termes de graphes. Il existe un isomorphisme d'algèbres de  $(S(\p)^\k,
\underset{CF}\star)$ dans $(S(\p_o)^{\k_o}, \underset{CF}\star)$ ;
c'est l'homomorphisme
d'Harish-Chandra.\\

\subsubsection{Calcul de l'élément de jauge}

On reprend ici les notations de \S~\ref{soussectionentrelacement}.
Le champ $v$ s'écrit $v=-V_a\partial_{P_a^*},$ où $(P_a^*)_a$ désigne une base
de $\k^\perp=\p^*$ et $V_a\in \k$.\\

L'opérateur qui réalise l'entrelacement est calculé de
manière récursive. Plus précisément on doit intégrer la composante de degré (total) $0$
 du champ  $[ DU_{\pi_t}(v), \bullet]_{GH}$. C'est en effet cet
 élément qui conjugue les différentiel\-les et réalise
 l'isomorphisme en cohomologie pour les star-produits, comme on l'a
 vu au Théorème~\ref{theodependance}.

 \begin{lem}

 L'élément du groupe
de jauge qui entrelace les star-produits pour deux choix de
supplémen\-taires, vivra dans le groupe associé à l'algèbre de Lie
engendrée par les coefficients de Taylor du champ $[
\left(DU_{\widehat{\, \pi_t}}(\,\widehat{v})\right)_0,
\bullet]_{GH}$.

\end{lem}

\noindent \textit{\textbf{Preuve :}} L'équation suivante sur les
composantes de degré~$0$ (dans le complexe de Hochschild) contrôle
la déformation des star-produits :

$$\frac{\partial (\mu_t)_0}{\partial t}=
[\left(DU_{\widehat{\pi_t}}(\widehat{v})\right)_0,(\mu_t)_0].$$

Pour simplifier, notons
$A(t)=[DU_{\widehat{\pi_t}}(\widehat{\,v})_0, \bullet]$ l'opérateur
linéaire d'ordre $0$ (total et de Hochschild) et
$y(t)=(\mu_t)_0$.\\

On doit donc résoudre une équation différentielle linéaire de la
forme $y'(t)=A(t)y(t)$. On a, comme série formelle en $t$ : $y(t)=y(0)+\sum_{n\geq 1} \frac{t^n}{n!}y_n$.
L'équation différentielle donne $y_1=A(0)y_0$. Posons
$y^{(1)}(t)=e^{-tA(t)}y(t)$. Alors $y^{(1)}(t)=y_0 {\pmod {t^2}}$ et
vérifie l'équation
\begin{equation}
\frac{\mathrm{d}}{\mathrm{d}t}y^{(1)}(t)=t \left( \frac{e^{-\ad(tA(t))}-1}{\ad(t A(t))}
A'(t)\right)y^{(1)}(t).
\end{equation}

On pose $A_1(t)=\left( \frac{e^{-\ad(tA(t))}-1}{\ad(t A(t))}
A'(t)\right),$
qui est dans l'algèbre de Lie engendrée par les
coefficients de Taylor de $A(t)$. Il vient l'équation différentielle
\[\frac{\mathrm{d}y^{(1)}(t)}{\mathrm{d}t}= tA_1(t)y^{(1)}(t).\]

  On
pose ensuite $y^{(2)}(t)=e^{-\frac{t^2}{2}A_1(t)} y^{(1)}(t)=e^{-\frac{t^2}{2}A_1(t)}e^{-tA(t)}y(t).$ Alors on a
$y^{(2)}(t)=y_0 {\pmod {t^3}}$ et on a, via la formule de
Campbell-Hausdorff,

\[y^{(2)}(t)=e^{\left(-\frac{t^2}{2}A_1(t) - t A(t) + \frac12[\frac{t^2}{2}A_1(t), tA(t)]
\ldots \right) } \, y(t).\]

Par récurrence on
trouve que $y(t)$ est l'image
 (comme série formelle en $t$) de $y(0)$ par l'action d'un élément  du groupe
 formel
associé à  l'algèbre de Lie des coefficients de Taylor
 de $A(t)$. On a de plus comme série formelle en $t$
 \[y(t)=e^{\Omega(t)} y(0)=e^{tA(t)}e^{\frac{t^2}{2}A_2(t)}\cdots y(0).\]
La résolution sous cette forme est connue en analyse numérique
\cite{iserles2} et
se trouve explicitement dans
 \cite{magnus}.\footnote{Référence que nous a indiquée D. Manchon et qui réfère à Zassenhaus.}

\fin

\subsubsection{Etude de l'opérateur $(DU_{\pi}(v))_0$}
On a $\pi_{t=0}=\pi$. Dans le cas des paires
symétriques, on montre maintenant que
l'opérateur $DU_{\pi_t}(v)_0$ est nul pour $t=0$.

\begin{lem}\label{lemmetaylor} Les graphes qui interviennent dans la composante de Hochschild de
degré $0$ de $DU_{\pi}(v)$ que l'on note $(DU_{\pi}(v))_0$ sont de
trois types : des graphes de type Bernoulli fermés par $v$ (\ref{Bernoulliv}) (dessin de gauche)
ou des graphes de type  Roue attachée à un
Bernoulli lequel s'attache à $v$ (\ref{Bernoulliv}) (dessin du milieu) ou des
graphes de type Roue pure attachée à $v$ (\ref{Bernoulliv}) (dessin de droite).
L'opérateur $DU_{\pi}(v)_0$  est nul mais
$(DU_{\pi}(v))_1$ n'est pas nul.
\end{lem}

\noindent \textit{\textbf{Preuve :}} Il suffit de regarder ce qui se
passe lorsqu'on applique l'opérateur sur des fonctions, car les
dérivées sortant de l'axe réel n'interviennent pas, sinon il
sortirait une arête d'un sommet aérien ce qui
est exclu (on aura alors $n\geq 2$). \\

L'opérateur n'agit donc que sur les coefficients de la fonction $f$,
placée au point terrestre. Le sommet o\`u on a placé le champ de
vecteurs $v$ doit être dérivé dans la direction de $\k^*$
(car les coefficients du champ $v$ sont dans $\k$). Cette arête provient :\\

-- soit d'un graphe de type Bernoulli lui même attaché ou non à
une roue,

-- soit d'une roue pure.\\

-- Dans le premier cas, l'arête issue de $v$, peut soit dériver la
racine du graphe de Bernoulli\footnote{Ce sommet doit être dérivé
car il est dans $\k$.}(\ref{Bernoulliv}) (dessin de gauche), soit dériver le sommet
terrestre si le graphe de Bernoulli est attaché à une roue
(\ref{Bernoulliv}) (dessin du milieu).

-- Dans le second cas, l'arête issue de $v$ doit dériver le sommet
sur l'axe réel (\ref{Bernoulliv}) (dessin de droite).\\

Il est facile de se convaincre en examinant tous les cas, qu'il y a
forcément un nombre $n$ impair  de sommets attachés au bi-vecteur
$\pi$. Par conséquent, la symétrie par rapport à l'axe vertical fait
apparaître un signe $(-1)^n$ dans son coefficient, ce qui montre les
coefficients de ces graphes sont nuls. On a
$(DU_\pi(v))_0=0$.\\

Le terme de degré $1$ (pour le complexe de Hochschild) de
l'opérateur $DU_\pi(v)=0$ n'est pas nul car le graphe avec un seul sommet aérien et deux points terrestres intervient de manière non triviale (le
coefficient vaut $\frac 12$).
\fin

\begin{equation}\label{Bernoulliv}
\begin{array}{lcc}
  \includegraphics[width=5cm]{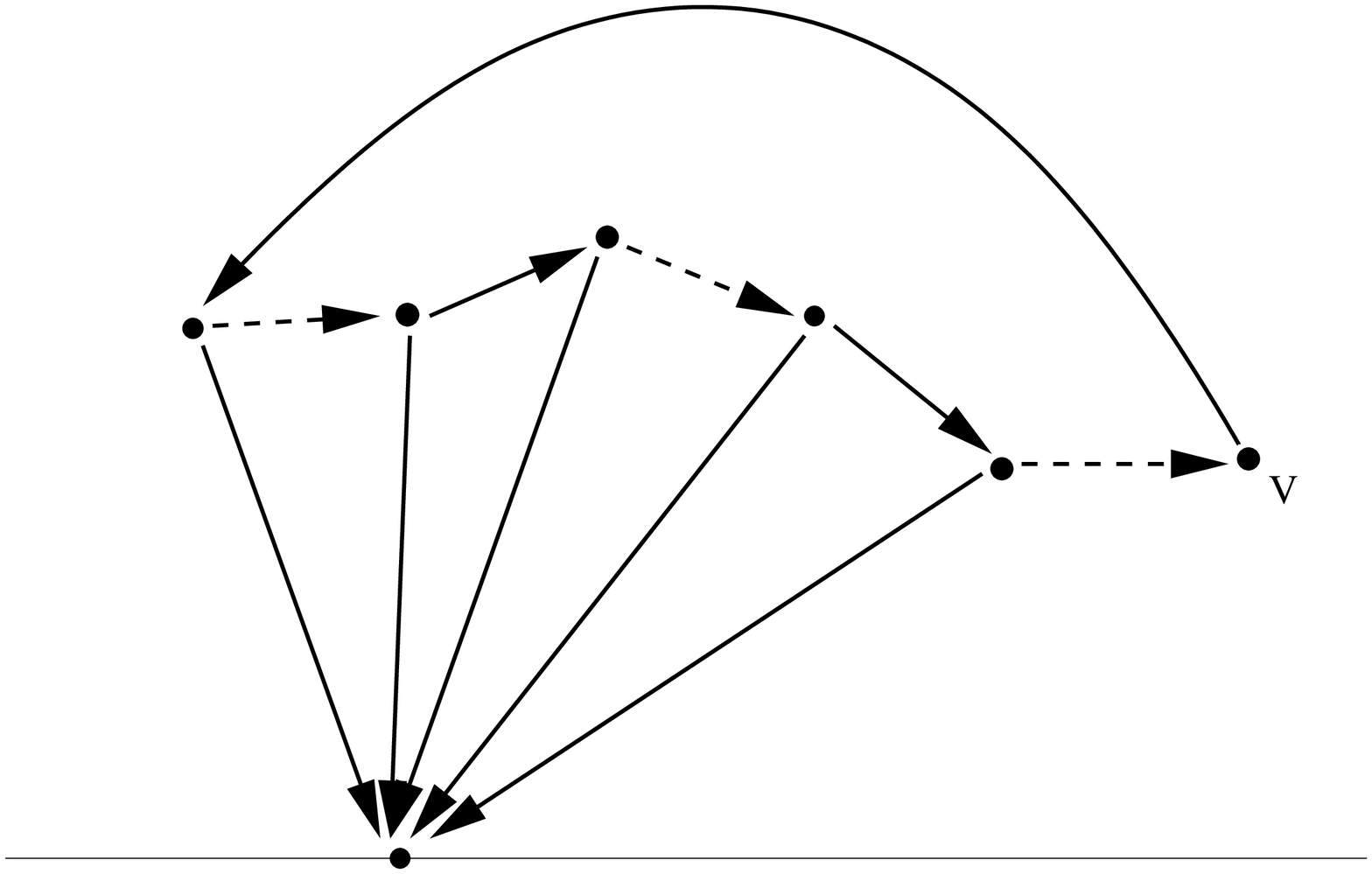} & \includegraphics[width=5cm]{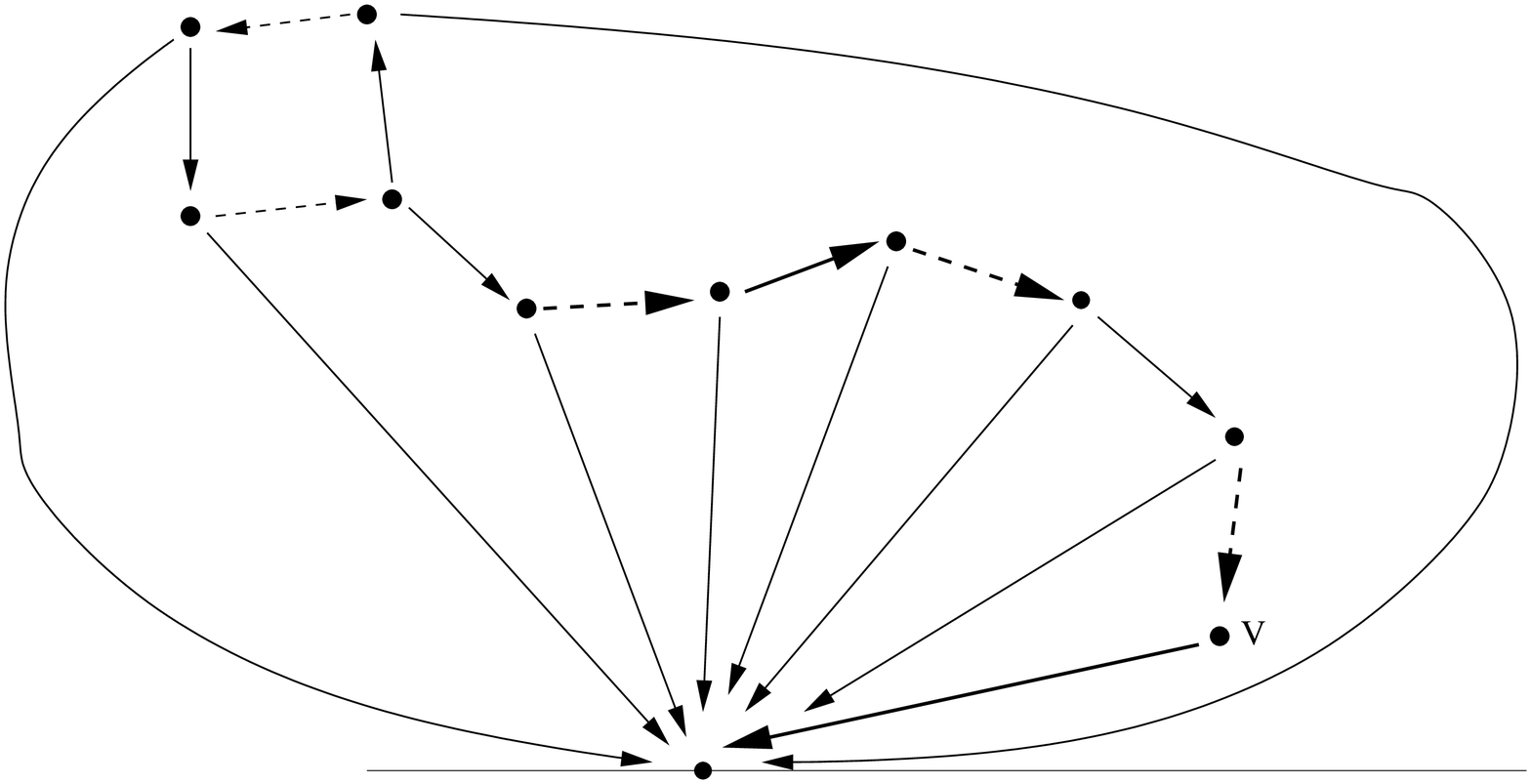} & \includegraphics[width=5cm]{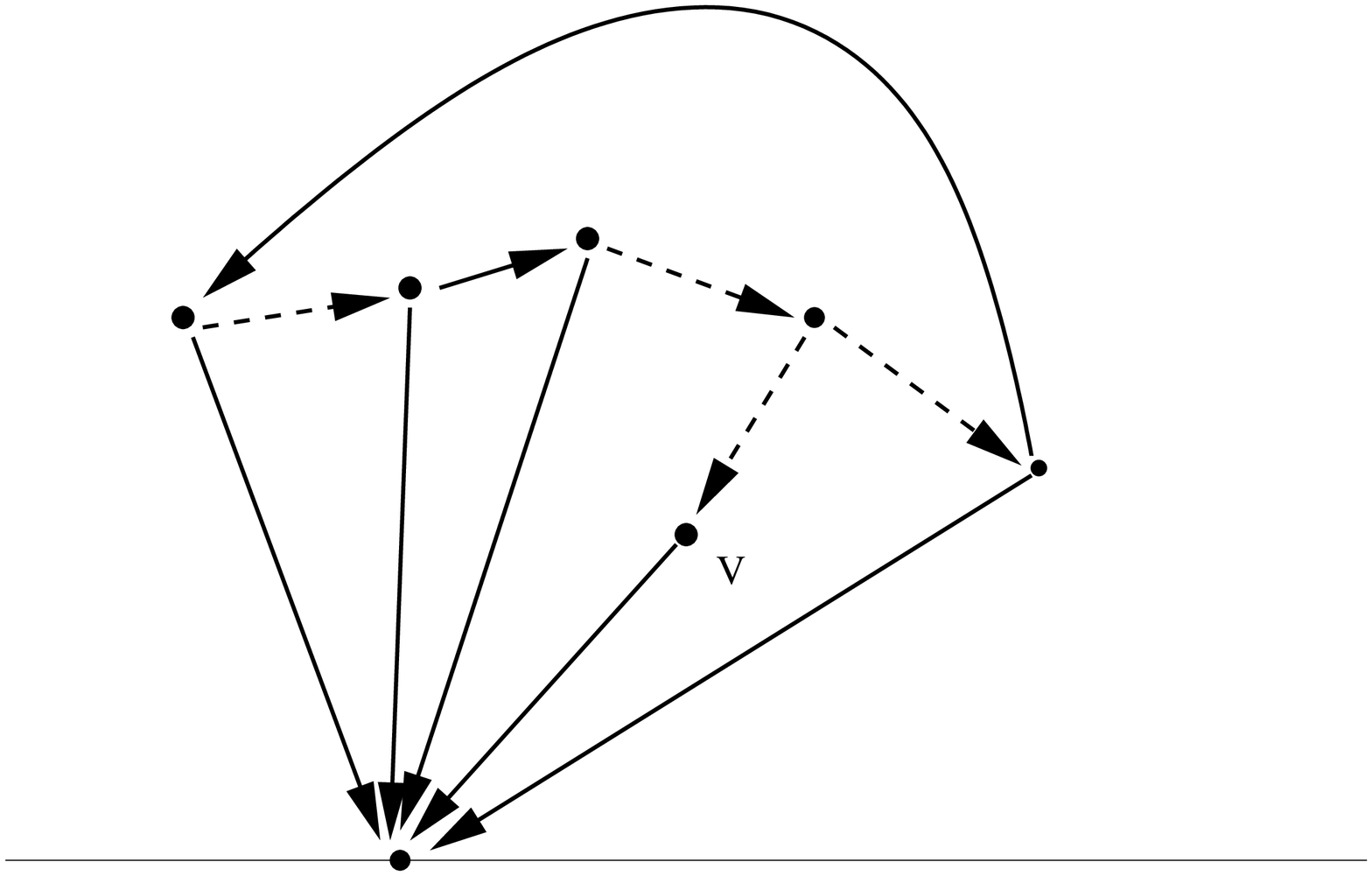} \\
 Bernoulli\; ferm\acute{e} \;  par \; v& Roue-Bernoulli \; attach\acute{e} \; \grave{a}\; v & Roue\; pure\; sur \; v \end{array}
\end{equation}

\subsubsection[Coefficients de Taylor]{Coefficients de Taylor et  invariance  par le groupe de Weyl généralisé} La
description diagrammatique de l'homomorphisme d'Harish-Chandra ne
 fait pas intervenir explicitement le groupe de Weyl généralisé.\\

Il faut  montrer que nos constructions sont indépendantes du choix
de $\n_+$ (choix d'un système de racines positives) et cela revient
essentiellement à montrer que ces constructions sont invariantes
lorsqu'on change $v$ en $-v$ (voir \cite{To3}).

\begin{prop}

Les coefficients de Taylor du champ $(DU_{\pi_t}(v))_0$ sont
invariants par le changement de champ $v \mapsto -v$.
L'homomorphisme d'Harish-Chandra ne dépend donc pas du choix de la
chambre de Weyl.
\end{prop}

\noindent \textit{\textbf{Preuve :}} On fait un inventaire des
graphes pour les coefficients de Taylor et on calcule le nombre
de fois que $v$ apparaît. \\

D'après les formules de \S~\ref{actionduchamp}, les bi-vecteurs $\pi$
et $[v,[v,\pi]]$ se comportent  de la même manière  vis à vis des\
couleurs, tandis que le
bi-vecteur  $[v, \pi]$  se comporte de manière opposée.\\

Les graphes, qui apparaissent dans le calcul des coefficients de
Taylor de $(DU_{\pi_t}(v))_0$ sont les graphes décrits dans le
Lemme~\ref{lemmetaylor} auquel il faut ajouter le graphe de type
Bernoulli non fermé (l'arête issue de $v$  dérive le sommet
terrestre et la racine du Bernoulli est dans $\p$). Les sommets, en dehors de $v$, sont
attachés aux bi-vecteurs $\pi, [v, \pi]$ ou  $[v, [v, \pi]]$.\\

Chaque sommet attaché à $\pi$ ou $[v,[v, \pi]]$ (on note $p$ leur
nombre) engendre un changement de couleur dans le cycle de la roue ou le brin de
Bernoulli, tandis que les sommets attachés à  $[v, \pi]$ (on note
$q$ leur nombre) ne change pas la couleur. On doit avoir $p+q+1$
impair pour
que le coefficient ne soit pas nul. Le nombre de fois que $v$ apparaît est congru à $q+1$.\\

i- Examinons le cas des graphes de type Bernoulli (fermé ou non):

Le nombre $p$ doit être impair. En effet si le graphe est de type
Bernoulli fermé alors l'arête issue de $v$ dérive la racine dans
la couleur $\p^*$. Si le graphe est de type Bernoulli non fermé les
arêtes issues de la racine sont colorées
par $\p^*$. Dans tous les cas $q+1$ est pair. Le nombre de fois que $v$ apparaît est donc pair.\\

ii- Examinons le cas des graphes de type Roue attachée à un Bernoulli :

On note $p_1$ (resp. $q_1$) le nombre de sommets, dans la roue,
attachés
 à $\pi$ ou $[v,[v, \pi]]$ (resp. $[v, \pi]$).
 On note $p_2$ (resp. $q_2$) le nombre de sommets, dans brin de Bernoulli, attachés
 à $\pi$ ou $[v,[v, \pi]]$ (resp. $[v, \pi]$).  Le nombre de changements de couleur dans le cycle
 est pair. On en déduit facilement que $p_1+p_2$ est impair. Donc $q_1+q_2+1$ est pair.

\fin

Rappelons que pour les paires symétriques, la décomposition
d'Iwasawa généralisée provient de l'action d'un tore $\s_f$
généralisant les sous-espaces de Cartan, c'est le sous-espace de
Cartan-Duflo (voir \S~\ref{sectioniwasawa}). Le groupe fini qui
remplace le groupe de Weyl est le groupe quotient
$N_K(\s_f)/Z_K(\s_f)$ du normalisateur sur le centralisateur.
Remplacer $v$ par $-v$ revient dans notre situation à remplacer
$\n_+$ par $\n_-$. On en déduit en faisant intervenir des
projections d'Harish-Chandra partielles que notre construction est
indépendante du choix du système de racines positives.\\

 On en déduit  le corollaire suivant:

\begin{cor}
L'homomorphisme d'Harish-Chandra est indépendant du choix du $\n_+$.
Il est invariant par le groupe de Weyl généralisé.

\end{cor}

\section{Construction de caractères}\label{sectionconstructioncaractere}
 On considère $(\g, \sigma)$ une paire symétrique et $\g=\k\oplus
\p$ sa décomposition de Cartan.\\

Dans cette section, on construit de manière systématique des
caractères pour l'algèbre des opérateurs différentiels invariants en
utilisant la  bi-quantifi\-cation (\textit{cf.}
\S\ref{soussectionbiquantification}) et des polarisations. Plus
précisément  on  applique le principe de bi-quanti\-fication aux
sous-variétés co-isotropes  $f +\b^\perp$ et $\k^\perp$ avec
 $f \in \k^\perp$ et $\b$ une polarisation en $f$. On renvoie au
  \S~\ref{sectioncasdespolarisations} pour la notion de polarisation.

\subsection[Construction de caractères]{Construction des caractères
pour des algèbres d'opéra\-teurs différen\-tiels invariants}

\paragraph{Principe de construction :} On  applique le principe de
bi-quantification aux sous-variétés co-isotropes  $f +\b^\perp$ et
$\k^\perp$ avec
 $f \in \k^\perp$ et $\b$ une polarisation en $f$ et on utilise la
 transmutation. D'après la Proposition~\ref{propreductionpolarisation} (\S~\ref{sectioncasdespolarisations})
le diagramme de bi-quantification fournit un caractère pour
l'algèbre de réduction~$(H^0_{\epsilon,\b}(\k^\perp),
\underset{CF}\star)$.
\\

En effet pour $P\in H^0_{\epsilon,\b}(\k^\perp)$ on calcule alors
$1\underset{1}\star P$. C'est une fonction polynomiale $K\cap B$~-invariante sur
$f+(\k+\b)^\perp$, où  $K$ et $B$ sont des groupes connexes d'algèbres
de Lie $\k$ et $\b$. Or $(K\cap B)\cdot f$ est ouvert dans
$f+(\k+\b)^\perp$, donc cette fonction polynomiale est constante.
L'application $P\mapsto 1\underset{1}\star P \in \R$ est donc un
caractère de l'algèbre associative~$H^0_{\epsilon,\b}(\k^\perp)$.\\

Cette algèbre de réduction dépend d'un choix de supplémentaire de
$\k$ dans $\g$, construit de manière compatible avec $\b$;  c'est-à-dire on doit choisir un supplémentaire  de $\b\cap \k$ dans $\b$
(que l'on pourrait noter abusivement $\b/(\b\cap \k)$), puis un
supplémentaire de $\k+\b $ dans $\g$ (que l'on pourrait noter
$\g/(\k+\b)$).\\

\begin{prop}\label{propconstructioncaractere} Le diagramme de Cattaneo-Felder appliqué au cas $f +\b^\perp$ et $\k^\perp$ avec
 $f \in \k^\perp$ et $\b$ une polarisation en $f$, fournit un caractère de l'algèbre
 de réduction $H^0_{\epsilon,\b}(\k^\perp)$.
 \end{prop}

Grâce au Théorème~\ref{theodependance}
(\S~\ref{soussectionentrelacement}), il existe un isomorphisme
d'algèbres entre $H^0_\epsilon(\k^\perp)$ et
$H^0_{\epsilon,\b}(\k^\perp)$. En composant par   cet isomorphisme,
on construit donc un caractère pour l'algèbre de réduction qui
dépend \textit{a priori} du choix de la polarisation $\b$
en~$f$.\footnote{On peut se demander si ce caractère est indépendant
du choix du supplémen\-taire de $\k+\b$ dans $\g$.}

\paragraph{Remarque 10 (importante) : }  Dans le cas d'une sous-algèbre $\h$,
cette méthode fournira un
caractère de l'algèbre de réduction $H^0_\epsilon(\h^\perp)$, pour
peu que l'on ait $(H\cap B) \cdot f$ ouvert dans $f+ (\h+\b)^\perp$.
Dans ce cas $H\cdot f$ est lagrangien dans $G\cdot f$ et on retrouve
une condition évoquée dans \cite{duflo-japon}, pour la commutativité
de l'algèbre l'algèbre des opérateurs différentiels invariants sur
$G/H$.

\subsection{Indépendance par rapport aux choix de la polarisation }
La méthode consiste à considérer une variation à
\textbf{$8$ couleurs} des graphes. On  retrouve aux bords les
graphes à $4$ couleurs, ce qui permet l'interpolation entre les
deux situations.

\subsubsection{Construction de la forme à $8$ couleurs}

 \begin{prop}
Il existe une $1$-forme à $8$ couleurs  qui interpole  les
$1$-formes à $4$ couleurs.
\end{prop}

\noindent \textit{\textbf{Preuve :}} Le reste de ce paragraphe  est
consacré à la construction de cette
$1$-forme à $8$ couleurs.\\

Considérons la demie-bande
\[
\hs:=\left\{\bx=(x,y)\in\bbR^2 : x\ge 0,\ -\hpi\le y\le
\hpi\right\}.
\]
Notons  $\hsb_i$, $i=1,2,3$, les trois composantes de son bord :
\begin{align*}
\hsb_1&\, :=\left\{\bx\in\hs : y=-\hpi\right\},\\
\hsb_2&\, :=\{\bx\in\hs : x=0\},\\
\hsb_3&\, :=\left\{\bx\in\hs : y=\hpi\right\}.
\end{align*}
Sur la compactification de l'espace de configurations  $C_2(\hs)$,
on considère l'involution
\[
p\colon\begin{array}[t]{ccc}
C_2(\hs) &\to &C_2(\hs)\\
(\bx_1,\bx_2)&\mapsto &(\bx_2,\bx_1).
\end{array}
\]
On veut définir huit  $1$\ndash formes fermées
$\theta_{j_1j_2j_3}$, $j_i\in\{1,2\}$, sur $C_2(\hs)$  ayant les
propriétés suivantes :
\begin{enumerate}
\item $\theta_{j_1j_2j_3}$ s'annulent quand $\bx_{j_i}$ approche de  $\hsb_i$, $i=1,2,3$;
\item sur la composante de  bord correspondant  au rapprochement de  $\bx_1$ et $\bx_2$  tout en restant
à l'intérieur de  $\hs$, $\theta_{j_1j_2j_3}$ s'identifie à la
forme volume normalisée et invariante sur $\mathbb{S}^1$;
\item pour  $j_i=1$ ($j_i=2$), $i=1,2,3$,
sur la composante de bord correspondant  au rapprochement de $\bx_1$
et $\bx_2$ près de l'intérieur de  $\hsb_i$,
$\theta_{j_1j_2j_3}$ ($p^*\theta_{j_1j_2j_3}$) est la $1$-forme de
Kontsevich;
\item sur le bord correspondant  au rapprochement de
 $\bx_1$ et $\bx_2$ près d'un coin,
$\theta_{j_1j_2j_3}$ est la $1$\ndash forme à $4$ couleurs de
Cattaneo-Felder~(\S~\ref{defi4couleurs}).
\end{enumerate}

\paragraph {Cas $j_1\not=j_3$ : } Soit $\theta$ la forme tautologique
sur $C_2(\bbR^2)$ correspondant à la métrique Euclidienne ;
c'est à dire, $\theta:=\phi^*\omega$, o\`u $\omega$ est la
$1$-forme volume normalisée et invariante sur $\mathbb{S}^1$ et
\[
\phi\colon\begin{array}[t]{ccc}
C_2(\bbR^2) &\to &\mathbb{S}^1\\
(\bx_1,\bx_2)&\mapsto &\frac{\bx_2-\bx_1}{||\bx_2-\bx_1||}
\end{array}
\]

Soit $\tau_i\colon\bbR^2\to\bbR^2$ ($i=1,2,3$) l'involution
correspondant à la réflexion par rapport aux droites supportant
$\hsb_i$; c'est à dire :
\begin{align*}
\tau_1(x,y) &= (x,-\pi-y),\\
\tau_2(x,y) &= (-x,y),\\
\tau_3(x,y) &= (x,\pi-y).
\end{align*}
Remarquons que l'on a  $\tau_1\circ\tau_2=\tau_2\circ\tau_1$ et
$\tau_3\circ\tau_2=\tau_2\circ\tau_3$.

Soit $C_2(\bbR^2)'$ le sous-ensemble de  $C_2(\bbR^2)$ défini par
 $\bx_1\not=\tau_i(\bx_2)$ pour $i=1,2,3$. Observons que l'on a  $C_2(\hs)\subset
C_2(\bbR^2)'$. Écrivons   $\tau_{ij}$ ($i=1,2,3$; $j=1,2$) pour
l'action de $\tau_i$ sur la  composante $j$ de  $C_2(\bbR^2)'$;
c'est à dire:
\begin{gather*}
\tau_{i1}\colon\begin{array}[t]{ccc}
C_2(\bbR^2)'&\to &C_2(\bbR^2)'\\
(\bx_1,\bx_2)&\mapsto &(\tau_i(\bx_1),\bx_2)
\end{array}\\
\tau_{i2}\colon\begin{array}[t]{ccc}
C_2(\bbR^2)'&\to &C_2(\bbR^2)'\\
(\bx_1,\bx_2)&\mapsto &(\bx_1,\tau_i(\bx_2))
\end{array}
\end{gather*}
Remarquons que l'on a
\begin{align*}
\tau_{1j_1}\circ\tau_{2j_2} &= \tau_{2j_2}\circ \tau_{1j_1}\ & \forall j_1,j_2,\\
\tau_{3j_3}\circ\tau_{2j_2} &= \tau_{2j_2}\circ \tau_{3j_3}\ & \forall j_1,j_3,\\
\tau_{1j_1}\circ\tau_{3j_3} &= \tau_{3j_3}\circ \tau_{1j_1}\ &
\forall j_1\not=j_3.
\end{align*}
Alors, si on fixe $j_1,j_2,j_3$ avec $j_1\not=j_3$, le groupe
engendré par  $\tau_{1j_1},\tau_{2j_2},\tau_{3j_3}$ est abélien
et fini (d'ordre $8$). Enfin, on définit
\[
\theta_{j_1j_2j_3} = \iota^*_{C_2(\hs)}\left( \sum_{a,b,c=0,1}
(-1)^{a+b+c}\;
(\tau_{1j_1}^*)^a\,(\tau_{2j_2}^*)^b\,(\tau_{3j_3}^*)^c\,
\iota^*_{C_2(\bbR^2)'}\theta\right),
\]
o\`u $\iota^*_{\bullet}$ est la restriction à $\bullet$. Il n'est
pas difficile de vérifier que les formes $\theta_{j_1j_2j_3}$,
$j_1\not=j_3$, satisfont aux conditions souhaitées.

\paragraph{Cas $j_1=j_3$ : } La construction précédente ne s'applique pas dans le cas
 $j_1=j_3$ car $\tau_1$ et $\tau_3$ ne commutent pas.
 Nous avons besoin d'une autre construction.

Soit $S$ la bande $\left\{\bx=(x,y)\in\bbR^2 : -\hpi\le y\le \hpi\right\}$.
Sur (l'intérieur de) $S$ on considère la métrique
\[
\dd s^2 = \frac{\dd x^2+\dd y^2}{\cos^2 y},
\]
qui approche la métrique de Poincaré sur chaque composante de
bord. Les géodésiques
\begin{figure}[h!]
\begin{center}
\includegraphics[width=5cm]{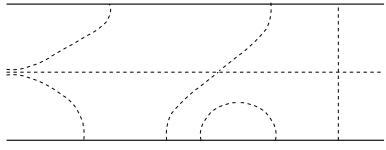} \caption{\footnotesize Géodésiques
dans la bande}
\end{center}
\end{figure}
sont, soit des lignes verticales, soit des courbes de la forme $ \sin y = A\ee^x+B\ee^{-x}$,
o\`u  $A$ et $B$ sont des paramètres.

Remarquons qu'il existe une et une seule géodésique passant par
deux points distincts de $S$.
On définit la fonction d'angle  $\phi(\bx_1,\bx_2)$
comme l'angle entre la géodésique verticale  passant par $\bx_1$
et la géodésique joignant $\bx_1$ à $\bx_2$. On définit
ensuite  $\vartheta$ comme $\dd\phi/(2\pi)$. Alors  $\vartheta$ est
une $1$-forme fermée sur $C_2(S)$ satisfaisant aux propriétés
suivantes :
\begin{enumerate}
\item $\vartheta$ s'annule quand $\bx_1$ s'approche du bord de $S$;
\item sur la composante de  bord, correspondant au rapprochement de $\bx_1$ et $\bx_2$
tout en restant dans l'intérieur de  $S$, $\vartheta$ s'identifie
à la forme volume normalisée et invariante sur $\mathbb{S}^1$;
\item sur la composante de  bord correspondant au rapprochement de $\bx_1$ et $\bx_2$
près du bord de  $S$, $\vartheta$ est la $1$-forme de Kontsevich.
\end{enumerate}
Soit $\sigma\colon S\to S$ l'involution correspondant à la
réflexion par rapport à $x=0$: c'est à dire,
$\sigma(x,y)=(-x,y)$. Soit $C_2(S)'$ le sous-ensemble de  $C_2(S)$
défini par  $\bx_1\not=\sigma(\bx_2)$. Remarquons que l'on a
$C_2(\hs)\subset C_2(S)'$. On écrit   $\sigma_{j}$ ($j=1,2$) pour
l'action de $\sigma$ sur la composante $j$ de $C_2(S)'$; c'est à
dire :
\begin{gather*}
\sigma_{1}\colon\begin{array}[t]{ccc}
C_2(S)'&\to &C_2(S)'\\
(\bx_1,\bx_2)&\mapsto &(\sigma(\bx_1),\bx_2)
\end{array}\\
\sigma_{2}\colon\begin{array}[t]{ccc}
C_2(S)'&\to &C_2(S)'\\
(\bx_1,\bx_2)&\mapsto &(\bx_1,\sigma(\bx_2))
\end{array}
\end{gather*}
Enfin, on définit
\[
\theta_{1j1} = \iota^*_{C_2(\hs)}\left( \sum_{a=0,1} (-1)^{a}\,
(\sigma_{j}^*)^a\, \iota^*_{C_2(S)'}\vartheta\right),
\]
et
\[
\theta_{2j2} = p^*\theta_{1,3-j,1}.
\]
Il n'est pas difficile de vérifier que les $1$-formes
$\theta_{j_1j_2j_3}$, $j_1=j_3$, satisfont aux conditions
désirées.

\subsubsection{Polarisations en position
d'intersection normale}

On utilise la forme à $8$ couleurs dans le cas de $3$
sous-algèbres $\k, \b_1, \b_2$ en position d'intersections normales.
c'est-à-dire il faut que l'on ait  (on peut intervertir les rôles de $\k, \b_1, \b_2$ ) :
 \[\k\cap (\b_1+ \b_2)=\k\cap \b_1+ \k\cap  \b_2.\]
   On peut alors trouver une base de $\g$ qui soit adaptée aux différentes
intersections. Dans ce cas, on peut adapter la forme à $8$ couleurs
en fonction des supplémentaires des  intersections possibles.

\begin{prop}\label{propindependance}

Soient $f\in \k^\perp$ et $\b_1$, $\b_2$ deux polarisations en $f$
telle que $\b_1, \b_2, \k$ soient en position d'intersections
normales. Alors le caractère cons\-truit est indépendant du choix de
la polarisation (une fois choisi un supplé\-men\-taire de
 $\k$ adapté à $\b_1$ et $\b_2$).
\end{prop}

\noindent \textit{\textbf{Preuve :}} Si les deux polarisations sont
en positions d'intersections normales on choisit un supplémentaire
de $\k$ adapté aux deux polarisations $\b_1$ et $\b_2$
simultanément. On dispose donc d'une algèbre de déformation
$H^0_{\epsilon, b_1,\b_2} (\k^\perp)$ (isomorphe
à l'algèbre $S(\p)^\k[[\epsilon]]$).\\

Considérons la forme à $8$ couleurs et les graphes de
Kontsevich colorés. Sur la variété de configurations des points dans
une demi-bande, plaçons aux coins
les fonctions $1$ et sur le coté borné la fonction $P\in
H^0_{\epsilon, b_1,\b_2} (\k^\perp)$ (\textit{cf.} (\ref{bandes}) dessin de gauche). Les contributions de tous les
graphes (pondérés par ces coefficients à $8$ couleurs)
représentent une fonction polynomiale sur $f+(\k+\b_1+\b_2)^\perp$.
Le résultat ne dépend pas de la position de la fonction  $P$.
Regardons les positions limites pour obtenir le résultat cherché.

\begin{equation}\label{bandes}
\begin{array}{lcr}
  \includegraphics[width=5cm]{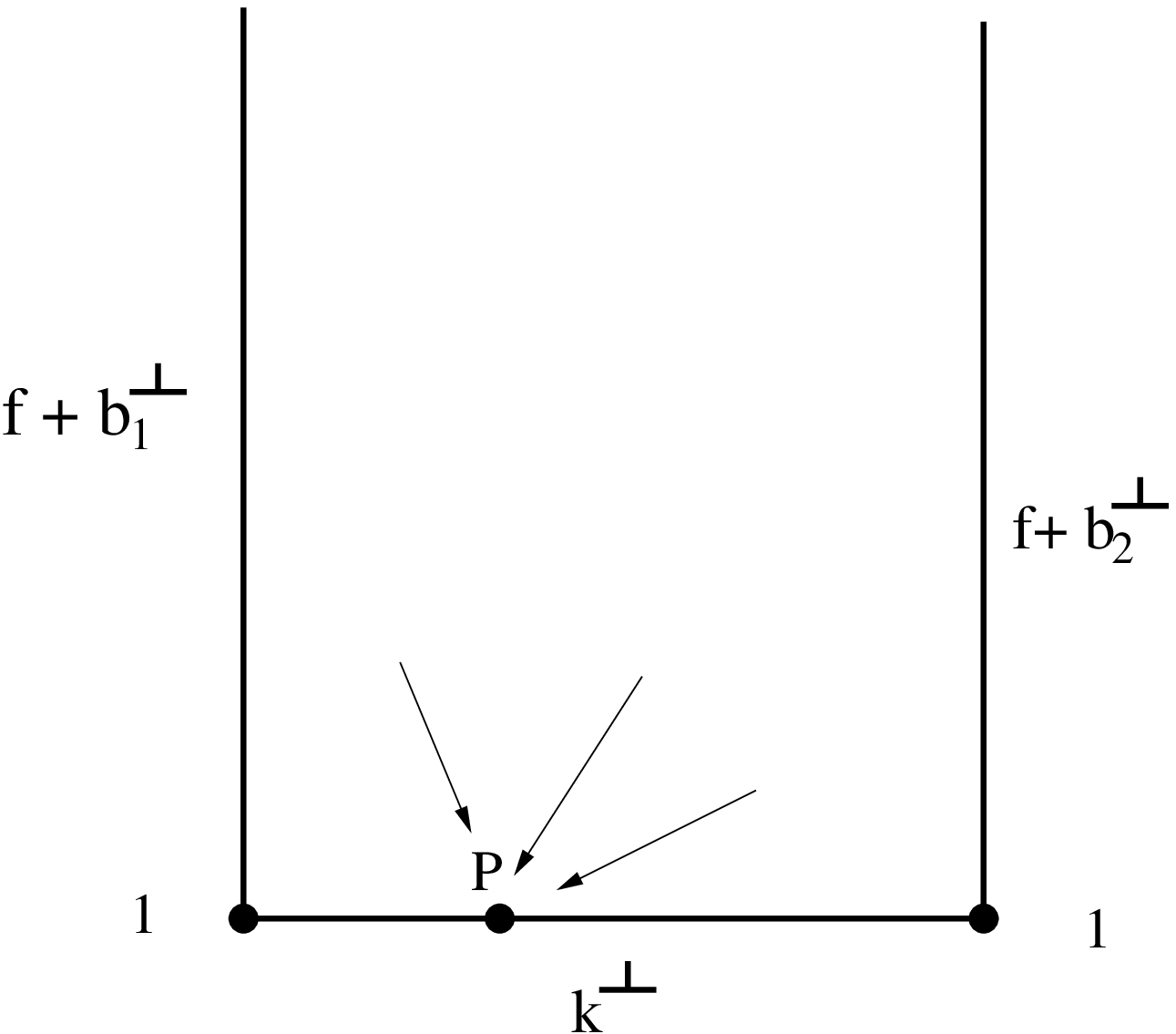} & \quad & \includegraphics[width=5cm]{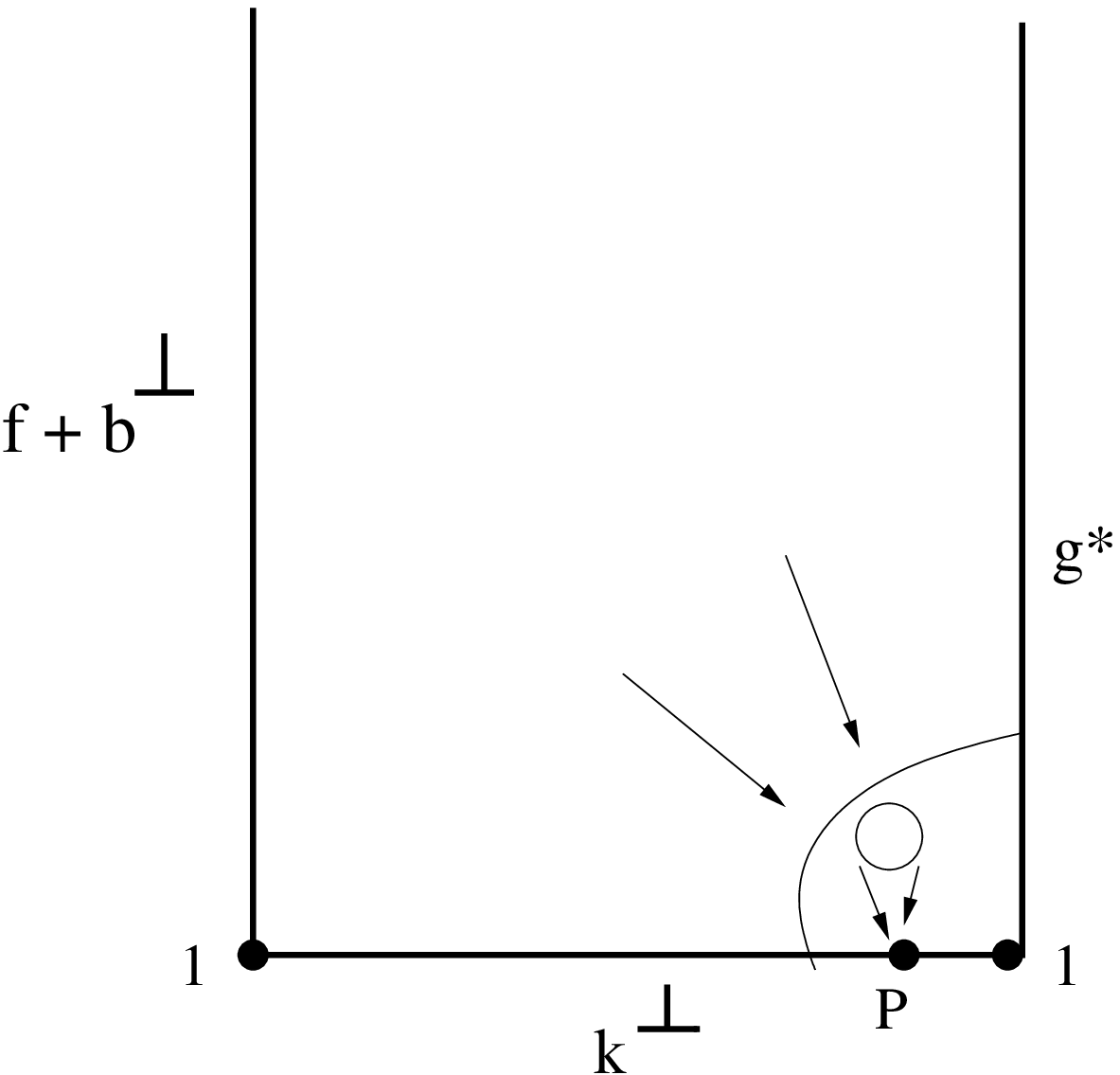} \\
  Ind\acute{e}pendance\; via \; la \;  forme \;
8 \; couleurs & \quad  & Calcul\; du \; caract\grave{e}re
\end{array}
\end{equation}

En s'approchant du coin de droite, on obtient, grâce aux
concentrations internes,
 un scalaire $1\underset{1}\star P$
(et donc il n'y a pas de contributions pour les concentrations
externes). On retrouve alors le caractère de la section précédente :

\[ H^0_{\epsilon, \b_1}(\k^\perp)
\rightarrow \R\]\[P\mapsto 1\underset{1}\star P.\] Lorsque la
fonction $P$ s'approche de l'autre coin on trouve le caractère
\[ H^0_{\epsilon, \b_2}(\k^\perp)
\rightarrow \R\]
\[P\mapsto  P\underset{2}\star 1.\]
On en déduit que les deux caractères construits sont identiques.
\fin

\paragraph{Remarque 11: }  Cette construction fonctionne si on choisit
comme sous-variétés co-isotropes $f_1+ \b_1^\perp$ et $f_2+
\b_2^\perp$, où $f_1, f_2$ sont  dans $\k^\perp$ et vérifient  la condition
d'intersection :

\[\big( f_1+ (\k+\b_1)^\perp\big) \cap \big(f_2+ (\k+\b_2)^\perp \big)\neq \emptyset.\]
Si l'intersection ci-dessus est non vide, alors $f_1$ et $f_2 $ sont
dans la même $K$-orbite pour peu que la condition de Pukanszky
soit vérifiée. En effet si la condition de Pukanszky (\textit{cf.}
 \S~\ref{sectioncasdespolarisations}) est vérifiée   on a
\[(K\cap B_1)\cdot f_1= f_1+ (\k+\b_1)^\perp \quad  \mathrm{ et }
\quad (K\cap B_1)\cdot f_2=f_2+ (\k+\b_2)^\perp. \]Alors $f_1$ et
$f_2$ sont conjuguées par $K$ dans~$\k^\perp$.

\subsection{Applications au cas des polarisations $\sigma$-stables}\label{sectionsigmastable}
On applique la construction ci-dessus, dans le cas des paires
symétriques qui admettent de manière générique des polarisations
$\sigma$-stables pour les formes linéaires dans $\k^\perp$.

\subsubsection{Exemples classiques}

 Les paires symétriques pour lesquelles on peut construire des
polarisations $\sigma$-stables sont proches heuristiquement des
algèbres de Lie. Voici quelques exemples (\textit{cf.} \cite{To1, To4}).
\begin{enumerate}
 \item Les petites paires symétriques comme dans \S~\ref{sectioniwasawa},
 \item  le cas
de paires symétriques nilpotentes, \item  le cas des paires
symétriques de Alekseev-Meinrenken, \item  le cas  des algèbres de
Lie considérées  comme des espaces symétriques,
\item les paires symétriques de  Takiff $\g \oplus \g\otimes T$.
\end{enumerate}

\subsubsection{Indépendance du caractère}
Lorsque l'on sait construire des polarisations $\sigma$-stables, les
supplémen\-taires s'imposent  pour $\b/(\b\cap \k)$ ainsi que
$\g/(\k+ \b)$ ;  on choisit $\b\cap \p $ et $\p/\p\cap \b$ (un
supplémentaire de $\p\cap \b$ dans $\p$). L'algèbre de réduction
pour $\k^\perp$ est alors tout simplement $S(\p)^\k$ munie du
produit de $\underset{CF}\star$, car tous les espaces
considérés sont $\sigma$-stables.\\

Différents choix de polarisations $\sigma$-stables vont être en
position d'intersection normale. Ainsi d'après la
Proposition~\ref{propindependance}, le caractère construit est
indépendant de la polarisation $\sigma$-stable; notons le
\[P\mapsto \gamma(P)(f).\]

Par ailleurs cette construction est équivariante par rapport à
l'action linéaire de $K$, donc le caractère construit ne dépend que
de la $K$-orbite de~$f$. On en déduit la proposition suivante (\textit{cf.}
\cite{To1, To4})

\begin{prop}\label{propcaractere} S'il existe génériquement des polarisations
$\sigma$-stables, alors la construction ci-dessus définit un
caractère de $\big(S(p)^\k, \underset{CF}\star\big)$. Ce caractère
ne dépend pas du choix de la polarisation $\sigma$-stable. Il est
constant sur la $K$-orbite et  polynomial en $f$.
\end{prop}

\noindent  \textit{\bf Preuve :} Le caractère ne dépend pas du choix
de la polarisation $\sigma$-stable, car on est automatiquement en
situation d'intersection normale et l'algèbre de réduction ne dépend
pas des
choix des  supplémentaires, car ils sont pris dans $\p$.\\

La dépendance polynomiale en $\b$ est claire, tant que l'on reste
sur des ouverts où $f$ admet des polarisations génériques. On
applique la méthode de Duflo-Conze \cite{du77} pour conclure. L'ensemble des $(f, \b)$ avec $f\in \k^\perp$ et $\b$ sous-algèbre
$\sigma$-stable subordonnée à $f$ et de dimension générique, tels
que $1\star P $
 soit scalaire dans $f+ (\k+\b)^\perp$ forme un
ensemble rationnel. La fibre étant projective, on conclut alors que le caractère rationnel est en fait polynomial
et $K$-invariant.

\fin

\subsubsection{Isomorphisme de Rouvière}

L'homomorphisme $\gamma$ est un isomorphisme de $S(\p)^\k $
muni du produit~$\underset{CF}\star$ sur $S(\p)^\k$ muni du produit
standard, car le terme de plus haut degré est clairement l'identité.\\

En fait on a le résultat plus fort suivant (démontré dans \cite{To4}
par la méthode des orbites, voir aussi~\cite{AM, To3}).

\begin{theo}\label{theosigmastable} Dans le cas paires symétriques, avec polarisations
géné\-riques $\sigma$-stables, l'homomorphisme $\gamma$ vaut
l'identité. En conséquence l'application\footnote{Ici $\beta$
désigne la symétrisation.}
\[P\mapsto \beta\left (\partial(J^{\frac 12}) P\right)\] est un
isomorphisme d'algèbres de $S(\p)^\k$ sur $\left(U(\g)/U(\g)\cdot
\k\right)^\k$.
\end{theo}
\noindent \textit{\bf Preuve :} Soit $f\in \k^\perp$,  $\b$
polarisation $\sigma$-stable en $f$. On considère dans le diagramme
de la bi-quantification la situation du triplet : $f+\b^\perp$,
 $\k^\perp$ et $\g^*$. On place les fonctions
$1$ aux deux coins et la fonction $P\in S(\p)^\k$ sur l'axe
horizontal (\textit{cf.} (\ref{bandes}) dessin de droite).\\

Lorsque $P$ se rapproche du coin $f+\b^\perp$ on retrouve le caractère recherché.\\

Lorsque $P$ se rapproche du coin correspondant à $\g^*$, on
retrouve dans les strates qui se concentrent (sur l'axe horizontal)
la situation des paires symétriques à savoir $P\star 1=B(P)$, où
$B$ désigne l'opérateur associé aux roues
pures attachées sur l'axe  $\k^\perp$ (\textit{cf.} \S\ref{sectionrouesHV}).\\

Considérons les contributions des graphes extérieurs et montrons qu'elles sont triviales. Les arêtes
qui arrivent sur le coin sont colorées par les couleurs $+++$ ou
$-++$. En fait pour la fonction d'angle à $8$ couleurs, on a
$\phi_{-++}(p, q)=0$ si $q$ est dans le coin. Par conséquent les
arêtes qui arrivent sur le coin sont d'une seule couleur : $+++$.
Cette couleur correspond à une dérivée dans la direction
$(\p/\p\cap \b)^*=(\k+\b)^\perp$. Or la fonction que l'on dérive est
$B(P)$, qui est $\k$-invariante et donc constante sur
$f+(\k+\b)^\perp$. De plus l'opérateur est évalué en $\xi \in f+
(\k+\b)^\perp$, par conséquent ses dérivées sont nulles. Le
caractère construit dans la Proposition~\ref{propcaractere} vaut donc $B(P)(f)$.\\

 \noindent La série $B$ est universelle et s'écrit sous la forme
:
\[\exp\left(\sum w_n \tr_\p(\ad X)^{2n}\right),\] où $w_n$ sont des
constantes universelles. On en déduit que l'application de
\[P\mapsto B(P)\]de $(S(\p)^\k,
\underset{CF}\star=\underset{Rou}\sharp)$ sur $(S(\p)^\k, \cdot)$
est un isomorphisme d'algèbres pour les paires symétriques  qui
admettent génériquement des polarisations $\sigma$-stables. Pour
terminer la preuve du théorème il suffit de montrer que l'on a
$B=1$.

\fin

\begin{prop}\label{B=1}La série universelle $B(X)$ qui intervient
pour les paires symé\-triques vaut~$1$.
\end{prop}

 \noindent \textit{\bf Preuve :} D'après \cite{To4} lorsque
génériquement on sait construire des polarisations $\sigma$-stables,
le produit
de Rouvière (qui vaut aussi le  star-produit $\underset{CF}\star$) et le produit standard coïncident sur $S(\p)^\k$.\\

 Si on avait
$B\neq 1$ alors on disposerait d'un isomorphisme d'algèbres non trivial, pour
$S(\p)^\k$ muni du produit standard (pour ce genre de  paires
symétriques). Il existerait $n>0$ tel que
$\tr_\p(\ad X)^{2n}$ serait une dérivation de $S(\p)^\k$, ce qui
n'est pas vrai sur des exemples \footnote{Par exemple pour $sl(2)$
considérée comme paire
symétrique}.\\

La série universelle $B$ vaut~$1$, ce qui justifie la formule pour
l'écriture des opérateurs différentiels invariants en coordonnées
exponentielles du \S~\ref{sectionOpdiffexp}.

\fin

\begin{paragraph}{Remarque 12 : }
 Dans le cas des paires symétriques de Takiff $\g \oplus
\g\otimes T$, avec $\g$ une paire symétrique  et $T^2=0$, on a
\[E_{\g \oplus \g\otimes T}(X+X'\otimes T,Y+Y'\otimes T)=E^2_\g(X,
Y).\] On sait toutefois que dans ce contexte, l'homomorphisme
d'Harish-Chandra coïncide avec l'isomorphisme de Rouvière car on
peut construire des polarisations $\sigma$-stables (\textit{cf.}
\S~\ref{sectionsigmastable}). Ceci montre que la formule de
Rouvière peut réaliser un isomorphisme d'algèbres  sans que la
fonction $E_{\g \oplus \g\otimes T}$ soit égale à $1$. La question
pertinente est donc de savoir si cette fonction est homotope à $1$
modulo les
 champs $\k\oplus \k_T$-adjoints.
 \end{paragraph}

\end{document}